\documentclass[a4paper, 12pt]{article}
\usepackage[T2A]{fontenc}
\usepackage[utf8x]{inputenc}
\usepackage[english,russian]{babel}
\usepackage{amsmath}
\usepackage{amsfonts}
\usepackage{amsbib}
\DeclareMathAlphabet{\mathpzc}{OT1}{pzc}{m}{it}

\newcommand{\sgn}{\operatorname{sign}}\newcommand{\loc}{\operatorname{loc}}

\newcommand{\e}{\mathfrak e}
\renewcommand{\l}{\mathfrak l}

\newcommand{\mn}{\operatorname{\mathfrak m}}

\newcommand{\s}{\mathpzc s}

\newcommand{\card}{\operatorname{card}}

\newcommand{\supp}{\operatorname{supp}}
\newcommand{\supvrai}{\operatornamewithlimits{sup\,vrai}}
\newcommand{\N}{\mathbb N}
\newcommand{\Z}{\mathbb Z}
\newcommand{\R}{\mathbb R}
\newcommand{\Nu}{\mathcal N}
\newcommand{\Rho}{\mathrm P}
\newcommand{\D}{\mathcal D}
\newcommand{\J}{\mathcal J}

\newcommand{\mes}{\operatorname{mes}}
\allowdisplaybreaks

\begin{document}

\author{ С. Н. Кудрявцев }
\title{Продолжение за пределы куба функций из пространств 
Никольского -- Бесова смешанной гладкости}
\date{}
\maketitle
\begin{abstract}
В статье рассмотрены пространства Никольского и Бесова 
с нормами, в определении которых вместо смешанных модулей непрерывности 
известных порядков определённых смешанных производных функций 
используются "$L_p$-усреднённые" смешанные модули непрерывности функций 
соответствующих порядков. Для таких пространств функций, заданных на кубе 
$ I^d, $ построены непрерывные линейные отображения их в обычные 
пространства Никольского и Бесова смешанной гладкости в $ \R^d, $ 
являющиеся операторами продолжения функций, что влечёт совпадение тех и
других пространств на кубе $ I^d. $
\end{abstract}

Ключевые слова: пространства Никольского -- Бесова смешанной гладкости,
продолжение функций, эквивалентные нормы

\bigskip

\centerline{Введение}
\bigskip

В работе рассмотрены некоторые задачи теории функциональных пространств 
для пространств Никольского и Бесова функций смешанной гладкости, заданных
на кубе $ I^d. $ Остановимся подробнее на содержании работы.

При $ d \in \N, \alpha \in \R_+^d, 1 \le p,\theta < \infty $ для области $ D
\subset \R^d $ вводятся в рассмотрение
пространства $ (S_{p,\theta}^\alpha B)^\prime(D) ((S_p^\alpha H)^\prime(D)) $
с нормами
\begin{multline*}
\| f \|_{(S_{p,\theta}^\alpha B)^\prime(D)} =\\
\max\biggl(\| f \|_{L_p(D)},
\max_{J \subset \{1,\ldots,d\}: J \ne \emptyset} 
\left(\int_0^\infty \ldots \int_0^\infty 
(\prod_{j \in J} t_j^{-1 -\theta \alpha_j}) 
(\Omega^{\prime l \chi_J}(f, t^J)_{L_p(D)})^{\theta} \,dt^J
\right)^{1/\theta}\biggr),
\end{multline*}
$$
\| f \|_{(S_p^\alpha H)^\prime(D)} = \max(\| f \|_{L_p(D)}, 
\max_{J \subset \{1, \ldots, d\}: J \ne \emptyset}
\sup_{t \in \R^d: t_j > 0, j \in J} (\prod_{j \in J} t_j^{-\alpha_j}) 
\Omega^{\prime l \chi_J}(f, t^J)_{L_p(D)}),
$$
где
\begin{multline*}
\Omega^{\prime l \chi_J}(f, t^J)_{L_p(D)} =
\biggl((\prod_{j \in J} (2 t_j)^{-1}) \int_{ \{\xi^J: |\xi_j| \le t_j, j \in J\}}
\| (\prod_{j \in J} \Delta_{\xi_j e_j}^{l_j}) f\|_{L_p(D_{\xi}^{l \chi_J})}^p d\xi^J\biggr)^{1 /p}, \\
l \in \N^d -- \text{ вектор с компонентами } \\
l_j = \min \{m \in \N: \alpha_j < m \}, j =1,\ldots,d, \\ 
\text{ а для } J = \{j_1, \ldots, j_k\}: 1 \le j_1 < \ldots < j_k \le d, \text{ и }
x \in \R^d \text{ вектор } x^J = (x_{j_1},\ldots,x_{j_k}).
\end{multline*}

Построены непрерывные линейные отображения пространств 
$ (S_{p,\theta}^\alpha B)^\prime(I^d) \ ((S_p^\alpha H)^\prime(I^d)) $ в 
пространства $ L_p(\R^d), $ являющиеся операторами продолжения
функций, значения которых обладают дифференциально-разностными свойствами, 
влекущими совпадение этих пространств с соответствующими обычными пространствами 
Бесова $ S_{p,\theta}^\alpha B(I^d) $
(Никольского $ S_p^\alpha H(I^d)) $ (см. п. 2.3.).
Публикации, в которых рассматривается такая задача, автору не известны. Из
близких работ по этой тематике приведём [1] -- [10] (см. также имеющуюся
там литературу), в которых изучается вопрос о продолжении за пределы
области определения с сохранением класса гладких функций из изотропных и 
неизотропных пространств, а также функций из пространств смешанной 
гладкости. Отметим, что средства построения операторов 
продолжения функций и схемы доказательства дифференциально-
разностных свойств значений таких операторов, применяемые ниже, отличаются от тех, что использовались в 
упомянутых работах. Добавим ещ\:е, что нормы, определяющие пространства функций 
смешанной гладкости, рассматриваемые в [10], отличаются от 
$ \| \cdot \|_{(S_{p,\theta}^\alpha B)^\prime(D)}, 
\| \cdot \|_{(S_p^\alpha H)^\prime(D)}. $

Статья состоит из введения и двух параграфов, в первом из которых
рассматриваются некоторые вспомогательные средства для решения 
объявленной задачи, а во втором --средства построения и конструкция 
подходящих операторов продолжения, а также доказательство их свойств.
\bigskip

\centerline{\S 1. Предварительные сведения и вспомогательные утверждения}
\bigskip

1.1. В этом пункте вводятся обозначения, относящиеся к 
функциональным пространствам, рассматриваемым в
настоящей работе, а также приводятся некоторые факты, необходимые
в дальнейшем.

Для $ d \in \N $ через $ \Z_+^d $ обозначим множество
$$
\Z_+^d = \{\lambda = (\lambda_1,  \ldots, \lambda_d) \in \Z^d:
\lambda_j \ge0, j=1, \ldots, d\}.
$$
 Обозначим также при  $ d \in \N $ для $ l \in \Z_+^d $ через
$ \Z_+^d(l) $ множество
 $$
 \Z_+^d(l) =\{ \lambda  \in \Z_+^d: \lambda_j \le l_j, j=1, \ldots,
 d\}.
 $$

Для $  d \in  \N, l \in \Z_+^d $ через $ \mathcal P^{d,l} $ будем
обозначать пространство вещественных  полиномов, состоящее из всех
функций $ f: \R^d \mapsto \R $ вида
$$
f(x) = \sum_{\lambda  \in \Z_+^d(l)}a_{\lambda}\cdot x^{\lambda},
x\in \R^d,
$$
где  $   a_{\lambda}   \in   \R,  x^{\lambda}  =x_1^{\lambda_1}
\ldots x_d^{\lambda_d}, \lambda \in \Z_+^d(l). $

При $ d \in \N, l \in \Z_+^d  $   для области $   D \subset \R^d $
через $\mathcal   P^{d,l}(D) $  обозначим пространство функций $
f, $ определённых в $  D, $   для каждой из  которых существует
полином $ g \in \mathcal P^{d,l} $
 такой, что  сужение $ g \mid_D = f.$

Для измеримого по Лебегу множества $ D \subset \R^d $ при $ 1 \le p \le \infty $ 
через $ L_p(D), $ как обычно, обозначается
пространство  всех  вещественных измеримых на $ D $ функций $f,$
для которых определена норма
$$
\|f\|_{L_p(D)} = \begin{cases} (\int_D |f(x)|^p dx)^{1/p}, 1 \le p < \infty;\\
\supvrai_{x \in D}|f(x)|, p = \infty. \end{cases}
$$

Обозначим через $ L^{\loc}(\R^d) $ линейное пространство всех локально
суммируемых в $ \R^d $ функций, принимающих вещественные значения.

Отметим здесь неоднократно используемое в дальнейшем неравенство
\begin{equation*} \tag{1.1.1}
| \sum_{j =1}^n x_j|^a \le \sum_{j =1}^n |x_j|^a, x_j \in \R, j =1,\ldots,n, 
n \in \N, 0 \le a \le 1.
\end{equation*}

Для $ x,y \in \R^d $ положим $ xy = x \cdot y = (x_1 y_1, \ldots, x_d y_d), $ 
а для $ x \in \R^d $ и $ A \subset \R^d $ определим
$$
x A = x \cdot A = \{xy: y \in A\}.
$$

Для $ x \in \R^d: x_j \ne 0, $ при $ j=1,\ldots,d,$ положим 
$ x^{-1} = (x_1^{-1},\ldots,x_d^{-1}). $

При $ d \in \N $ для $ x,y \in \R^d $ будем писать $ x \le y (x < y), $
если для каждого $ j=1,\ldots,d $ выполняется неравенство $ x_j
\le y_j (x_j < y_j). $

Обозначим через $ \R_+^d $ множество $ x \in \R^d, $ для которых 
$ x_j >0 $ при $ j=1,\ldots,d,$ и для $ a \in \R_+^d, x \in \R^d $
положим $ a^x = a_1^{x_1} \ldots a_d^{x_d}.$

При $ d \in \N $ определим множества
$$
I^d = \{x \in \R^d: 0 < x_j < 1,j=1,\ldots,d\},
$$
$$
\overline I^d = \{x \in \R^d: 0 \le x_j \le 1,j=1,\ldots,d\},
$$
$$
B^d = \{x \in \R^d: -1 \le x_j \le 1,j=1,\ldots,d\}.
$$

Через $ \e $ будем обозначать вектор в $ \R^d, $ задаваемый
равенством $ \e = (1,\ldots,1). $

При $ d \in \N $ для $ \lambda \in \Z_+^d $ через $ \D^\lambda $
будем обозначать оператор дифференцирования $ \D^\lambda =
\frac{\D^{|\lambda|}} {\D x_1^{\lambda_1} \ldots \D x_d^{\lambda_d}}, $ 
где $ |\lambda| = \sum_{j=1}^d \lambda_j. $

Теперь приведём некоторые факты, относящиеся к полиномам, которыми
мы будем пользоваться ниже.

В [11] содержится такое утверждение.

    Лемма 1.1.1  

Пусть $ d \in \N, l \in \Z_+^d, \lambda \in \Z_+^d, 1 \le p,q \le \infty, 
\rho, \sigma \in \R_+^d. $ Тогда существует константа
$ c_1(d,l,\lambda,\rho, \sigma) >0 $ такая, что для любых
измеримых по Лебегу множеств $ D,Q \subset \R^d, $  для которых
можно найти $ \delta \in \R_+^d $ и $ x^0 \in \R^d $ такие, что 
$ D \subset (x^0 +\rho \delta B^d) $ и $ (x^0 +\sigma \delta I^d) \subset Q, $ 
для любого полинома $ f \in \mathcal P^{d,l} $
выполняется неравенство
    \begin{equation*} \tag{1.1.2}
\| \D^\lambda f\|_{L_q(D)} \le c_1 \delta^{-\lambda -p^{-1} \e +q^{-1} \e} 
\|f\|_{L_p(q)}.
   \end{equation*}

Далее, напомним, что для открытого множеста $ D \subset \R^d $ и вектора 
$ h \in \R^d $ через $ D_h $ обозначается множество
$$
D_h = \{x \in D: x +th \in D \ \forall t \in \overline I\}.
$$

Для функции $ f, $ заданной на открытом множестве $ D \subset \R^d, $ и
вектора $ h \in \R^d $ определим на $ D_h $ её разность $ \Delta_h f $ 
с шагом $ h, $ полагая
$$
(\Delta_h f)(x) = f(x+h) -f(x), x \in D_h,
$$
а для $ l \in \N: l \ge 2, $ на $ D_{lh} $ определим $l$-ую
разность $ \Delta_h^l f $ функции $ f $ с шагом $ h $ равенством
$$
(\Delta_h^l f)(x) = (\Delta_h (\Delta_h^{l-1} f))(x), x \in
D_{lh},
$$
положим также $ \Delta_h^0 f = f. $

Как известно, справедливо равенство
$$
(\Delta_h^l f)(\cdot) = \sum_{k=0}^l C_l^k (-1)^{l-k} f(\cdot +kh), 
C_l^k =\frac{l!} {k! (l-k)!}.
$$

При $ d \in \N $ для $ j=1,\ldots,d$ через $ e_j $ будем
обозначать вектор $ e_j = (0,\ldots,0,1_j,0,\ldots,0).$

Как показано в [12], справедлива следующая лемма.

   Лемма 1.1.2

   Пусть $ d \in \N, l \in \Z_+^d. $ Тогда
для любого $ \delta \in \R_+^d $ и $ x^0 \in \R^d $ для $ Q = x^0 +\delta I^d $ 
существует единственный линейный оператор 
$ P_{\delta, x^0}^{d,l}: L_1(Q) \mapsto \mathcal P^{d,l}, $
обладающий следующими свойствами:

1) для $ f \in \mathcal P^{d,l} $ имеет место равенство
\begin{equation*} \tag{1.1.3}
P_{\delta, x^0}^{d,l}(f \mid_Q) = f,
  \end{equation*}

2)
\begin{equation*} 
\ker P_{\delta,x^0}^{d,l} = \biggl\{\,f \in L_1(Q): 
\int \limits_{Q} f(x) g(x) \,dx =0\ \forall g \in \mathcal P^{d,l}\,\biggr\},
\end{equation*}

причём, существуют константы $ c_2(d,l) >0 $ и $ c_3(d,l) >0 $
такие, что

   3) при $ 1 \le p \le \infty $ для $ f \in L_p(Q) $ справедливо неравенство
  \begin{equation*} \tag{1.1.4}
\|P_{\delta, x^0}^{d,l}f \|_{L_p(Q)} \le c_2 \|f\|_{L_p(Q)},
  \end{equation*}

4) при $ 1 \le p < \infty $ для $ f \in L_p(Q) $ выполняется неравенство
   \begin{equation*} \tag{1.1.5}
  \|f -P_{\delta, x^0}^{d,l}f \|_{L_p(Q)} \le c_3 \sum_{j=1}^d
\delta_j^{-1/p} \biggl(\int_{\delta_j B^1} \int_{Q_{(l_j +1) \xi e_j}}
|\Delta_{\xi e_j}^{l_j +1} f(x)|^p dx d\xi\biggr)^{1/p}.
\end{equation*}

Теперь определим пространства функций, изучаемые в настоящей работе (ср. с  
[13], [14]). Но прежде введём некоторые обозначения.

При $ d \in \N $ для $ x \in \R^d $ через $\s(x) $ обозначим
множество $ \s(x) = \{j =1,\ldots,d: x_j \ne 0\}, $ а для
множества $ J \subset \{1,\ldots,d\} $ через $ \chi_J $ обозначим
вектор из $ \R^d $ с компонентами
$$
(\chi_J)_j = \begin{cases} 1, &\text{ для } j \in J; \\ 
0, &\text{ для } j \in (\{1,\ldots,d\} \setminus J). \end{cases}
$$

При $ d \in \N $ для $ x \in \R^d $ и $ J = \{j_1,\ldots,j_k\}
\subset \N: 1 \le j_1 < j_2 < \ldots < j_k \le d, $ через $ x^J $
обозначим вектор $ x^J = (x_{j_1},\ldots,x_{j_k}) \in \R^k, $ а
для множества $ A \subset \R^d $ положим $ A^J = \{x^J: x \in A\}. $

Для открытого множества $ D \subset \R^d $ и векторов $ h \in \R^d $ и $ l \in
\Z_+^d $ через $ D_h^l $ обозначим множество
\begin{multline*}
D_h^l = (\ldots (D_{l_d h_d e_d})_{l_{d-1} h_{d-1} e_{d-1}}
\ldots)_{l_1 h_1 e_1} = \{ x \in D: x +tlh \in D \ \forall t \in
\overline I^d\} = \\ \{ x \in D: (x +\sum_{j \in \s(l)} t_j l_j h_j e_j) \in D \  
\forall t^{\s(l)} \in (\overline I^d)^{\s(l)} \}.
\end{multline*}

Пусть $ d \in \N, D $ -- открытое множество в $ \R^d $ и $ 1 \le p \le \infty. $ 
Тогда для $ f \in L_p(D), h \in \R^d $ и $ l \in \Z_+^d $ определим в $ D_h^l $ 
смешанную разность функции $ f $ порядка $ l, $ соответствующую вектору $ h, $ 
равенством
\begin{multline*}
(\Delta_h^l f)(x) = \biggl(\biggl(\prod_{j=1}^d \Delta_{h_j e_j}^{l_j}\biggr) f\biggr)(x)
= \biggl(\biggl(\prod_{j \in \s(l)} \Delta_{h_j e_j}^{l_j}\biggr) f\biggr)(x) = \\ 
\sum_{k \in \Z_+^d(l)} (-\e)^{l-k} C_l^k f(x+kh), x \in D_h^l,
\end{multline*}
где $ C_l^k = \prod_{j=1}^d C_{l_j}^{k_j}. $

Имея в виду, что для $ f \in L_p(D), l \in \Z_+^d $ и векторов 
$ h,h^\prime \in \R^d: h^{\s(l)} = (h^\prime)^{\s(l)}, $ соблюдается
соотношение
$$
\| \Delta_h^l f\|_{L_p(D_h^l)} = \| \Delta_{h^\prime}^l
f\|_{L_p(D_{h^\prime}^l)}, 1 \le p \le \infty,
$$
определим при $ 1 \le p \le \infty $ для функции $ f \in L_p(D) $ смешанный 
модуль непрерывности в $ L_p(D) $ порядка $ l \in \Z_+^d $ равенством
$$
\Omega^l (f,t^{\s(l)})_{L_p(D)} = \supvrai_{ \{ h \in \R^d:
h^{\s(l)} \in t^{\s(l)} (B^d)^{\s(l)} \}} \| \Delta_h^l f\|_{L_p(D_h^l)}, 
t^{\s(l)} \in (\R_+^d)^{\s(l)}.
$$
Кроме того, при тех же условиях введём в рассмотрение для функции $ f $
"усреднённый" смешанный модуль непрерывности в $ L_p(D) $ порядка $ l, $
полагая   
\begin{multline*}
\Omega^{\prime l} (f, t^{\s(l)})_{L_p(D)} = \begin{cases}
\biggl((2 t^{\s(l)})^{-\e^{\s(l)}}
\int_{ t^{\s(l)} (B^d)^{\s(l)}} \| \Delta_\xi^l f\|_{L_p(D_\xi^l)}^p 
d \xi^{\s(l)}\biggr)^{1 /p} = \\
\biggl((2 t^{\s(l)})^{-\e^{\s(l)}}
\int_{ (t B^d)^{\s(l)}} \int_{D_\xi^{l \chi_{\s(l)}}} 
| \Delta_\xi^{l \chi_{\s(l)}} f(x)|^p dx 
d \xi^{\s(l)}\biggr)^{1 /p}, p \ne \infty; \\
\Omega^l (f,t^{\s(l)})_{L_p(D)}, p = \infty, 
\end{cases}\\
t^{\s(l)} \in (\R_+^d)^{\s(l)}.
\end{multline*}

Из приведенных определений видно, что
\begin{multline*} \tag{1.1.6}
\Omega^{\prime l} (f, t^{\s(l)})_{L_p(D)} \le
\Omega^l (f, t^{\s(l)})_{L_p(D)}, t^{\s(l)} \in (\R_+^d)^{\s(l)}, \\
f \in L_p(D), 1 \le p \le \infty, l \in \Z_+^d, D \text{-- произвольное 
открытое множество в } \R^d.
\end{multline*}

Пусть теперь $ d \in \N, \alpha \in \R_+^d, 1 \le p \le \infty, D $ --
область в $ \R^d $ и вектор $ \l \in \Z_+^d $ такой, что $ \l < \alpha. $ 
Тогда зададим вектор $ l = l(\alpha) \in \N^d, $ полагая 
$ l_j = (l(\alpha))_j = \min \{m \in \N: \alpha_j < m \}, j =1,\ldots,d, $ 
и через $ (S_p^\alpha H)^{\l}(D) $ обозначим пространство всех функций 
$ f \in L_p(D), $ обладающих тем свойством, что для любого непустого 
множества $ J \subset \{1,\ldots,d\} $ обобщённая производная 
$ \D^{\l \chi_J} f \in L_p(D) $ и выполняется условие
\begin{multline*}
\sup_{t^J \in (\R_+^d)^J} (t^J)^{-(\alpha^J -\l^J)}
\Omega^{(l -\l) \chi_J}(\D^{\l \chi_J} f,t^J)_{L_p(D)} =\\ 
\sup_{t^J \in (\R_+^d)^J}
\biggl(\prod_{j \in J} t_j^{-(\alpha_j -\l_j)}\biggr) \Omega^{(l -\l) \chi_J}(\D^{\l \chi_J} f,
t^{\s((l -\l) \chi_J)})_{L_p(D)} < \infty.
\end{multline*}

В пространстве $ (S_p^\alpha H)^{\l}(D) $ вводится норма 
\begin{multline*}
\| f \|_{(S_p^\alpha H)^{\l}(D)} =\\ 
\max\biggl(\| f \|_{L_p(D)}, \max_{J \subset \{1, \ldots, d\}: 
J \ne \emptyset} \sup_{t^J \in (\R_+^d)^J} (t^J)^{-(\alpha^J -\l^J)}
\Omega^{(l -\l) \chi_J}(\D^{\l \chi_J} f,t^J)_{L_p(D)}\biggr), \\
f \in (S_p^\alpha H)^{\l}(D).
\end{multline*} 
При тех же условиях на $ \alpha, p, D $ обозначим через 
$ (S_p^\alpha H)^\prime(D) $ пространство всех функций 
$ f \in L_p(D), $ обладающих тем свойством, что для любого непустого 
множества $ J \subset \{1,\ldots,d\} $ соблюдается условие
$$
\sup_{t^J \in (\R_+^d)^J} (t^J)^{-\alpha^J}
\Omega^{\prime l \chi_J}(f, t^J)_{L_p(D)} = \sup_{t^J \in (\R_+^d)^J}
\biggl(\prod_{j \in J} t_j^{-\alpha_j}\biggr) \Omega^{\prime l \chi_J}(f,
t^{\s(l \chi_J)})_{L_p(D)} < \infty.
$$

В пространстве $ (S_p^\alpha H)^\prime(D) $ задаётся норма 
\begin{multline*}
\| f \|_{(S_p^\alpha H)^\prime(D)} =\\ 
\max\biggl(\| f \|_{L_p(D)}, \max_{J \subset \{1, \ldots, d\}: 
J \ne \emptyset} \sup_{t^J \in (\R_+^d)^J} (t^J)^{-\alpha^J}
\Omega^{\prime l \chi_J}(f, t^J)_{L_p(D)}\biggr), \\
f \in (S_p^\alpha H)^\prime(D).
\end{multline*} 

Для области $ D \subset \R^d $ при $ \alpha \in \R_+^d,\ 1 \le p \le \infty,
\theta \in \R: 1 \le \theta < \infty, \l \in \Z_+^d: \l < \alpha, $ 
полагая, как и выше, $ l = l(\alpha) \in \N^d, $ 
через $ (S_{p,\theta}^\alpha B)^{\l}(D) $ обозначим пространство всех функций 
$ f \in L_p(D), $ у которых для
любого непустого множества $ J \subset \{1,\ldots,d\} $ обобщённая производная 
$ \D^{\l \chi_J} f \in L_p(D) $ и соблюдается условие
\begin{multline*}
\int_{(\R_+^d)^J} (t^J)^{-\e^J -\theta (\alpha^J -\l^J)} 
(\Omega^{(l -\l) \chi_J}(\D^{\l \chi_J} f, t^J)_{L_p(D)})^\theta dt^J = \\
\int_{(\R_+^d)^J} \biggl(\prod_{j \in J} t_j^{-1 -\theta (\alpha_j -\l_j)}\biggr) 
(\Omega^{(l -\l) \chi_J}(\D^{\l \chi_J} f, t^{\s((l -\l) \chi_J)})_{L_p(D)})^\theta 
\prod_{j \in J} dt_j < \infty. 
\end{multline*}

В пространстве $ (S_{p,\theta}^\alpha B)^{\l}(D) $ фиксируется норма
\begin{multline*}
\| f \|_{(S_{p,\theta}^\alpha B)^{\l}(D)} = \\
\max\biggl(\| f\|_{L_p(D)}, \max_{J \subset \{1, \ldots, d\}: J \ne \emptyset}   
\biggl(\int_{(\R_+^d)^J} (t^J)^{-\e^J -\theta (\alpha^J -\l^J)} 
(\Omega^{(l -\l) \chi_J}(\D^{\l \chi_J} f, t^J)_{L_p(D)})^\theta dt^J\biggr)^{1/\theta}\biggr), \\
f \in (S_{p,\theta}^\alpha B)^{\l}(D).
\end{multline*}

При $ \theta = \infty $ положим $ (S_{p,\infty}^\alpha B)^{\l}(D) =
(S_p^\alpha H)^{\l}(D). $

При тех же условиях на параметры обозначим через 
$ (S_{p,\theta}^\alpha B)^\prime(D) $ пространство всех функций 
$ f \in L_p(D), $ которые для любого непустого множества $ J \subset \{1,\ldots,d\} $ 
подчинены условию
\begin{multline*}
\int_{(\R_+^d)^J} (t^J)^{-\e^J -\theta \alpha^J} 
(\Omega^{\prime l \chi_J}(f, t^J)_{L_p(D)})^\theta dt^J = \\
\int_{(\R_+^d)^J} \biggl(\prod_{j \in J} t_j^{-1 -\theta \alpha_j}\biggr) 
(\Omega^{\prime l \chi_J}(f, t^{\s(l \chi_J)})_{L_p(D)})^\theta 
\prod_{j \in J} dt_j < \infty. 
\end{multline*}

В пространстве $ (S_{p,\theta}^\alpha B)^\prime(D) $ определяется норма
\begin{multline*}
\| f \|_{(S_{p,\theta}^\alpha B)^\prime(D)} = \\
\max\biggl(\| f\|_{L_p(D)}, 
\max_{J \subset \{1, \ldots, d\}: J \ne \emptyset}   
\biggl(\int_{(\R_+^d)^J} (t^J)^{-\e^J -\theta \alpha^J} 
(\Omega^{\prime l \chi_J}(f, t^J)_{L_p(D)})^\theta dt^J\biggr)^{1/\theta}\biggr),\\
 \ f \in (S_{p,\theta}^\alpha B)^\prime(D).
\end{multline*} 
При $ \theta = \infty $ положим $ (S_{p,\infty}^\alpha B)^\prime(D) =
(S_p^\alpha H)^\prime(D). $

В случае, когда вектор $ \l = \l(\alpha) \in \Z_+^d $ имеет компоненты
$ (\l(\alpha))_j = \max\{m \in \Z_+: m < \alpha_j\}, j =1,\ldots,d,$
пространство $ (S_{p,\theta}^\alpha B)^{\l}(D) ((S_p^\alpha H)^{\l}(D))$ обычно
обозначается $ S_{p,\theta}^\alpha B(D) (S_p^\alpha H(D)).$

Будет полезна

Лемма 1.1.3

Пусть $ d \in \N, l \in \Z_+^d, 1 \le p < \infty, D $ -- открытое множество в $ \R^d. $ 
Тогда для функции $ f \in L_p(D), $ у которой обобщённая производная 
$ \D^l f \in L_p(D), $ при $ h \in \R^d $ верно неравенство
\begin{equation*} \tag{1.1.7}
\| \Delta_h^l f \|_{L_p(D_h^l)} \le \biggl(\prod_{j \in \s(l)} |h_j|^{l_j}\biggr) 
\| \D^l f \|_{L_p(D)}.
\end{equation*}

Доказательство.

Прежде всего отметим, что, пользуясь формулой Ньютона-Лейбница, нетрудно 
показать, что при $ l \in \Z_+^d $ для любой функции $ f \in C^\infty(\R^d) $ 
при $ x, h \in \R^d $ выполняется равенство
\begin{multline*}
(\Delta_h^l f )(x) = \biggl(\prod_{j \in \s(l)} h_j^{l_j}\biggr) 
\int_{\prod_{j \in \s(l)} (\prod_{i_j =1}^{l_j} I)}
\D^l f \biggl(x +\sum_{j \in \s(l)} (\sum_{i_j =1}^{l_j} t_{j, i_j}) h_j e_j\biggr)
\prod_{j \in \s(l)} \biggl(\prod_{i_j =1}^{l_j} dt_{j, i_j}\biggr).
\end{multline*}
Отсюда заключаем, что при $ l \in \Z_+^d $ для любого ограниченного открытого 
множества $ D \subset \R^d $ и функции $ f \in C^\infty(\R^d) $ при $ h \in \R^d $
имеет место неравенство 
\begin{multline*}
\|(\Delta_h^l f )(\cdot) \|_{L_p(D_h^l)} = \\
\biggl(\prod_{j \in \s(l)} |h_j|^{l_j}\biggr) 
\biggl\|\int_{\prod_{j \in \s(l)} (\prod_{i_j =1}^{l_j} I)}
\D^l f \biggl(\cdot +\sum_{j \in \s(l)} (\sum_{i_j =1}^{l_j} t_{j, i_j}) h_j e_j\biggr)
\prod_{j \in \s(l)} \biggl(\prod_{i_j =1}^{l_j} dt_{j, i_j}\biggr) \biggr\|_{L_p(D_h^l)} \le \\
\biggl(\prod_{j \in \s(l)} |h_j|^{l_j}\biggr) 
\int_{\prod_{j \in \s(l)} (\prod_{i_j =1}^{l_j} I)}
\biggl\|\D^l f \biggl(\cdot +\sum_{j \in \s(l)} (\sum_{i_j =1}^{l_j} t_{j, i_j}) h_j e_j\biggr)
\biggr\|_{L_p(D_h^l)} \prod_{j \in \s(l)} \biggl(\prod_{i_j =1}^{l_j} dt_{j, i_j}\biggr) \le \\
\biggl(\prod_{j \in \s(l)} |h_j|^{l_j}\biggr) 
\int_{\prod_{j \in \s(l)} (\prod_{i_j =1}^{l_j} I)}
\|\D^l f \|_{L_p(D)} \prod_{j \in \s(l)} \biggl(\prod_{i_j =1}^{l_j} dt_{j, i_j}\biggr) = \\
\biggl(\prod_{j \in \s(l)} |h_j|^{l_j}\biggr) \|\D^l f \|_{L_p(D)},
\end{multline*}
что совпадает с (1.1.7).  

Далее, в условиях леммы возьмём неотрицательную функцию 
$ w \in C^\infty(\R^d), $ носитель которой $ \supp w \subset B^d, $ а 
$ \int_{\R^d} w(y) dy =1, $ и для функции $ f \in L_p(D), $ производная 
которой $ \D^l f \in L_p(D), $ при $ n \in \N $ определим функцию 
$ f_n \in C^\infty(\R^d) $ равенством
$$
f_n(x) = \int_D f(y) w_n(y -x) dy, x \in \R^d,
$$
где $ w_n(x) = n^d w(n x), x \in \R^d. $

Пусть компакт $ K \subset D. $ Тогда существует $ n_0 \in \N, $
для которого при $ n > n_0 $ для любого $ x \in K $ носитель 
$ \supp w_n(\cdot -x) \subset (x +n^{-1} B^d) \subset D, $ и, следовательно, 
при $ n > n_0 $ для любого $ x \in K $ справедливо равенство

\begin{multline*} 
\D^l f_n(x) = \int_D f(y) \frac{ \D^{|l|}} {\D x^{l}} (w_n(y -x)) dy = 
(-1)^{ |l|} \int_D f(y) \frac{ \D^{|l|}} {\D y^{l}} (w_n(y -x)) dy \\
= \int_D \D^{l} f(y) w_n(y -x) dy = \int_{ x +n^{-1} B^d} \D^{l} f(y) w_n(y -x) dy = \\
\int_{ n^{-1}  B^d} \D^{l} f(x +z) w_n(z) dz.
\end{multline*}
Отсюда, учитывая, что при $ n \in \N $ интеграл $ \int_{ n^{-1} B^d} w_n(z) dz 
= \int_{ \R^d} n^d w(n z) dz = \int_{ \R^d} w(z) dz =1, $
получаем, что для $ n > n_0 $ соблюдается соотношение
\begin{multline*} 
\| \D^{l} f -\D^{ l} f_n \|_{ L_p(K)}^p =
\int_K | \D^{ l} f(x) -\D^{ l} f_n(x)|^p dx \\
= \int_K \biggl| \int_{ n^{-1} B^d} \D^{ l} f(x) w_n(z) dz -
\int_{ n^{-1} B^d} \D^{ l} f(x +z) w_n(z) dz\biggr|^p dx = \\
\int_K \biggl|\int_{ n^{-1} B^d} (\D^{ l} f(x) -\D^{ l} f(x +z))
w_n(z) dz \biggr|^p dx \\
\le \int_K \biggl(\int_{ n^{-1} B^d} | \D^{ l} f(x) -\D^{ l} f(x +z)|^p w_n(z) dz\biggr) 
\biggl(\int_{ n^{-1} B^d} w_n(z)
 dz\biggr)^{p/p^{\prime}} dx \\
= \int_{ n^{-1} B^d} w_n(z) \int_K
 |\D^{ l} f(x) -\D^{ l} f(x +z)|^p dx dz\\
\le \biggl(\int_{ n^{-1} B^d} w_n(z) dz \biggr) \sup_{ z \in n^{-1} B^d} 
\| \D^{ l} f(\cdot) -\D^{ l} f(\cdot +z) \|_{ L_p(K)}^p \\
\le \sup_{ z \in n^{-1} B^d} \| \D^{ l} f(\cdot)
-\D^{ l} f(\cdot +z) \|_{ L_p(D_z)}^p \to 0, \text{ при } n \to \infty,
\text{ (см. [15])}.
\end{multline*}
Точно так же проверяется, что последовательность $ \{f_n, n \in \N\} $ 
сходится к $ f $ в $ L_p(K). $ 

Теперь с учётом сказанного выше, в условиях леммы для любого ограниченного 
открытого множества $ G, $ замыкание которого лежит в $ D, $ применяя (1.1.7) 
к каждой функции $ f_n,\ n \in \N,\ $ с множеством $ G $ вместо $ D, $ а затем
переходя к пределу при $ n \to \infty, $ приходим к выводу, что для функции $ f, $
удовлетворяющей условиям леммы, и любого ограниченного открытого множества $ G, $ 
замыкание которого лежит в $ D, $ имеет место (1.1.7) с заменой $ D $ на $ G. $
Наконец, беря последовательность ограниченных открытых множеств 
$ \{G_n :\ n \in \N\} $ со следующими свойствами:  
$ \bigcup_{n =1}^\infty G_n = D,\ G_n \subset G_{n +1}, $ замыкание $ G_n $ 
лежит в $ D,\ n \in \N, $ и применяя к функции $ f \in L_p(D), $ производная 
которой $ \D^l f \in L_p(D), $ неравенство (1.1.7) с $ G_n $ вместо 
$ D, $ а затем, переходя к пределу при  $ n \to \infty, $ приходим к (1.1.7)
в условиях леммы. $ \square $

В тех же обозначениях, что и выше, принимая во внимание то обстоятельство,
что для $ f \in (S_{p,\theta}^\alpha B)^{\l}(D) $ для $ J \subset 
\{1,\ldots,d\}: J \ne \emptyset, $ и $ t^J \in (\R_+^d)^J $
справедлива оценка (см., например, [16])
\begin{multline*}
(t^J)^{-(\alpha^J -\l^J)} \Omega^{(l -\l) \chi_J}(\D^{\l \chi_J} f, t^J)_{L_p(D)} = \\
\biggl(\int_{\{\tau^J \in (\R_+^d)^J: t_j < \tau_j < 2t_j, j \in
J\}} (t^J)^{-\e^J -\theta (\alpha^J -\l^J)} (\Omega^{(l -\l) \chi_J}(\D^{\l \chi_J} f,
t^J)_{L_p(D)})^\theta d\tau^J\biggr)^{1/\theta}  \le \\
\biggl(\int_{\{\tau^J \in (\R_+^d)^J:  t_j < \tau_j < 2t_j, j \in
J\}} (\prod_{j \in J} 2^{1+\theta (\alpha_j -\l_j)}) (\tau^J)^{-\e^J
-\theta (\alpha^J -\l^J)} (\Omega^{(l -\l) \chi_J}(\D^{\l \chi_J} f, \tau^J)_{L_p(D)})^\theta
d\tau^J\biggr)^{1/\theta} \le \\
(\prod_{j \in J} 2^{\alpha_j -\l_j +1/\theta}) \biggl(\int_{(\R_+^d)^J}
(\tau^J)^{-\e^J -\theta (\alpha^J -\l^J)} (\Omega^{(l -\l) \chi_J}(\D^{\l \chi_J} f,
\tau^J)_{L_p(D)})^\theta d\tau^J\biggr)^{1/\theta},
\end{multline*}
заключаем, что
\begin{multline*}
(S_{p,\theta}^\alpha B)^{\l}(D) \subset (S_p^\alpha H)^{\l}(D), \\
D \text{ -- область в } \R^d, \alpha \in \R_+^d, 1 \le p \le \infty,
1 \le \theta < \infty, \l \in \Z_+^d: \l < \alpha.
\end{multline*}

Учитывая, что для $ f \in (S_{p,\theta}^\alpha B)^{\prime}(D), 
t^J \in (\R_+^d)^J (J \subset \{1, \ldots, d\}: J \ne \emptyset) $ 
выполняется неравенство
\begin{multline*}
(t^J)^{-\alpha^J} \Omega^{\prime l \chi_J}(f, t^J)_{L_p(D)} = \\
\biggl(\int_{t^J +t^J (I^d)^J} (t^J)^{-\e^J -\theta \alpha^J} (\Omega^{\prime l \chi_J}(f,
t^J)_{L_p(D)})^\theta d\tau^J\biggr)^{1/\theta} = \\
\biggl(\int_{t^J +t^J (I^d)^J} (t^J)^{-\e^J -\theta \alpha^J} ((2 t^J)^{-\e^J}
\int_{ (t B^d)^J} \| \Delta_\xi^{l \chi_J} f\|_{L_p(D_\xi^{l \chi_J})}^p 
d \xi^J)^{\theta /p} d\tau^J\biggr)^{1/\theta} \le \\
\biggl(\int_{t^J +t^J (I^d)^J} (\prod_{j \in J} 2^{1 +\theta \alpha_j}) 
(\tau^J)^{-\e^J -\theta \alpha^J} ((\tau^J)^{-\e^J}
\int_{ (\tau B^d)^J} \| \Delta_\xi^{l \chi_J} f\|_{L_p(D_\xi^{l \chi_J})}^p 
d \xi^J)^{\theta /p} d\tau^J\biggr)^{1/\theta} = \\
\biggl(\int_{t^J +t^J (I^d)^J} (\prod_{j \in J} 2^{1 +\theta \alpha_j}) 
(\tau^J)^{-\e^J -\theta \alpha^J} (\prod_{j \in J} 2^{\theta /p})
(\Omega^{\prime l \chi_J}(f, \tau^J)_{L_p(D)})^\theta
d\tau^J\biggr)^{1/\theta} \le \\
(\prod_{j \in J} 2^{\alpha_j +1/\theta +1 /p}) \biggl(\int_{(\R_+^d)^J}
(\tau^J)^{-\e^J -\theta \alpha^J} (\Omega^{\prime l \chi_J}(f,
\tau^J)_{L_p(D)})^\theta d\tau^J\biggr)^{1/\theta},
\end{multline*}
заключаем, что
\begin{equation*} \tag{1.1.8}
(S_{p, \theta}^\alpha B)^\prime(D) \subset 
(S_p^\alpha H)^\prime(D),
\end{equation*}
и
\begin{equation*} 
\| f\|_{(S_p^\alpha H)^\prime(D)} \le c_4(\alpha)
\| f\|_{(S_{p, \theta}^\alpha B)^\prime(D)}, 
\end{equation*}
где $ c_4(\alpha) = \prod_{j=1}^d 2^{2+\alpha_j}. $

Из (1.1.6) и (1.1.7) следует, что для $ f \in (S_{p, \theta}^\alpha B)^{\l}(D), $
при $ \l \in \Z_+^d: \l < \alpha, l = l(\alpha), $ для 
$ J \subset \{1,\ldots,d\}: J \ne \emptyset, $ при $ t^J \in (\R_+^d)^J, $
справедливо неравенство

\begin{multline*}
\Omega^{\prime l \chi_J}(f, t^J)_{L_p(D)} \le
\Omega^{l \chi_J}(f, t^J)_{L_p(D)} = \\
\supvrai_{\xi \in \R^d: \xi^J \in 
(t B^d)^J} \| \Delta_\xi^{\l \chi_J +(l -\l) \chi_J} f
\|_{L_p(D_\xi^{\l \chi_J +(l -\l) \chi_J})} \le \\
\supvrai_{\xi \in \R^d: \xi^J \in 
(t B^d)^J} \| \Delta_\xi^{\l \chi_J}(\Delta_\xi^{(l -\l) \chi_J} f)
\|_{L_p((D_\xi^{(l -\l) \chi_J})_\xi^{\l \chi_J})} \le \\
\supvrai_{\xi \in \R^d: \xi^J \in (t B^d)^J} (\prod_{j \in \s(\l \chi_J)} |\xi_j|^{\l_j})
\| \D^{\l \chi_J} \Delta_\xi^{(l -\l) \chi_J} f
\|_{L_p(D_\xi^{(l -\l) \chi_J})} = \\
\supvrai_{\xi \in \R^d: \xi^J \in (t B^d)^J} (\prod_{j \in \s(\l \chi_J)} |\xi_j|^{\l_j})
\| \Delta_\xi^{(l -\l) \chi_J} \D^{\l \chi_J} f
\|_{L_p(D_\xi^{(l -\l) \chi_J})} \le \\
(\prod_{j \in J} t_j^{\l_j}) \supvrai_{\xi \in \R^d: \xi^J \in (t B^d)^J} 
\| \Delta_\xi^{(l -\l) \chi_J} \D^{\l \chi_J} f
\|_{L_p(D_\xi^{(l -\l) \chi_J})} = \\
(\prod_{j \in J} t_j^{\l_j}) \Omega^{(l -\l) \chi_J}(\D^{\l \chi_J} f, t^J)_{L_p(D)} =
(t^J)^{\l^J} \Omega^{(l -\l) \chi_J}(\D^{\l \chi_J} f, t^J)_{L_p(D)},
\end{multline*}
и, значит,
\begin{equation*} \tag{1.1.9}			
(S_{p, \theta}^\alpha B)^{\l}(D) \subset (S_{p, \theta}^\alpha B)^\prime(D) 
\end{equation*}
и
\begin{multline*} \tag{1.1.10}
\| f\|_{(S_{p, \theta}^\alpha B)^\prime(D)} \le
\| f\|_{(S_{p, \theta}^\alpha B)^{\l}(D)}, \\
\ f \in (S_{p, \theta}^\alpha B)^{\l}(D), \ 
\alpha \in \R_+^d, 1 \le p < \infty, 1 \le \theta \le \infty, \\
D \text{ -- произвольная область в } \R^d. \l \in \Z_+^d: \l < \alpha.
\end{multline*}
 
Обозначим через $ C^\infty(D) $ пространство бесконечно
дифференцируемых функций в области $ D \subset \R^d, $ а через 
$ C_0^\infty(D) $ -- пространство функций $ f \in C^\infty(\R^d), $
у которых носитель $ \supp f \subset D. $

В заключение этого пункта введём ещё несколько обозначений.

Для банахова пространства $ X $ (над $ \R$) обозначим $ B(X) = \{x
\in X: \|x\|_X \le 1\}. $

Для банаховых пространств $ X,Y $ через $ \mathcal B(X,Y) $
обозначим банахово пространство, состоящее из непрерывных линейных
операторов $ T: X \mapsto Y, $ с нормой
$$
\|T\|_{\mathcal B(X,Y)} = \sup_{x \in B(X)} \|Tx\|_Y.
$$
Отметим, что если $ X=Y,$ то $ \mathcal B(X,Y) $ является
банаховой алгеброй.
\bigskip

1.2. В этом пункте приведём некоторые вспомогательные утверждения,
которые используются в следующем пункте и далее.

В [16] содержится следующее утверждение.

   Лемма 1.2.1

Пусть $ d \in \N, 1 \le p < \infty. $ Тогда

1) при  $ j=1,\ldots,d$ для любого непрерывного линейного
оператора $ T: L_p(I) \mapsto L_p(I) $ существует единственный
непрерывный линейный оператор $ \mathcal T^j: L_p(I^d) \mapsto
L_p(I^d), $ для которого для любой функции $ f \in L_p(I^d) $
почти для всех $ (x_1,\ldots,x_{j-1},x_{j+1}, \ldots,x_d) \in
I^{d-1} $ в $ L_p(I) $ выполняется равенство
\begin{equation*} \tag{1.2.1}
(\mathcal T^j f)(x_1,\ldots, x_{j-1},\cdot,x_{j+1}, \ldots,x_d) =
(T(f(x_1,\ldots,x_{j-1},\cdot,x_{j+1}, \ldots,x_d)))(\cdot),
\end{equation*}

2) при этом, для каждого $ j=1,\ldots,d $ отображение $ V_j^{L_p},$ которое 
каждому оператору $ T \in \mathcal B(L_p(I), L_p(I)) $
ставит в соответствие оператор $ V_j^{L_p}(T) = \mathcal T^j \in
\mathcal B(L_p(I^d), L_p(I^d)), $ удовлетворяющий (1.2.1),
является непрерывным гомоморфизмом банаховой алгебры 
$ \mathcal B(L_p(I), L_p(I)) $ в банахову алгебру $ \mathcal B(L_p(I^d),
L_p(I^d)), $

3) причём для любых операторов $ S,T \in \mathcal B(L_p(I), L_p(I)) $ 
при любых $ i,j =1,\ldots,d: i \ne j, $ соблюдается равенство
\begin{equation*} \tag{1.2.2}
(V_i^{L_p}(S) V_j^{L_p}(T))f = (V_j^{L_p}(T) V_i^{L_p}(S))f, f \in
L_p(I^d).
\end{equation*}

Аналогичное утверждение установлено в [17].

Лемма 1.2.2

Пусть $ d \in \N, 1 \le p < \infty. $ Тогда

1) при  $ j=1,\ldots,d$ для любого непрерывного линейного оператора 
$ T: L_p(\R) \mapsto L_p(\R) $ существует единственный непрерывный 
линейный оператор $ \mathcal T^j: L_p(\R^d) \mapsto L_p(\R^d), $ для которого 
для любой функции $ f \in L_p(\R^d) $ почти для всех $ (x_1,\ldots,x_{j-1}, 
x_{j+1},\ldots,x_d) \in \R^{d-1} $ в $ L_p(\R) $ имеет место равенство
\begin{equation*} \tag{1.2.3}
(\mathcal T^j f)(x_1,\ldots, x_{j-1},\cdot,x_{j+1},\ldots,x_d) = 
(T(f(x_1,\ldots,x_{j-1},\cdot,x_{j+1},
\ldots,x_d)))(\cdot),
\end{equation*}

2) при этом, для каждого $ j=1,\ldots,d $ отображение $ V_j^{L_p}, $ которое
каждому оператору $ T \in \mathcal B(L_p(\R), L_p(\R)) $ ставит в соответствие
оператор $ V_j^{L_p}(T) = \mathcal T^j \in \mathcal B(L_p(\R^d), L_p(\R^d)), $
удовлетворяющий (1.2.3), является непрерывным гомоморфизмом банаховой алгебры
$ \mathcal B(L_p(\R), L_p(\R)) $ в банахову алгебру $ \mathcal B(L_p(\R^d),
L_p(\R^d)), $

3) причём для любых операторов $ S,T \in \mathcal B(L_p(\R), L_p(\R)) $
при любых $ i,j =1,\ldots,d: i \ne j, $ выполняется равенство
\begin{equation*} \tag{1.2.4}
(V_i^{L_p}(S) V_j^{L_p}(T))f = (V_j^{L_p}(T) V_i^{L_p}(S))f, f \in L_p(\R^d).
\end{equation*}

Отметим, что в дальнейшем символы $ L_p $ в обозначениях операторов $ V_j^{L_p} $
будем опускать.

Так же, как соответствующие усверждения в [16], устанавливаются 
леммы 1.2.3, 1.2.4.
    
Лемма 1.2.3

Пусть $ d \in \N, l \in \N^d, 1 \le p < \infty. $ Тогда существует константа
$ c_1(d,l) >0 $ такая, что для любых $ x^0 \in \R^d, \delta \in \R_+^d $ 
таких, что $ Q = (x^0 +\delta I^d) \subset I^d, $ для любой функции 
$ f \in L_p(I^d), $ для любого $ \xi \in \R^d, $ для любых множеств 
$ J, J^\prime \subset \{1,\ldots,d\} $ имеет место неравенство
\begin{multline*} \tag{1.2.5}
\biggl\| \Delta_\xi^{l \chi_J} ((\prod_{j \in J^\prime} V_j(E
-P_{\delta_j, x_j^0}^{1,l_j -1})) f)\biggr\|_{L_p(Q_\xi^{l \chi_J})} \le\\
c_1 \biggl(\prod_{j \in J^\prime \setminus J} \delta_j^{-1/p}\biggr)
\biggl(\int_{(\delta B^d)^{J^\prime \setminus J}} \int_{ Q_\xi^{l
\chi_{J \cup J^\prime}}} |(\Delta_\xi^{l \chi_{J \cup J^\prime}}
f)(x)|^p dx d\xi^{J^\prime \setminus J}\biggr)^{1/p},
\end{multline*}
где $ E $ -- тождественный оператор.

   Лемма 1.2.4

Пусть $ d \in \N, l \in \Z_+^d, 1 \le p < \infty $ и $ \rho, \sigma \in \R_+^d. $ 
Тогда существует константа $ c_2(d,l,\rho,\sigma) >0 $ такая, что для любых 
$ x^0, y^0 \in \R^d $ и $ \delta \in \R_+^d, $ для которых 
$ (x^0 +\sigma \delta I^d) \subset (y^0
+\rho \delta I^d) \subset I^d, $ для любого множества $ J \subset
\{1,\ldots,d\}, $ для $ f \in L_p(I^d) $ соблюдается неравенство
\begin{equation*} \tag{1.2.6}
\|(\prod_{j \in J} V_j(E -P_{\sigma_j \delta_j,x_j^0}^{1,l_j}))f\|_{L_p(D)} \le 
c_2 \|f\|_{L_p(D)},
\end{equation*}
где $ D = y^0 +\rho \delta I^d. $

В [12] доказана такая лемма.

Лемма 1.2.5

Пусть $ d \in \N, l \in \Z_+^d, $ а $ \delta \in \R_+^d, x^0 \in \R^d $ и 
$ Q = x^0 +\delta I^d $ таковы, что $ Q \subset I^d. $
Тогда для любой функции $ f \in L_1(I^d) $ почти  для всех $ x \in I^d $ имеет 
место равенство
\begin{equation*} \tag{1.2.7}
(P_{\delta, x^0}^{d,l}(f \mid_Q))(x) = ((\prod_{j=1}^d
V_j(P_{\delta_j, x_j^0}^{1,l_j}))f)(x).
\end{equation*}

Подобно доказательству леммы 1.2.5 проводится доказательство леммы 1.2.6,
для формулировки которой введём обозначения.
Для $ \phi \in L_\infty(\R) $ через $ M_\phi $ обозначим линейный оператор 
в пространстве локально суммируемых функций на $ \R, $ определяемый
равенством $ M_\phi g = \phi g, $ где $ g $ -- локально суммируемая на $ \R $ функция.
Через $ \chi_A $ будем обозначать характеристическую функцию множества 
$ A \subset \R^d. $

Лемма 1.2.6

Пусть $ d \in \N, l \in \Z_+^d, \delta, \Delta \in \R_+^d, x^0, X^0 \in \R^d. $ 
Тогда при $ 1 \le p < \infty $ для любой функции $ f \in L_p(\R^d) $ почти для 
всех $ x \in \R^d $ справедливо равенство
\begin{equation*} \tag{1.2.8}
\chi_{X^0 +\Delta I^d}(x) (P_{\delta, x^0}^{d,l}(f \mid_{x^0 +\delta I^d}))(x) 
= ((\prod_{j=1}^d V_j(M_{\chi_{X_j^0 +\Delta_j I}} 
P_{\delta_j, x_j^0}^{1,l_j}))f)(x).
\end{equation*}

Доказательство.

Прежде всего, построим ортонормированную в $ L_2(I) $ систему
полиномов $ \{\pi_m \in \mathcal P^{1,m}, m \in \Z_+\},$ полагая 
$ \pi_0 =1, $ а при $ m \in \N $ выбирая $ \pi_m \in \mathcal P^{1,m} $ так, 
чтобы $ \| \pi_m \|_{L_2(I)} =1 $ и для любого полинома $ f \in 
\mathcal P^{1, m-1} $ соблюдалось равенство 
$ \int_I \pi_m f dx =0. $

Затем при $ d \in \N, \lambda \in \Z_+^d $ определим
$$
\pi_\lambda^d(x) = \prod_{j=1}^d \pi_{\lambda_j}(x_j), x \in \R^d.
$$

Тогда при $ d \in \N, l \in \Z_+^d $ для $ \delta \in \R_+^d, x^0
\in \R^d, Q = x^0 +\delta I^d $ линейный оператор $ P_{\delta, x^0}^{d,l}: 
L_1(Q) \mapsto \mathcal P^{d,l}, $ задаётся равенством
\begin{multline*} \tag{1.2.9}
(P_{\delta, x^0}^{d,l}f)(x) = \delta^{-\e} \sum_{\lambda \in
\Z_+^d(l)} \biggl(\int_Q f(y) \pi_\lambda^d(\delta^{-1}(y -x^0)) dy\biggr)
\pi_\lambda^d(\delta^{-1}(x -x^0)),\\ 
f \in L_1(Q), x \in \R^d, \text{ (см. 
доказательство леммы 1.1.1 в [12])}.
\end{multline*}

Далее, введём следующее обозначение. При $ d \in \N $ для $ J \subset
\{1,\ldots,d\} $ обозначим через $ \eta_J: \R^d \times \R^d
\mapsto \R^d $ отображение, у которого
$$
(\eta_J(\xi,x))_j = \begin{cases} \xi_j, j \in J; \\
x_j, j \in \{1,\ldots,d\} \setminus J,\end{cases} \ \xi,x \in \R^d.
$$

Учитывая (1.2.9) и соотношение
$$
\chi_{X^0 +\Delta I^d}(x) = \prod_{j =1}^d \chi_{X_j^0 +\Delta_j I}(x_j), \ x \in \R^d,
$$
для доказательства (1.2.8) достаточно показать, что в условиях
леммы для $ f \in L_p(\R^d) $ для любого непустого множества 
$ J \subset \{1,\ldots,d\} $ почти для всех $ x \in \R^d $ справедливо
равенство
\begin{multline*} \tag{1.2.10}
((\prod_{j \in J} V_j(M_{\chi_{X_j^0 +\Delta_j I}} 
P_{\delta_j, x_j^0}^{1,l_j}))f)(x) = \\
(\prod_{j \in J} \chi_{X_j^0 +\Delta_j I}(x_j)) 
(\prod_{j \in J} \delta_j^{-\!1})\!\sum_{\lambda^J \in \Z_+^k(l^J)}
\biggl(\int_{Q^J} (\prod_{j \in J} \pi_{\lambda_j}(\delta_j^{-1} (y_j\!-\!x_j^0)))\,\cdot \\
\cdot f(\eta_J(y,x)) dy^J\biggr) (\prod_{j \in J}
\pi_{\lambda_j}(\delta_j^{-1} (x_j\!-\!x_j^0))),
\end{multline*}

где $ k = \card J. $ 

Равенство (1.2.10) установим по индукции относительно $ \card J. $ 
В случае, когда $ \card J =1, $ и, следовательно, $ J = \{j\}, $ 
где $ j \in \{1, \ldots, d\}, $ согласно (1.2.3), (1.2.9) при $ d =1, $ для 
$ f \in L_p(\R^d) $ почти для всех $ (x_1, \ldots, x_{j -1}, x_{j +1}, \ldots, x_d) \in \R^{d -1} $ 
для почти всех $ x_j \in \R $ имеет место равенство
\begin{multline*} 
(V_j(M_{\chi_{X_j^0 +\Delta_j I}} P_{\delta_j, x_j^0}^{1, l_j}) 
f)(x_1, \ldots, x_{j -1}, x_j, x_{j +1}, \ldots, x_d)   = \\
((M_{\chi_{X_j^0 +\Delta_j I}} P_{\delta_j, x_j^0}^{1, l_j})
f(x_1, \ldots, x_{j-1}, \cdot, x_{j+1}, \ldots, x_d))(x_j)  = \\
\chi_{X_j^0 +\Delta_j I}(x_j) (P_{\delta_j, x_j^0}^{1, l_j}
(f(x_1, \ldots, x_{j-1}, \cdot, x_{j+1}, \ldots, x_d) \mid_{x_j^0 +\delta_j I}))(x_j) = \\
\chi_{X_j^0 +\Delta_j I}(x_j) 
\delta_j^{-\!1} \!\sum_{\lambda_j =0}^{l_j}
\biggl(\int_{Q^{\{j\}}} \pi_{\lambda_j}(\delta_j^{-\!1} (y_j\!-\!x_j^0))\,\cdot\\
\cdot f(x_1, \ldots, x_{j -1}, y_j, x_{j +1}, \ldots, x_d) dy_j\biggr) 
\pi_{\lambda_j}(\delta_j^{-1} (x_j\!-\!x_j^0)),
\end{multline*}
которое влечёт (1.2.10) при $ J = \{j\}. $

Предположим теперь, что равенство (1.2.10) справедливо 
для любого множества $ J \subset \{1,\ldots,d\}, $ у которого 
$ \card J \le k, (1 \le k \le d-1), $ и покажем, что тогда оно верно для любого 
множества $ J \subset \{1, \ldots, d\}, $ у которого $ \card J = k+1. $ 

В этой ситуации, представляя $ J $ в виде $ J = J^\prime \cup \{j\}, $
где $ j \notin J^\prime, \card J^\prime = k,$ с учётом (1.2.4) на основании 
предположения индукции с использованием теоремы Фубини для $ f \in L_p(\R^d) $ 
почти для всех $ x \in \R^d $ получаем
\begin{multline*}
((\prod_{i \in J}
V_i(M_{\chi_{X_i^0 +\Delta_i I}} P_{\delta_i, x_i^0}^{1, l_i}))f)(x) = \\
(V_j(M_{\chi_{X_j^0 +\Delta_j I}} P_{\delta_j, x_j^0}^{1, l_j})
((\prod_{i \in J^\prime}
V_i(M_{\chi_{X_i^0 +\Delta_i I}} P_{\delta_i, x_i^0}^{1, l_i}))f))(x) = \\
(V_j(M_{\chi_{X_j^0 +\Delta_j I}} P_{\delta_j, x_j^0}^{1, l_j})
((\prod_{i \in J^\prime} \chi_{X_i^0 +\Delta_i I}(z_i)) 
(\prod_{i \in J^\prime} \delta_i^{-\!1})\\
\times\sum_{\lambda^{J^\prime}  \in
\Z_+^k(l^{J^\prime})}
\biggl(\int_{Q^{J^\prime}} (\prod_{i \in J^\prime} \pi_{\lambda_i}(\delta_i^{-\!1} (y_i\!-\!x_i^0)))
 f(\eta_{J^\prime}(y,z)) dy^{J^\prime}\biggr) (\prod_{i \in J^\prime}
\pi_{\lambda_i}(\delta_i^{-1} (z_i\!-\!x_i^0)))))(x) = \\
(\prod_{i \in J^\prime} \delta_i^{-\!1})\!\sum_{\lambda^{J^\prime}  \in
\Z_+^k(l^{J^\prime})} (V_j(M_{\chi_{X_j^0 +\Delta_j I}} 
P_{\delta_j, x_j^0}^{1, l_j})
\!((\prod_{i \in J^\prime} \chi_{X_i^0 +\Delta_i I}(z_i)) \\
\times\biggl(\int_{Q^{J^\prime}} (\prod_{i \in J^\prime} \pi_{\lambda_i}(\delta_i^{-\!1} (y_i\!-\!x_i^0)))
 f(\eta_{J^\prime}(y,z)) dy^{J^\prime}\biggr) (\prod_{i \in J^\prime}
\pi_{\lambda_i}(\delta_i^{-1} (z_i\!-\!x_i^0)))))(x) = \\
(\prod_{i \in J^\prime} \delta_i^{-\!1})\!\sum_{\lambda^{J^\prime}  \in
\Z_+^k(l^{J^\prime})} 
\chi_{X_j^0 +\Delta_j I}(x_j) 
\delta_j^{-\!1} \!\sum_{\lambda_j =0}^{l_j}
\biggl(\int_{Q^{\{j\}}} \pi_{\lambda_j}(\delta_j^{-\!1} (y_j\!-\!x_j^0))\,\cdot\\
\! (\prod_{i \in J^\prime} \chi_{X_i^0 +\Delta_i I}(x_i)) 
\biggl(\int_{Q^{J^\prime}} (\prod_{i \in J^\prime} \pi_{\lambda_i}(\delta_i^{-\!1} (y_i\!-\!x_i^0)))\,\cdot\\
\cdot f(\eta_{J^\prime}(y,\eta_{\{j\}}(y,x))) dy^{J^\prime}\biggr) (\prod_{i \in J^\prime}
\pi_{\lambda_i}(\delta_i^{-1} (x_i\!-\!x_i^0))) dy_j\biggr) 
\pi_{\lambda_j}(\delta_j^{-1} (x_j\!-\!x_j^0)) = \\
(\prod_{i \in J^\prime} \chi_{X_i^0 +\Delta_i I}(x_i))
\chi_{X_j^0 +\Delta_j I}(x_j) 
(\prod_{i \in J^\prime} \delta_i^{-\!1})\!
\delta_j^{-\!1} \\
\times\sum_{\lambda^{J^\prime} \in \Z_+^k(l^{J^\prime})} 
\sum_{\lambda_j =0}^{l_j}
\biggl(\int_{Q^{\{j\}}} \! \int_{Q^{J^\prime}} \pi_{\lambda_j}(\delta_j^{-\!1} (y_j\!-\!x_j^0))
(\prod_{i \in J^\prime} \pi_{\lambda_i}(\delta_i^{-\!1} (y_i\!-\!x_i^0)))\,\cdot\\
\cdot f(\eta_{J^\prime}(y,\eta_{\{j\}}(y,x))) dy^{J^\prime} dy_j\biggr)  
(\prod_{i \in J^\prime} \pi_{\lambda_i}(\delta_i^{-1} (x_i\!-\!x_i^0)))
\pi_{\lambda_j}(\delta_j^{-1} (x_j\!-\!x_j^0)) = \\
(\prod_{i \in J} \chi_{X_i^0 +\Delta_i I}(x_i))
(\prod_{i \in J} \delta_i^{-\!1})\\
\times\sum_{\lambda^J \in \Z_+^{k +1}(l^J)} 
\biggl(\int_{Q^J} (\prod_{i \in J} \pi_{\lambda_i}(\delta_i^{-\!1} (y_i\!-\!x_i^0)))
 f(\eta_J(y,x)) dy^J\biggr) (\prod_{i \in J} \pi_{\lambda_i}(\delta_i^{-1} (x_i\!-\!x_i^0))),
\end{multline*}
что и завершает вывод (1.2.10), которое. как отмечалось выше, влечёт (1.2.8). 
$ \square $
\bigskip

1.3. В этом пункте будут построены средства приближения для
функций из рассматриваемых нами пространств, на которые опирается вывод
основных результатов работы. Но сначала приведём некоторые вспомогательные 
сведения.

Для $ d \in \N, y \in \R^d $ положим
$$
\mn(y) = \min_{j=1,\ldots,d} y_j
$$
и для банахова пространства $ X, $ вектора $ x \in X $ и семейства
$ \{x_\kappa \in X, \kappa \in \Z_+^d\} $ будем писать $ x =
\lim_{ \mn(\kappa) \to \infty} x_\kappa, $ если для любого $ \epsilon >0 $ 
существует $ n_0 \in \N $ такое, что для любого $ \kappa \in \Z_+^d, $ 
для которого $ \mn(\kappa) > n_0, $
справедливо неравенство $ \|x -x_\kappa\|_X < \epsilon. $

Пусть $ X $ -- банахово пространство (над $ \R $), $ d \in \N $ и
$ \{ x_\kappa \in X: \kappa \in \Z_+^d\} $ -- семейство векторов.
Тогда под суммой ряда $ \sum_{\kappa \in \Z_+^d} x_\kappa $ будем
понимать вектор $ x \in X, $ для которого выполняется равенство 
$ x = \lim_{\mn(k) \to \infty} \sum_{\kappa \in \Z_+^d(k)} x_\kappa. $

При $ d \in \N $ через $ \Upsilon^d $ обозначим множество
$$
\Upsilon^d = \{ \epsilon \in \Z^d: \epsilon_j \in \{0,1\},
j=1,\ldots,d\}.
$$

Как показано в [12], имеет место следующая лемма.

   Лемма 1.3.1

Пусть $ X $ -- банахово пространство, а вектор $ x \in X $ и
семейство $ \{x_\kappa \in X: \kappa \in \Z_+^d\} $ таковы, что 
$ x = \lim_{ \mn(\kappa) \to \infty} x_\kappa, $ Тогда для семейства 
$ \{ \mathcal X_\kappa \in X, \kappa \in \Z_+^d \}, $ элементы которого определяются
равенством
$$
\mathcal X_\kappa = \sum_{\epsilon \in \Upsilon^d: \s(\epsilon)
\subset \s(\kappa)} (-\e)^\epsilon x_{\kappa -\epsilon}, \kappa
\in \Z_+^d,
$$
справедливо равенство
$$
x = \sum_{\kappa \in \Z_+^d} \mathcal X_\kappa.
$$

Замечание

Как отмечалось в [12], если для семейства векторов $ \{x_\kappa \in X, 
\kappa \in \Z_+^d\} $ банахова пространства $ X $ ряд 
$ \sum_{\kappa \in \Z_+^d} \| x_\kappa \|_X $ сходится, то для любой 
последовательности подмножеств $ \{Z_n \subset \Z_+^d, n \in \Z_+\}, $
таких, что $ \card Z_n < \infty, Z_n \subset Z_{n+1}, n \in \Z_+, $
и $ \cup_{ n \in \Z_+} Z_n = \Z_+^d, $
в $ X $ соблюдается равенство
$ \sum_{\kappa \in \Z_+^d} x_\kappa =
\lim_{ n \to \infty} \sum_{\kappa \in Z_n} x_\kappa. $

Введём в рассмотрение систему разбиений единицы на кубе $ I^d, $
используемую для построения средств приближения функций из изучаемых 
пространств. Для этого обозначим через $ \psi^{1,0} $ характеристическую функцию
интервала $ I, $ т.е. функцию, определяемую равенством
$$
\psi^{1,0}(x) = \begin{cases} 1, & \text{ для } x \in I; \\
0, & \text{ для } x \in \R \setminus I.
\end{cases}
$$
При $ m \in \N $ положим
$$
\psi^{1,m}(x) = \int_I \psi^{1, m-1}(x-y) dy \ (\text{см. } [18]),
$$
а для $ d \in \N, m \in \Z_+^d $ определим
$$
\psi^{d,m}(x) = \prod_{j=1}^d \psi^{1,m_j}(x_j), x =
(x_1,\ldots,x_d) \in \R^d.
$$

Для $ d \in \N, m,n \in \Z^d: m \le n, $ обозначим
$$
\Nu_{m,n}^d = \{ \nu \in \Z^d: m \le \nu \le n \} = \prod_{j=1}^d
\Nu_{m_j,n_j}^1.
$$

Опираясь на определения, используя индукцию, нетрудно проверить
следующие свойства функций $ \psi^{d,m}, d \in \N, m \in \Z_+^d. $

1) При $ d \in \N, m \in \Z_+^d $
$$
\sgn \psi^{d,m}(x) = \begin{cases} 1, \text{ для } x \in ((m+\e) I^d); \\
0, \text{ для } x \in \R^d \setminus ((m+\e) I^d),
\end{cases}
$$

2) при $ d \in \N, m \in \Z_+^d $ для каждого $ \lambda \in
\Z_+^d(m) $ (обобщённая) производная $ \D^\lambda \psi^{d,m} \in
L_\infty(\R^d), $

3) при $ d \in \N, m \in \Z_+^d $ почти для всех $ x \in \R^d $
справедливо равенство
$$
\sum_{\nu \in \Z^d} \psi^{d,m}(x -\nu) =1,
$$

4) при $ m \in \N $ для всех $ x \in \R $ (при $ m =0 $ почти для всех $ x \in \R $) 
имеет место равенство
\begin{equation*} \tag{1.3.1}
\psi^{1,m}(x) = \sum_{\mu \in \Nu_{0, m+1}^1} a_{\mu}^m
\psi^{1,m}(2x -\mu),
\end{equation*}
где $ a_\mu^m = 2^{-m} C_{m+1}^\mu, \mu \in \Nu_{0, m+1}^1. $

При $ d \in \N $ для $ t \in \R^d $ через $ 2^t $ будем обозначать
вектор $ 2^t = (2^{t_1}, \ldots, 2^{t_d}). $

Для $ d \in \N, m,\kappa \in \Z_+^d, \nu \in \Z^d $ обозначим
$$
g_{\kappa, \nu}^{d,m}(x) = \psi^{d,m}(2^\kappa x -\nu) =
\prod_{j=1}^d \psi^{1,m_j}( 2^{\kappa_j} x_j -\nu_j), x \in \R^d,
$$
$$
Q_{\kappa, \nu}^d = 2^{-\kappa} \nu +2^{-\kappa} I^d.
$$

Из первого среди приведенных выше свойств функций $ \psi^{d,m} $
следует, что при $ d \in \N, m,\kappa \in \Z_+^d, \nu \in \Z^d $
носитель $ \supp g_{\kappa,\nu}^{d,m} = 
2^{-\kappa} \nu +2^{-\kappa} (m+\e) \overline I^d. $

Отметим некоторые полезные для нас свойства носителей функций 
$ g_{\kappa,\nu}^{d,m}. $

Для $ d \in \N, m,\kappa \in \Z_+^d $ справедливо следующее:
\begin{equation*} \tag{1.3.2}
1) \{\nu \in \Z^d: \supp g_{\kappa,\nu}^{d,m} \cap I^d \ne \emptyset \} = 
\Nu_{-m, 2^\kappa -\e}^d,
\end{equation*}

2) для каждого $ \nu^\prime \in \Z^d $ имеет место равенство
\begin{equation*} \tag{1.3.3}
\{ \nu \in \Z^d: Q_{\kappa, \nu^\prime}^d \cap 
\supp g_{\kappa, \nu}^{d,m} \ne \emptyset\} = \nu^\prime +\Nu_{-m,0}^d.
\end{equation*}

Из свойства 3) функций $ \psi^{d,m} $ и равенства (1.3.2)
вытекает, что при $ d \in \N, m,\kappa \in \Z_+^d $ почти для всех
$ x \in I^d $ соблюдается равенство
\begin{equation*} \tag{1.3.4}
\sum_{ \nu \in \Nu_{-m, 2^\kappa -\e}^d} g_{\kappa, \nu}^{d,m}(x) =1.
\end{equation*}

Используя разложение Ньютона для $ (1+1)^{m+1} $ и $ (-1+1)^{m+1}, $ легко
проверить, что при $ m \in \Z_+ $ выполняются равенства
\begin{equation*} 
\sum_{\mu \in \Nu_{0,m +1}^1 \cap (2 \Z)} a_\mu^m =1,
\sum_{\mu \in \Nu_{0,m +1}^1 \cap (2 \Z +1)} a_\mu^m =1.
\end{equation*}
Отсюда, принимая во внимание, что при $ \kappa \in \N, m \in \Z_+, $ для каждого
$ \nu \in \Nu_{-m, 2^\kappa -1}^1 $ верно равенство
\begin{multline*}
\{\mu \in \Nu_{0, m+1}^1: (\nu -\mu)/2 \in \Nu_{-m, 2^{\kappa -1} -1}^1 \} = \\
\{\mu \in \Nu_{0, m+1}^1: (\nu -\mu)/2 \in \Z\} = \Nu_{0, m+1}^1 \cap (\nu +(2 \Z)) =\\
\begin{cases} 
\Nu_{0, m+1}^1 \cap (2 \Z), \text{ если } \nu \in (2 \Z); \\
\Nu_{0, m+1}^1 \cap ((2 \Z) +1), \text{ если } \nu \in ((2 \Z) +1),
\end{cases}
\end{multline*}
заключаем, что при указанных условиях имеет место равенство
\begin{equation*} \tag{1.3.5}
\sum_{\mu \in \Nu_{0, m+1}^1: (\nu
-\mu)/2 \in \Nu_{-m, 2^{\kappa -1} -1}^1 } a_\mu^m = \sum_{ \rho
\in \Nu_{-m, 2^{\kappa -1} -1}^1: (\nu -2\rho) \in \Nu_{0, m+1}^1}
a_{\nu -2\rho}^m =1.
\end{equation*}

При $ d \in \N $ для $ x,y \in \R^d $ будем обозначать
$$
(x,y) = \sum_{j=1}^d x_j y_j.
$$

Имея в виду свойство 2) функций $ \psi^{d,m}, $ отметим, что при 
$ d \in \N, m,\kappa \in \Z_+^d, \nu \in \Z^d, \lambda \in \Z_+^d(m)$ 
выполняется равенство
\begin{multline*} \tag{1.3.6}
\| \D^\lambda g_{\kappa, \nu}^{d,m} \|_{L_\infty (\R^d)} =
2^{(\kappa, \lambda)} \| \D^\lambda \psi^{d,m} \|_{L_\infty(\R^d)} = 
c_1(d,m,\lambda) 2^{(\kappa, \lambda)}.
\end{multline*}

Для формулировки леммы 1.3.2 введём следующие обозначения.

При $ d \in \N $ для $ x \in \R^d $ положим
$$
x_+ = ((x_1)_+, \ldots, (x_d)_+),
$$
где $ t_+ = \frac{1} {2} (t +|t|), t \in \R. $

При $ d \in \N, m, \kappa \in \Z_+^d, \nu \in \Nu_{0,2^\kappa -\e}^d $ положим
$$
x_{\kappa, \nu}^{d,m} = 2^{-\kappa} (\nu -m)_+,
$$
а для каждого $ j=1,\ldots, d $ зададим $ j $-ую координату
вектора $ \delta_{\kappa, \nu}^{d,m} $ равенством
$$
(\delta_{\kappa, \nu}^{d,m})_j = \min(2^{\kappa_j}, m_j +1)
2^{-\kappa_j},
$$
и обозначим через
$$
D_{\kappa, \nu}^{d,m} = x_{\kappa, \nu}^{d,m} +\delta_{\kappa, \nu}^{d,m} I^d.
$$

Заметим, что при $ d \in \N, \kappa \in \Z_+^d, \nu \in \Nu_{0,2^\kappa -\e}^d $ 
верно равенство $ D_{\kappa, \nu}^{d, 0} =
Q_{\kappa, \nu}^d, $ а при $ d \in \N, m,\kappa \in \Z_+^d,
\nu \in \Nu_{0, 2^\kappa -\e}^d $ имеет место включение
\begin{equation*} \tag{1.3.7}
D_{\kappa, \nu}^{d,m} \subset I^d.
\end{equation*}

Для $ d \in \N, l,m,\kappa \in \Z_+^d, \nu \in \Nu_{0, 2^\kappa -\e}^d $ 
определим непрерывный линейный оператор 
$ S_{\kappa,\nu}^{d,l,m}: L_1(I^d) \mapsto \mathcal P^{d,l}(I^d) \cap 
L_\infty(I^d), $ полагая для $ f \in L_1(I^d) $ значение
$$
(S_{\kappa, \nu}^{d,l,m} f)(x) = (P_{\delta, x^0}^{d,l} (f \mid_{(x^0
+\delta I^d)}))(x), x \in I^d,
$$
при $ \delta = \delta_{\kappa, \nu}^{d,m}, x^0 = x_{\kappa,
\nu}^{d,m} $ (см. лемму 1.1.2 и (1.3.7)).

При $ d \in \N, l, \kappa \in \Z_+^d, \nu \in \Nu_{0, 2^\kappa -\e}^d $ 
обозначим
$$
\mathcal S_{\kappa, \nu}^{d,l} = S_{\kappa, \nu}^{d,l,0}.
$$

Далее, при $ d \in \N, m, \kappa \in \Z_+^d $ для $ \nu \in \Nu_{-m, 2^\kappa -\e}^d $
зададим $ \nu_\kappa^{d,m}(\nu) = \nu_\kappa(\nu), $ полагая
$$
\nu_\kappa^{d,m}(\nu) = (2^\kappa -m -\e)_+ -(2^\kappa -m -\e -\nu_+)_+,
$$
и определим при $ d \in \N, l,m,\kappa \in \Z_+^d $ линейный непрерывный
оператор $ E_\kappa^{d,l,m}: L_1(I^d) \mapsto L_\infty(I^d) $
равенством
\begin{equation*}
E_\kappa^{d,l,m} f = \sum_{\nu \in \Nu_{-m, 2^\kappa -\e}^d}
(\mathcal S_{\kappa, \nu_\kappa^{d,m}(\nu)}^{d,l} f) g_{\kappa, \nu}^{d,m} \mid_{I^d}, 
f \in L_1(I^d).
\end{equation*}

Аналогом соответствующего утверждения из [16] является лемма 1.3.2.

Лемма 1.3.2

Пусть $ d \in \N, l \in \N^d, m \in \Z_+^d, 1 \le p < \infty. $
Тогда существуют константы $ c_2(d,l,m) >0 $ и $ c_3(d,m) >0 $ 
такие, что для любой функции $ f \in L_p(I^d) $ при $ \kappa \in
\Z_+^d $ справедливо неравенство
\begin{equation*} \tag{1.3.8}
\|f -E_\kappa^{d,l -\e,m}f \|_{L_p(I^d)} \le c_2 \sum_{j=1}^d
2^{\kappa_j /p} \biggl(\int_{ c_3 2^{-\kappa_j} B^1} 
\int_{(I^d)_{l_j \xi e_j}} |\Delta_{\xi e_j}^{l_j} f(x)|^p dx d\xi\biggr)^{1/p}.
\end{equation*}

Доказательство.

Сначала отметим свойства множеств $ D_{\kappa, \nu}^{d,m}, $
которые понадобятся для доказательства леммы.

При $ d \in \N, m,\kappa \in \Z_+^d, \nu \in \Nu_{0, 2^\kappa -\e}^d $ 
справедливо включение
\begin{equation*} \tag{1.3.9}
D_{\kappa, \nu}^{d,0} \subset D_{\kappa, \nu}^{d,m}.
\end{equation*}

Нетрудно видеть, что при $ d \in \N, m \in \Z_+^d $ существует
константа $ c_4(d,m) >0 $ такая, что для любого $ \kappa \in \Z_+^d $ для 
каждого $ x \in I^d $ число
\begin{equation*} \tag{1.3.10}
\card \{ \nu \in \Nu_{0, 2^\kappa -\e}^d: x \in D_{\kappa,\nu}^{d,m} \} \le c_4.
\end{equation*}

С помощью (1.3.3) несложно установить, что если при $ d \in \N,
m,\kappa \in \Z_+^d $ для $ \nu \in \Nu_{-m, 2^\kappa -\e}^d $ и 
$ \nu^\prime \in \Nu_{0, 2^\kappa -\e}^d $ пересечение 
$ \supp g_{\kappa, \nu}^{d,m} \cap Q_{\kappa, \nu^\prime}^d \ne \emptyset, $ то
\begin{equation*} \tag{1.3.11}
D_{\kappa, \nu_\kappa^{d,m}(\nu)}^{d,0} \subset D_{\kappa, \nu^\prime}^{d,m}.
\end{equation*}

Из (1.3.3) также вытекает, что
\begin{multline*} \tag{1.3.12}
\card \{ \nu \in \Nu_{-m, 2^\kappa -\e}^d: \supp g_{\kappa,\nu}^{d,m} \cap 
Q_{\kappa, \nu^\prime}^d \ne \emptyset \} \le \\
c_5(d,m), \ d \in \N, m,\kappa \in \Z_+^d, \nu^\prime \in \Nu_{0, 2^\kappa -\e}^d.
\end{multline*}

Пусть $ f \in L_p(I^d) $ и $ \kappa \in \Z_+^d. $ Тогда в виду (1.3.4) имеем
\begin{multline*} \tag{1.3.13}
\|f -E_\kappa^{d,l -\e,m}f \|_{L_p(I^d)}^p \\
\le \sum_{ \nu^\prime \in \Nu_{0, 2^\kappa -\e}^d} \biggl(\sum_{\nu \in 
\Nu_{-m, 2^\kappa -\e}^d: \supp g_{\kappa, \nu}^{d,m} \cap
Q_{\kappa, \nu^\prime}^d \ne \emptyset} 
\|(f -\mathcal S_{\kappa, \nu_\kappa(\nu)}^{d, l-\e} f) g_{\kappa, \nu}^{d,m} 
\|_{L_p(Q_{\kappa, \nu^\prime}^d)} \biggr)^p.
\end{multline*}

Далее, для $ f \in L_p(I^d) $ при $ \kappa \in \Z_+^d $ для $ \nu
\in \Nu_{-m, 2^\kappa -\e}^d $ и $ \nu^\prime \in \Nu_{0, 2^\kappa -\e}^d $ 
таких, что $ \supp g_{\kappa, \nu}^{d,m} \cap Q_{\kappa, \nu^\prime}^d 
\ne \emptyset, $ выводим
\begin{multline*} \tag{1.3.14}
\|(f -\mathcal S_{\kappa, \nu_\kappa(\nu)}^{d, l-\e} f) g_{\kappa,\nu}^{d,m} 
\|_{L_p( Q_{\kappa, \nu^\prime}^d)} \le \\
\|f -S_{\kappa, \nu^\prime}^{d, l-\e,m}f \|_{L_p( Q_{\kappa, \nu^\prime}^d)} + 
\|S_{\kappa, \nu^\prime}^{d, l-\e,m}f
-\mathcal S_{\kappa, \nu_\kappa(\nu)}^{d, l-\e} f \|_{L_p( Q_{\kappa, \nu^\prime}^d)}.
\end{multline*}

Учитывая (1.3.7), (1.3.9), на основании (1.1.5) для $ f \in L_p(I^d) $ 
при $ \kappa \in \Z_+^d $ и $ \nu^\prime \in \Nu_{0, 2^\kappa -\e}^d $ 
получаем
\begin{multline*} \tag{1.3.15}
\|f -S_{\kappa, \nu^\prime}^{d, l-\e,m}f \|_{L_p( Q_{\kappa, \nu^\prime}^d)} \\
\le c_6 \sum_{j=1}^d 2^{\kappa_j /p} \biggl(\int_{ c_3 2^{-\kappa_j} B^1} 
\int_{ (D_{\kappa, \nu^\prime}^{d,m})_{l_j \xi e_j}} 
|\Delta_{\xi e_j}^{l_j} f(x)|^p dx d\xi\biggr)^{1/p}.
\end{multline*}

Принимая во внимание (1.3.9), (1.3.11),  (1.1.2),
(1.1.3), (1.1.4) и снова (1.3.11), а также учитывая (1.1.5),
для $ f \in L_p(I^d) $ при $ \kappa \in \Z_+^d $ для $ \nu \in
\Nu_{-m, 2^\kappa -\e}^d $ и $  \nu^\prime \in \Nu_{0, 2^\kappa -\e}^d: 
\supp g_{\kappa, \nu}^{d,m} \cap Q_{\kappa, \nu^\prime}^d \ne \emptyset, $ 
находим, что
\begin{multline*} \tag{1.3.16}
\|S_{\kappa, \nu^\prime}^{d, l-\e,m}f -
\mathcal S_{\kappa,\nu_\kappa(\nu)}^{d, l-\e} f \|_{L_p( Q_{\kappa, \nu^\prime}^d)} \le
c_7 \|\mathcal S_{\kappa, \nu_\kappa(\nu)}^{d, l-\e} 
(S_{\kappa,\nu^\prime}^{d, l-\e,m}f -f) \|_{L_p( D_{\kappa, \nu_\kappa(\nu)}^{d, 0})}
\\ \le c_8 \sum_{j=1}^d 2^{\kappa_j /p} \biggl(\int_{ c_3 2^{-\kappa_j} B^1} 
\int_{ (D_{\kappa, \nu^\prime}^{d,m})_{l_j \xi e_j}} 
|\Delta_{\xi e_j}^{l_j} f(x)|^p dx d\xi\biggr)^{1/p}.
\end{multline*}

Соединяя (1.3.14), (1.3.15) и (1.3.16), для $ f \in L_p(I^d),
\kappa \in \Z_+^d, \nu \in \Nu_{-m, 2^\kappa -\e}^d, \nu^\prime
\in \Nu_{0, 2^\kappa -\e}^d: \supp g_{\kappa, \nu}^{d,m} \cap
Q_{\kappa, \nu^\prime}^d \ne \emptyset, $ получаем
$$
\|(f -\mathcal S_{\kappa, \nu_\kappa(\nu)}^{d, l-\e} f) g_{\kappa,\nu}^{d,m} 
\|_{L_p( Q_{\kappa, \nu^\prime}^d)} \le c_9 \sum_{j=1}^d 2^{\kappa_j /p} 
\biggl(\int_{ c_3 2^{-\kappa_j} B^1} 
\int_{ (D_{\kappa, \nu^\prime}^{d,m})_{l_j \xi e_j}}
|\Delta_{\xi e_j}^{l_j} f(x)|^p dx d\xi\biggr)^{1/p}.
$$

Подставляя эту оценку в (1.3.13) и применяя (1.3.12) и неравенство
Гёльдера, а затем используя (1.3.10) и (1.1.1) при $ a = 1/p, $ для 
$ f \in L_p(I^d) $ при $ \kappa \in \Z_+^d $ выводим
\begin{multline*}
\|f -E_\kappa^{d,l -\e,m}f \|_{L_p(I^d)}^p \le \\
c_{10}^p \sum_{j=1}^d 2^{\kappa_j } \int_{ c_3 2^{-\kappa_j} B^1} 
\int_{ (I^d)_{l_j \xi e_j}} \biggl(\sum_{ \nu \in \Nu_{0, 2^\kappa -\e}^d} 
\chi_{D_{\kappa, \nu}^{d,m}}(x)\biggr) |\Delta_{\xi e_j}^{l_j}f(x)|^p dx d\xi \\
\le \biggl(c_2 \sum_{j=1}^d 2^{\kappa_j /p} \biggl(\int_{c_3 2^{-\kappa_j} B^1} 
\int_{ (I^d)_{l_j \xi e_j}} |\Delta_{\xi e_j}^{l_j} f(x)|^p dx d\xi\biggr)^{1/p}\biggr)^p, 
\end{multline*}
что влечёт (1.3.8). $ \square $

Как известно, имеет место

    Лемма 1.3.3

Пусть $ d \in \N, \lambda \in \Z_+^d, D $ --- область в $ \R^d $ и
функция $ f \in C^\infty(D), $ а  $  g  \in  L_1(D), $ причём,
для каждого $ \mu \in \Z_+^d(\lambda) $ обобщённая производная 
$ \D^\mu g \in L_1(D). $ Тогда в пространстве обобщённых функций в
области $ D $ справедливо равенство
\begin{equation*} \tag{1.3.17}
\D^\lambda (fg) = \sum_{ \mu  \in  \Z_+^d(\lambda)} C_\lambda^\mu
\D^{\lambda -\mu}f \D^\mu g.
\end{equation*}

Для формулировки следующего утверждения введём обозначение.

При $ d \in \N, l,m,\kappa \in \Z_+^d $ определим линейный непрерывный оператор 
$ \mathcal E_\kappa^{d,l,m}: L_1(I^d) \mapsto L_\infty(I^d), $ полагая
$$
\mathcal E_\kappa^{d,l,m} = \sum_{\epsilon \in \Upsilon^d:
\s(\epsilon) \subset \s(\kappa)} (-\e)^\epsilon E_{\kappa -\epsilon}^{d,l,m}.
$$

Преобразуем выражение, задающее оператор $ \mathcal E_\kappa^{d,l,m}, $ к виду, 
подходящему для получения интересующих нас оценок.
Для этого при $ d \in \N, l,m,\kappa \in \Z_+^d $ для $ f \in L_1(I^d) $ 
почти для всех $ x \in I^d $ имеем
\begin{equation*} \tag{1.3.18}
(\mathcal E_\kappa^{d,l,m} f)(x) = \sum_{\epsilon \in \Upsilon^d:
\s(\epsilon) \subset \s(\kappa)} (-\e)^\epsilon \sum_{\rho \in
\Nu_{-m, 2^{\kappa -\epsilon} -\e}^d} (\mathcal S_{\kappa -\epsilon, 
\nu_{\kappa -\epsilon}^{d,m}(\rho)}^{d,l} f)(x) g_{\kappa -\epsilon,
\rho}^{d,m}(x).
\end{equation*}

Используя (1.3.1), при $ d \in \N, m,\kappa \in \Z_+^d, \epsilon
\in \Upsilon^d: \s(\epsilon) \subset \s(\kappa),  $ и $ \rho \in
\Nu_{-m, 2^{\kappa -\epsilon} -\e}^d $ почти для всех $ x \in I^d $ получаем
\begin{multline*}
g_{\kappa -\epsilon, \rho}^{d,m}(x) = \\
 \sum_{ \mu^{\s(\epsilon)} \in (\Nu_{0, m+\e}^d)^{\s(\epsilon)}}
(\prod_{j \in \s(\epsilon)} a_{\mu_j}^{m_j}) (\prod_{j \in
\Nu_{1,d}^1 \setminus \s(\epsilon)} \psi^{1, m_j}(2^{\kappa_j} x_j
-\rho_j)) (\prod_{j \in \s(\epsilon)} \psi^{1, m_j}(2^{\kappa_j}
x_j -2 \rho_j -\mu_j)).
\end{multline*}

Отсюда, учитывая (1.3.2), при $ d \in \N, l, m,\kappa \in \Z_+^d,
\epsilon \in \Upsilon^d: \s(\epsilon) \subset \s(\kappa),  $ почти
для всех $ x \in I^d $ имеем
\begin{multline*} \tag{1.3.19}
\sum_{\rho \in \Nu_{-m, 2^{\kappa -\epsilon} -\e}^d} 
(\mathcal S_{\kappa -\epsilon, \nu_{\kappa -\epsilon}^{d,m}(\rho)}^{d,l} f)(x) 
g_{\kappa -\epsilon, \rho}^{d,m}(x) =\\
 \sum_{ \nu \in \Nu_{-m, 2^\kappa -\e}^d} g_{\kappa, \nu}^{d,m}(x)
\sum_{ \rho \in \Rho_{\kappa,\nu,\epsilon}^{d,m}} (\prod_{j \in
\s(\epsilon)} a_{\nu_j -2\rho_j}^{m_j}) (\mathcal S_{\kappa -\epsilon, 
\nu_{\kappa -\epsilon}^{d,m}(\rho)}^{d,l} f)(x), 
\end{multline*}
где
$$
\Rho_{\kappa,\nu,\epsilon}^{d,m} = \{ \rho \in \Nu_{-m, 2^{\kappa -\epsilon} -\e}^d: 
\rho_j = \nu_j, j \in \Nu_{1,d}^1 \setminus \s(\epsilon);  
(\nu_j -2 \rho_j) \in \Nu_{0, m_j +1}^1, j \in \s(\epsilon) \}.
$$

Подставляя (1.3.19) в (1.3.18), выводим
\begin{equation*} \tag{1.3.20} 
(\mathcal E_\kappa^{d,l,m} f)(x)
= \sum_{ \nu \in \Nu_{-m, 2^\kappa -\e}^d} g_{\kappa, \nu}^{d,m}(x)
(U_{\kappa,\nu}^{d,l,m} f)(x), \text{ почти для всех } x \in I^d,
\end{equation*}

$ d \in \N, l,m,\kappa \in \Z_+^d, f \in L_1(I^d), $ 
где $ U_{\kappa,\nu}^{d,l,m}: L_1(I^d) \mapsto \mathcal P^{d,l}(I^d) $ -- 
линейный оператор, определяемый при $ d \in \N, l,m,\kappa \in
\Z_+^d, \nu \in \Nu_{-m, 2^\kappa -\e}^d $ для $ f \in L_1(I^d) $ равенством
\begin{equation*} \tag{1.3.21}
U_{\kappa,\nu}^{d,l,m} f = \sum_{\epsilon \in \Upsilon^d:
\s(\epsilon) \subset \s(\kappa)} (-\e)^\epsilon \sum_{ \rho \in
\Rho_{\kappa,\nu,\epsilon}^{d,m}} (\prod_{j \in \s(\epsilon)}
a_{\nu_j -2\rho_j}^{m_j}) (\mathcal S_{\kappa -\epsilon,
\nu_{\kappa -\epsilon}^{d,m}(\rho)}^{d,l} f).
\end{equation*}

В силу равенства (см. (1.2.7))
$$
\mathcal S_{\kappa, \nu}^{d,l} = 
\prod_{j=1}^d V_j(S_{\kappa_j, \nu_j}^{1,l_j,0}), \  
d \in \N, l,\kappa \in \Z_+^d, \nu \in \Nu_{0, 2^\kappa -\e}^d,
$$
благодаря п. 2 леммы 1.2.1 и (1.2.2), при $ d \in \N, l,m,\kappa \in \Z_+^d, 
\nu \in \Nu_{-m, 2^\kappa -\e}^d, \epsilon \in \Upsilon^d:
\s(\epsilon) \subset \s(\kappa),$ 
для $ \rho \in \Rho_{\kappa,\nu,\epsilon}^{d,m} $ выполняется равенство
\begin{multline*} \tag{1.3.22}
\mathcal S_{\kappa -\epsilon, \nu_{\kappa -\epsilon}^{d,m}(\rho)}^{d,l} = 
\prod_{j=1}^d V_j(S_{\kappa_j -\epsilon_j, 
(\nu_{\kappa -\epsilon}^{d,m}(\rho))_j}^{1,l_j,0}) = \\
\prod_{j=1}^d (E -V_j(E -S_{\kappa_j -\epsilon_j, 
\nu_{\kappa_j -\epsilon_j}^{1,m_j}(\rho_j)}^{1,l_j,0}))\\
= \sum_{ \gamma \in \Upsilon^d} (-\e)^\gamma \prod_{j \in \s(\gamma)} 
V_j(E -S_{\kappa_j -\epsilon_j, \nu_{\kappa_j -\epsilon_j}^{1,m_j}(\rho_j)}^{1,l_j,0}). 
\end{multline*}

Подставляя (1.3.22) в (1.3.21), находим, что при $ d \in \N,
l,m,\kappa \in \Z_+^d, \nu \in \Nu_{-m, 2^\kappa -\e}^d $ имеет место равенство
\begin{multline*} \tag{1.3.23}
U_{\kappa, \nu}^{d,l,m} = \\
\sum_{\gamma \in \Upsilon^d: \s(\kappa) \setminus \s(\gamma) =
\emptyset} (-\e)^\gamma \sum_{\epsilon \in \Upsilon^d:
\s(\epsilon) \subset \s(\kappa)} (-\e)^\epsilon \sum_{ \rho \in
\Rho_{\kappa,\nu,\epsilon}^{d,m}} (\prod_{j \in \s(\epsilon)}
a_{\nu_j -2\rho_j}^{m_j}) (\prod_{j \in \s(\gamma)} 
V_j(E -S_{\kappa_j -\epsilon_j, \nu_{\kappa_j -\epsilon_j}^{1,m_j}(\rho_j)}^{1,l_j,0})) +\\
\sum_{\gamma \in \Upsilon^d: \s(\kappa) \setminus \s(\gamma) \ne
\emptyset} (-\e)^\gamma \sum_{\epsilon \in \Upsilon^d:
\s(\epsilon) \subset \s(\kappa)} (-\e)^\epsilon \sum_{ \rho \in
\Rho_{\kappa,\nu,\epsilon}^{d,m}} (\prod_{j \in \s(\epsilon)}
a_{\nu_j -2\rho_j}^{m_j}) (\prod_{j \in \s(\gamma)} 
V_j(E -S_{\kappa_j -\epsilon_j, \nu_{\kappa_j -\epsilon_j}^{1,m_j}(\rho_j)}^{1,l_j,0})). 
\end{multline*}

При $ d \in \N, l,m,\kappa \in \Z_+^d, \nu \in \Nu_{-m, 2^\kappa -\e}^d $ 
для $ \gamma \in \Upsilon^d: \s(\kappa) \setminus \s(\gamma) \ne \emptyset, $ 
соблюдается равенство
\begin{multline*} \tag{1.3.24}
\sum_{\epsilon \in \Upsilon^d: \s(\epsilon) \subset \s(\kappa)}
(-\e)^\epsilon \sum_{ \rho \in \Rho_{\kappa,\nu,\epsilon}^{d,m}}
(\prod_{j \in \s(\epsilon)} a_{\nu_j -2\rho_j}^{m_j}) (\prod_{j \in \s(\gamma)} 
V_j(E -S_{\kappa_j -\epsilon_j, \nu_{\kappa_j -\epsilon_j}^{1,m_j}(\rho_j)}^{1,l_j,0})) =\\
\sum_{\substack{\epsilon \in \Upsilon^d:\\ \s(\epsilon) \subset
(\s(\kappa) \cap \s(\gamma))}} (-\e)^\epsilon
\sum_{\substack{\epsilon^\prime \in \Upsilon^d:\\
\s(\epsilon^\prime) \subset (\s(\kappa) \setminus \s(\gamma))}}
(-\e)^{\epsilon^\prime} \sum_{ \rho \in \Rho_{\kappa,\nu,\epsilon
+\epsilon^\prime}^{d,m}} (\prod_{j \in \s(\epsilon)} a_{\nu_j -2\rho_j}^{m_j}) 
(\prod_{j \in \s(\epsilon^\prime)} a_{\nu_j -2\rho_j}^{m_j})\\ \times 
(\prod_{j \in \s(\gamma)} V_j(E -S_{\kappa_j -\epsilon_j, 
\nu_{\kappa_j -\epsilon_j}^{1,m_j}(\rho_j)}^{1,l_j,0})).
\end{multline*}

Далее, при  $ d \in \N, l,m,\kappa \in \Z_+^d, \nu \in \Nu_{-m, 2^\kappa -\e}^d, 
\gamma \in \Upsilon^d: \s(\kappa) \setminus \s(\gamma) \ne \emptyset, 
\epsilon \in \Upsilon^d: \s(\epsilon) \subset (\s(\kappa) \cap \s(\gamma)), 
\epsilon^\prime \in \Upsilon^d: \s(\epsilon^\prime) \subset (\s(\kappa) 
\setminus \s(\gamma)), $ обозначая для $ j \in \Nu_{1,d}^1 $ и $ \rho_j \in
(\Rho_{\kappa,\nu,\epsilon +\epsilon^\prime}^{d,m})^{\{j\}} = 
\Rho_{\kappa_j,\nu_j,\epsilon_j +\epsilon_j^\prime}^{1,m_j} $
через $ A_{\rho_j}^j $ оператор
$$
A_{\rho_j}^j = \begin{cases} a_{\nu_j -2\rho_j}^{m_j} V_j(E -S_{\kappa_j -1, 
\nu_{\kappa_j -1}^{1,m_j}(\rho_j)}^{1,l_j,0}), j \in \s(\epsilon); \\
V_j(E -S_{\kappa_j, \nu_{\kappa_j}^{1,m_j}(\rho_j)}^{1,l_j,0}), 
j \in \s(\gamma) \setminus \s(\epsilon); \\
a_{\nu_j -2\rho_j}^{m_j} E, j \in \s(\epsilon^\prime); \\
E, j \in (\Nu_{1,d}^1 \setminus \s(\gamma)) \setminus \s(\epsilon^\prime),
\end{cases}
$$
и учитывая (1.2.2), (1.3.5), выводим
\begin{multline*} \tag{1.3.25}
\sum_{ \rho \in \Rho_{\kappa,\nu,\epsilon +\epsilon^\prime}^{d,m}}
(\prod_{j \in \s(\epsilon)} a_{\nu_j -2\rho_j}^{m_j}) (\prod_{j
\in \s(\epsilon^\prime)} a_{\nu_j -2\rho_j}^{m_j}) (\prod_{j \in \s(\gamma)} 
V_j(E -S_{\kappa_j -\epsilon_j, \nu_{\kappa_j -\epsilon_j}^{1,m_j}(\rho_j)}^{1,l_j,0}))\\
= \sum_{ \rho \in \Z^d: \rho_j \in
(\Rho_{\kappa,\nu,\epsilon +\epsilon^\prime}^{d,m})^{\{j\}}, j =1,\ldots,d} 
(\prod_{j=1}^d A_{\rho_j}^j)
= \prod_{j \in \Nu_{1,d}^1} 
(\sum_{ \rho_j \in (\Rho_{\kappa,\nu,\epsilon +\epsilon^\prime}^{d,m})^{\{j\}}} A_{\rho_j}^j)\\
= \biggl(\prod_{j \in \s(\gamma) \setminus \s(\epsilon)} 
V_j(E -S_{\kappa_j, \nu_{\kappa_j}^{1,m_j}(\nu_j)}^{1,l_j,0})\biggr)\\
\times \biggl(\prod_{j \in \s(\epsilon)} \biggl(\sum_{\rho_j \in
\Nu_{-m_j, 2^{\kappa_j -1} -1}^1: (\nu_j -2\rho_j) \in \Nu_{0,m_j +1}^1} 
a_{\nu_j -2\rho_j}^{m_j} V_j(E -S_{\kappa_j -1, \nu_{\kappa_j -1}^{1,m_j}(\rho_j)}^{1,l_j,0})\biggr)\biggr).
\end{multline*}

Подставляя (1.3.25) в (1.3.24), приходим к равенству
\begin{multline*} \tag{1.3.26}
\sum_{\epsilon \in \Upsilon^d: \s(\epsilon) \subset \s(\kappa)}
(-\e)^\epsilon \sum_{ \rho \in \Rho_{\kappa,\nu,\epsilon}^{d,m}}
(\prod_{j \in \s(\epsilon)} a_{\nu_j -2\rho_j}^{m_j}) (\prod_{j
\in \s(\gamma)} V_j(E -S_{\kappa_j -\epsilon_j, 
\nu_{\kappa_j -\epsilon_j}^{1,m_j}(\rho_j)}^{1,l_j,0}))\\
= \sum_{\epsilon \in \Upsilon^d: \s(\epsilon) \subset (\s(\kappa)
\cap \s(\gamma))} (-\e)^\epsilon \sum_{\epsilon^\prime \in
\Upsilon^d: \s(\epsilon^\prime) \subset (\s(\kappa) \setminus \s(\gamma))} 
(-\e)^{\epsilon^\prime} (\prod_{j \in \s(\gamma) \setminus \s(\epsilon)} 
V_j(E -S_{\kappa_j, \nu_{\kappa_j}^{1,m_j}(\nu_j)}^{1,l_j,0}))\\
\times \biggl(\prod_{j \in \s(\epsilon)} \biggl(\sum_{\rho_j \in
\Nu_{-m_j, 2^{\kappa_j -1} -1}^1: (\nu_j -2\rho_j) \in \Nu_{0,m_j +1}^1} 
a_{\nu_j -2\rho_j}^{m_j} V_j(E -S_{\kappa_j -1, \nu_{\kappa_j -1}^{1,m_j}(\rho_j)}^{1,l_j,0})\biggr)\biggr) =0
\end{multline*}
при $ d \in \N, l,m,\kappa \in \Z_+^d, \nu \in \Nu_{-m, 2^\kappa -\e}^d, 
\gamma \in \Upsilon^d: \s(\kappa) \setminus \s(\gamma) \ne \emptyset, $ ибо 
в этом случае
$$
\sum_{\epsilon^\prime \in \Upsilon^d: \s(\epsilon^\prime) \subset
(\s(\kappa) \setminus \s(\gamma))} (-\e)^{\epsilon^\prime} =0.
$$

Соединяя (1.3.26) с (1.3.23), заключаем, что при $ d \in \N, l,m,\kappa \in \Z_+^d, 
\nu \in \Nu_{-m, 2^\kappa -\e}^d $ справедливо равенство
\begin{multline*} \tag{1.3.27}
U_{\kappa,\nu}^{d,l,m} = \sum_{ \gamma \in \Upsilon^d: \s(\kappa)
\subset \s(\gamma)} (-\e)^\gamma \sum_{\epsilon \in \Upsilon^d:
\s(\epsilon) \subset \s(\kappa)} (-\e)^\epsilon \sum_{ \rho \in
\Rho_{\kappa,\nu,\epsilon}^{d,m}}\\
(\prod_{j \in \s(\epsilon)} a_{\nu_j -2\rho_j}^{m_j}) (\prod_{j \in \s(\gamma)} 
V_j(E -S_{\kappa_j -\epsilon_j, \nu_{\kappa_j -\epsilon_j}^{1,m_j}(\rho_j)}^{1,l_j,0})).
\end{multline*}

Лемма 1.3.4

Пусть $ d \in \N, l \in \N^d, m \in \Z_+^d, 1 \le p < \infty, 1 \le q \le \infty $ 
и $ \lambda \in \Z_+^d(m). $ Тогда существуют
константы $ c_{11}(d,l,m,p,q,\lambda) >0 $ и $ c_{12}(d,m) >0 $
такие, что для любой функции $ f \in L_p(I^d) $ при $ \kappa \in
\Z_+^d $ выполняется неравенство
\begin{multline*} \tag{1.3.28}
\| \D^\lambda \mathcal E_\kappa^{d,l -\e,m} f \|_{L_q(I^d)} \le
c_{11} (\prod_{j \in \s(\kappa)} 2^{\kappa_j (\lambda_j +p^{-1}
+(p^{-1} -q^{-1})_+)})\\
\times \biggl(\int_{ (c_{12} 2^{-\kappa} B^d)^{\s(\kappa)}} 
\int_{(I^d)_\xi^{l \chi_{\s(\kappa)}}} |\Delta_\xi^{l \chi_{\s(\kappa)}}
f(x)|^p dx d\xi^{\s(\kappa)} \biggr)^{1/p}. 
\end{multline*}

Доказательство.

Сначала, в дополнение к уже имеющимся,  построим некоторые вспомогательные 
объекты и приведём их свойства, которые нам понадобятся при 
доказательстве.

При $ d \in \N, m,\kappa \in \Z_+^d, \nu \in \Nu_{0, 2^\kappa -\e}^d $ зададим 
точку $ x_{\kappa, \nu}^{\prime d,m} \in \R^d $ и вектор 
$ \delta_{\kappa, \nu}^{\prime d,m} \in \R_+^d, $ полагая
\begin{multline*}
( x_{\kappa, \nu}^{\prime d,m})_j = 2^{-\kappa_j} \min((\nu_j -2m_j -1)_+, 
(2^{\kappa_j} -2m_j -3)_+),
 (\delta_{\kappa, \nu}^{\prime d,m})_j \\
= 2^{-\kappa_j} \min(2^{\kappa_j}, 2m_j +3), j=1,\ldots,d, \ 
\end{multline*} 
и определим клетку $ D_{\kappa, \nu}^{\prime d,m} $
равенством $ D_{\kappa, \nu}^{\prime d,m} = x_{\kappa, \nu}^{\prime d,m} + 
\delta_{\kappa, \nu}^{\prime d,m} I^d.$

Из определения следует, что при $ d \in \N, m,\kappa \in \Z_+^d,
\nu \in \Nu_{0, 2^\kappa -\e}^d $ множество
\begin{equation*} \tag{1.3.29}
 D_{\kappa, \nu}^{\prime d,m} \subset I^d, \ 
\end{equation*}
а также
\begin{equation*} \tag{1.3.30}
Q_{\kappa,\nu}^d = D_{\kappa,\nu}^{d,0} \subset
D_{\kappa, \nu}^{\prime d,m}.
\end{equation*}

С помощью (1.3.30) легко проверить, что при $ d \in \N, m \in \Z_+^d $ 
существует константа $ c_{13}(d,m) >0 $ такая, что для $ \kappa \in \Z_+^d $ и 
$ x \in I^d $ верно неравенство
\begin{equation*} \tag{1.3.31}
\card \{ \nu \in \Nu_{0, 2^\kappa -\e}^d: x \in D_{\kappa, \nu}^{\prime d,m} \} 
\le c_{13}.
\end{equation*}

Пользуясь тем, что для $ t \in \R, a \in \R_+ $ справедливы
соотношения $ (at)_+ = at_+ $ и $ t_+ \le (t+a)_+ \le t_+ +a, $ а
также используя (1.3.3), несложно показать, что при $ d \in \N,
m,\kappa \in \Z_+^d, $ для $ \nu^\prime \in \Nu_{0, 2^\kappa -\e}^d, 
\nu \in \Nu_{-m, 2^\kappa -\e}^d: Q_{\kappa, \nu^\prime}^d \cap 
\supp g_{\kappa, \nu}^{d,m} \ne \emptyset, \epsilon \in \Upsilon^d: 
\s(\epsilon) \subset \s(\kappa), \rho \in \Rho_{\kappa, \nu, \epsilon}^{d,m} $ 
имеет место включение
\begin{equation*} \tag{1.3.32}
D_{\kappa -\epsilon, \nu_{\kappa -\epsilon}^{d,m}(\rho)}^{d,0} \subset 
D_{\kappa, \nu^\prime}^{\prime d,m}.
\end{equation*}

Отметим ещё, что
\begin{equation*} \tag{1.3.33}
\card \Rho_{\kappa,\nu,\epsilon}^{d,m} \le c_{14}(d,m), \ 
\end{equation*}
$ d \in \N, m,\kappa \in \Z_+^d, \nu \in \Nu_{-m, 2^\kappa -\e}^d,
\epsilon \in \Upsilon^d: \s(\epsilon) \subset \s(\kappa). $

Определим ещё при $ l, m, \kappa \in \Z_+, \nu \in \Nu_{0, 2^\kappa -1}^1 $ 
непрерывный линейный оператор $ S_{\kappa, \nu}^{\prime l, m}: 
L_1(I) \mapsto \mathcal P^{1, l}(I) \cap L_\infty (I), $ полагая для 
$ f \in L_1(I) $ значение
$$
(S_{\kappa, \nu}^{\prime l, m} f)(x) = 
(P_{\delta, x^0}^{1, l}(f \mid_{(x^0 +\delta I)}))(x), x \in I,
$$
при $ \delta = 2^{-\kappa} \min(2^\kappa, 2 m +3), 
x^0 = 2^{-\kappa} \min((\nu -2 m -1)_+, (2^\kappa -2 m -3)_+). $

Теперь перейдём к выводу оценки (1.3.28).

В условиях леммы 1.3.4 ввиду (1.3.20), (1.3.17) имеем
\begin{multline*} \tag{1.3.34}
\| \D^\lambda \mathcal E_\kappa^{d, l-\e,m}f \|_{L_q(I^d)} \le\\ 
\sum_{\mu \in \Z_+^d(\lambda)} C_\lambda^\mu \| \sum_{ \nu \in
\Nu_{-m, 2^\kappa -\e}^d} \D^\mu (U_{\kappa, \nu}^{d, l-\e,m}f)
\D^{\lambda -\mu} g_{\kappa,\nu}^{d,m} \|_{L_q(I^d)}.
\end{multline*}

Оценивая правую часть (1.3.34), при $ \mu \in \Z_+^d(\lambda) $
находим, что
\begin{multline*} \tag{1.3.35}
\biggl\| \sum_{ \nu \in \Nu_{-m, 2^\kappa -\e}^d} \D^\mu
(U_{\kappa, \nu}^{d, l-\e,m}f) \D^{\lambda -\mu}
g_{\kappa,\nu}^{d,m} \biggr\|_{L_q(I^d)}^q \le\\ 
\sum_{\nu^\prime \in \Nu_{0, 2^\kappa -\e}^d} 
\biggl( \sum_{ \nu \in \Nu_{-m, 2^\kappa -\e}^d: \supp g_{\kappa, \nu}^{d,m} 
\cap Q_{\kappa, \nu^\prime}^d \ne \emptyset} \|\D^\mu (U_{\kappa, \nu}^{d,
l-\e,m}f) \D^{\lambda -\mu} g_{\kappa,\nu}^{d,m} 
\|_{L_q(Q_{\kappa, \nu^\prime}^d)}\biggr)^q.
\end{multline*}

Используя (1.3.6), (1.1.2) и (1.3.30), выводим
\begin{multline*} \tag{1.3.36}
\|\D^\mu (U_{\kappa, \nu}^{d, l-\e,m}f) \D^{\lambda -\mu}
g_{\kappa,\nu}^{d,m} \|_{L_q( Q_{\kappa, \nu^\prime}^d)} \le
c_{15} 2^{(\kappa, \lambda +p^{-1} \e -q^{-1} \e)} 
\| U_{\kappa, \nu}^{d, l-\e,m}f \|_{L_p( D_{\kappa, \nu^\prime}^{\prime d,m})} \ 
\end{multline*}
при $ \nu^\prime \in \Nu_{0, 2^\kappa -\e}^d, \nu \in \Nu_{-m, 2^\kappa -\e}^d: 
\supp g_{\kappa, \nu}^{d,m} \cap Q_{\kappa, \nu^\prime}^d \ne \emptyset. $

Из (1.3.27) для $ \nu^\prime \in \Nu_{0, 2^\kappa -\e}^d, \nu \in
\Nu_{-m, 2^\kappa -\e}^d: \supp g_{\kappa, \nu}^{d,m} \cap
Q_{\kappa, \nu^\prime}^d \ne \emptyset, $ имеем неравенство
\begin{multline*} \tag{1.3.37}
\| U_{\kappa, \nu}^{d, l-\e,m}f \|_{L_p( D_{\kappa, \nu^\prime}^{\prime d,m})} \le\\ 
\sum_{ \gamma \in \Upsilon^d: \s(\kappa) \subset \s(\gamma)}
\sum_{\epsilon \in \Upsilon^d: \s(\epsilon) \subset \s(\kappa)}
\sum_{ \rho \in \Rho_{\kappa,\nu,\epsilon}^{d,m}} (\prod_{j \in
\s(\epsilon)} a_{\nu_j -2\rho_j}^{m_j}) \|(\prod_{j \in \s(\gamma)} 
V_j(E -S_{\kappa_j -\epsilon_j, \nu_{\kappa_j -\epsilon_j}^{1,m_j}(\rho_j)}^{1,l_j
-1,0}))f \|_{L_p( D_{\kappa, \nu^\prime}^{\prime d,m})}.
\end{multline*}

Далее, для $ \nu^\prime \in \Nu_{0, 2^\kappa -\e}^d, \nu \in
\Nu_{-m, 2^\kappa -\e}^d: \supp g_{\kappa, \nu}^{d,m} \cap
Q_{\kappa, \nu^\prime}^d \ne \emptyset, \gamma \in \Upsilon^d:
\s(\kappa) \subset \s(\gamma), \epsilon \in \Upsilon^d:
\s(\epsilon) \subset \s(\kappa), \rho \in
\Rho_{\kappa,\nu,\epsilon}^{d,m}, $ используя (1.2.2), а затем с
учётом (1.3.32), (1.3.29) применяя (1.2.6), получаем
\begin{multline*} \tag{1.3.38}
\|(\prod_{j \in \s(\gamma)} V_j(E -S_{\kappa_j -\epsilon_j,
\nu_{\kappa_j -\epsilon_j}^{1,m_j}(\rho_j)}^{1,l_j -1,0}))f 
\|_{L_p( D_{\kappa, \nu^\prime}^{\prime d,m})} \le\\ 
c_{16} \|(\prod_{j \in \s(\kappa)}
V_j(E -S_{\kappa_j -\epsilon_j, \nu_{\kappa_j -\epsilon_j}^{1,m_j}(\rho_j)}^{1,l_j -1,0}))f
\|_{L_p( D_{\kappa, \nu^\prime}^{\prime d,m})} \  
\end{multline*} 
и, пользуясь тем, что при $ j =1,\ldots,d $ в виду (1.1.3) соблюдено
равенство
$$ 
(E -S_{\kappa_j -\epsilon_j, \nu_{\kappa_j -\epsilon_j}^{1,m_j}(\rho_j)}^{1,l_j -1,0}) =
(E -S_{\kappa_j -\epsilon_j, \nu_{\kappa_j -\epsilon_j}^{1,m_j}(\rho_j)}^{1,l_j -1,0}) 
(E -S_{\kappa_j, \nu_j^\prime}^{\prime l_j -1, m_j}),
$$
благодаря п.2 леммы 1.2.1 и (1.2.2), а также в силу (1.3.32),
(1.3.29), (1.2.6) выводим
\begin{multline*} \tag{1.3.39}
\|(\prod_{j \in \s(\kappa)} V_j(E -S_{\kappa_j -\epsilon_j,
\nu_{\kappa_j -\epsilon_j}^{1,m_j}(\rho_j)}^{1,l_j -1,0}))f 
\|_{L_p( D_{\kappa, \nu^\prime}^{\prime d,m})} \le\\ 
c_{17} \| (\prod_{j \in \s(\kappa)} V_j(E -S_{\kappa_j,
\nu_j^\prime}^{\prime l_j -1, m_j}))f \|_{L_p( D_{\kappa, \nu^\prime}^{\prime d,m})}. 
\end{multline*}

Объединяя (1.3.37), (1.3.38), (1.3.39) и учитывая (1.3.33), на
основании (1.2.5) (при $ J = \emptyset, J^\prime = \s(\kappa)$)
находим, что для $ \nu^\prime \in \Nu_{0, 2^\kappa -\e}^d, \nu \in
\Nu_{-m, 2^\kappa -\e}^d: \supp g_{\kappa, \nu}^{d,m} \cap
Q_{\kappa, \nu^\prime}^d \ne \emptyset, $ выполняется неравенство
\begin{multline*} \tag{1.3.40}
\| U_{\kappa, \nu}^{d, l-\e,m}f \|_{L_p( D_{\kappa, \nu^\prime}^{\prime d,m})} \le\\ 
c_{18} (\prod_{j \in \s(\kappa)} 2^{\kappa_j /p}) 
\biggl(\int_{ (c_{12} 2^{-\kappa} B^d)^{\s(\kappa)}} 
\int_{ (D_{\kappa, \nu^\prime}^{\prime d,m})_\xi^{l \chi_{\s(\kappa)}}} 
|\Delta_\xi^{l \chi_{\s(\kappa)}} f(x)|^p dx d\xi^{\s(\kappa)}\biggr)^{1/p}. 
\end{multline*}

Из (1.3.36) и (1.3.40) следует,что
\begin{multline*} \tag{1.3.41}
\|\D^\mu (U_{\kappa, \nu}^{d, l-\e,m}f) \D^{\lambda -\mu}
g_{\kappa,\nu}^{d,m} \|_{L_q( Q_{\kappa, \nu^\prime}^d)} \le\\ 
c_{19} 2^{(\kappa, \lambda +p^{-1} \e -q^{-1} \e)} (\prod_{j \in \s(\kappa)} 
2^{\kappa_j /p}) \biggl(\int_{ (c_{12} 2^{-\kappa} B^d)^{\s(\kappa)} } 
\int_{ (D_{\kappa, \nu^\prime}^{\prime d,m})_\xi^{l \chi_{\s(\kappa)}}} 
|\Delta_\xi^{l \chi_{\s(\kappa)}} f(x)|^p dx d\xi^{\s(\kappa)}\biggr)^{1/p}, \ 
\end{multline*}
при $ \nu^\prime \in \Nu_{0, 2^\kappa -\e}^d, \nu \in \Nu_{-m, 2^\kappa -\e}^d: 
\supp g_{\kappa, \nu}^{d,m} \cap Q_{\kappa, \nu^\prime}^d \ne \emptyset, 
\mu \in \Z_+^d(\lambda). $
Подставляя (1.3.41) в (1.3.35), и учитывая (1.3.12), а затем
применяя неравенство (1.1.1) при $ a =  p/q \le 1 $ и
используя (1.3.31), получаем, что при $ p \le q $ и $ \mu \in
\Z_+^d(\lambda) $ справедливо неравенство
\begin{multline*} \tag{1.3.42}
\| \sum_{ \nu \in \Nu_{-m, 2^\kappa -\e}^d} \D^\mu (U_{\kappa, \nu}^{d, l-\e,m}f) 
\D^{\lambda -\mu} g_{\kappa,\nu}^{d,m} \|_{L_q(I^d)}^q \le\\ 
(c_{20} 2^{(\kappa, \lambda
+2p^{-1} \e -q^{-1} \e)})^q \biggl(\sum_{ \nu^\prime \in \Nu_{0, 2^\kappa -\e}^d} 
\int_{ (c_{12} 2^{-\kappa} B^d)^{\s(\kappa)}}
\int_{ D_{\kappa, \nu^\prime}^{\prime d,m} \cap (I^d)_\xi^{l \chi_{\s(\kappa)}}} 
|\Delta_\xi^{l \chi_{\s(\kappa)}} f(x)|^p dx d\xi^{\s(\kappa)}\biggr)^{q/p} \le\\
\biggl(c_{21} 2^{(\kappa, \lambda +2p^{-1} \e -q^{-1} \e)}
\biggl(\int_{ (c_{12} 2^{-\kappa} B^d)^{\s(\kappa)}} \int_{(I^d)_\xi^{l \chi_{\s(\kappa)}}} 
|\Delta_\xi^{l \chi_{\s(\kappa)}}
f(x)|^p dx d\xi^{\s(\kappa)}\biggr)^{1/p}\biggr)^q.
\end{multline*}

Соединяя (1.3.42) с (1.3.34), приходим к (1.3.28) при $ p \le q. $
Для получения (1.3.28) при $ q < p $ достаточно заметить, что в этом случае
$$
\| \D^\lambda \mathcal E_\kappa^{d,l -\e,m}f \|_{L_q(I^d)} \le 
\| \D^\lambda \mathcal E_\kappa^{d,l -\e,m}f \|_{L_p(I^d)}, \ 
$$
и применить (1.3.28) при $ q = p. \square $

Следствие

В условиях леммы 1.3.4 существует константа $ c_{22}(d,l,m,\lambda,p,q) > 0 $ 
такая, что для $ f \in L_p(I^d) $ при $ \kappa \in \Z_+^d $ имеет место 
неравенство
\begin{equation*} \tag{1.3.43}
\| \D^\lambda \mathcal E_\kappa^{d,l-\e,m}f \|_{L_q(I^d)} \le
c_{22} 2^{(\kappa, \lambda +(p^{-1} -q^{-1})_+ \e)} 
\Omega^{\prime l \chi_{\s(\kappa)}}(f, (c_{12} 2^{-\kappa})^{\s(\kappa)})_{L_p(I^d)}.
\end{equation*}

Доказательство.

В самом деле, для $ f \in L_p(I^d) $ при $ \kappa \in \Z_+^d $ в силу (1.3.28) 
имеем
\begin{multline*}
\| \D^\lambda \mathcal E_\kappa^{d,l -\e,m}f \|_{L_q(I^d)} \le
c_{11} (\prod_{j \in \s(\kappa)} 2^{\kappa_j (\lambda_j +(p^{-1} -q^{-1})_+)}) \\
\times \biggl((\prod_{j \in \s(\kappa)} 2^{\kappa_j}) \int_{ (c_{12} 2^{-\kappa} B^d)^{\s(\kappa)}} 
\int_{(I^d)_\xi^{l \chi_{\s(\kappa)}}} |\Delta_\xi^{l \chi_{\s(\kappa)}}
f(x)|^p dx d\xi^{\s(\kappa)} \biggr)^{1/p} \le \\
c_{22} 2^{(\kappa, \lambda +(p^{-1} -q^{-1})_+ \e)} 
\Omega^{\prime l \chi_{\s(\kappa)}}(f, (c_{12} 2^{-\kappa})^{\s(\kappa)})_{L_p(I^d)}. \square
\end{multline*}

Теорема 1.3.5

Пусть $ d \in \N, l \in \N^d, 1 \le p < \infty, 1 \le q \le \infty, 
m \in \Z_+^d, \lambda \in \Z_+^d(m). $ Тогда если для функции $ f \in L_p(I^d) $ 
и любого непустого множества $ J \subset \Nu_{1,d}^1 $ функция
\begin{equation*} \tag{1.3.44}
\biggl(\prod_{j \in J} t_j^{-\lambda_j -p^{-1} -(p^{-1} -q^{-1})_+ -1}\biggr) 
\biggl(\int_{(c_{12} t B^d)^J} \int_{ (I^d)_\xi^{l \chi_J}} 
|\Delta_\xi^{l \chi_J} f(x)|^p dx d\xi^J\biggr)^{1/p} \in L_1((I^d)^J), \ 
\end{equation*}
то в $ L_q(I^d) $ имеет место равенство
\begin{equation*} \tag{1.3.45}
\D^\lambda f = \sum_{\kappa \in \Z_+^d} \D^\lambda \mathcal E_\kappa^{d, l -\e,m} f.
\end{equation*}

Доказательство.

В условиях теоремы 1.3.5, прежде всего, заметим, что в силу (1.3.8)
справедливо соотношение
\begin{equation*}
\|f -E_\kappa^{d,l -\e,m}f \|_{L_p(I^d)} \le
c_{23} \sum_{j=1}^d \Omega^{l e_j} (f, c_{3} 2^{-\kappa_j} )_{L_p(I^d)} \to 0 
\text{ при } \mn(\kappa) \to \infty. 
\end{equation*}

Отсюда на основании леммы 1.3.1 заключаем, что в $ L_p(I^d) $ имеет 
место равенство
\begin{equation*} \tag{1.3.46}
f = \sum_{ \kappa \in \Z_+^d} \mathcal E_\kappa^{d, l -\e,m}f.
\end{equation*}

Далее, для функции $ f \in L_p(I^d), $ подчинённой условию (1.3.44), и 
каждого множества $ J \subset \Nu_{1,d}^1: J \ne \emptyset, $ при 
$ \kappa \in \Z_+^d: \s(\kappa) = J, $ ввиду (1.3.28) выполняется неравенство
\begin{multline*} \tag{1.3.47}
\| \D^\lambda \mathcal E_\kappa^{d, l -\e,m}f \|_{L_q(I^d)} \le \\
c_{11} \biggl(\prod_{j \in J} 2^{\kappa_j (\lambda_j +p^{-1}
+(p^{-1} -q^{-1})_+)}\biggr) \biggl(\int_{ (c_{12} 2^{-\kappa}
B^d)^J} \int_{ (I^d)_\xi^{l \chi_J}} |\Delta_\xi^{l \chi_J}
f(x)|^p dx d\xi^J\biggr)^{1/p} \le\\
c_{24} \int_{ (2^{-\kappa} +2^{-\kappa} I^d)^J}
\biggl(\prod_{j \in J} t_j^{-\lambda_j -p^{-1} -(p^{-1} -q^{-1})_+
-1}\biggr) \biggl(\int_{ (c_{12} t B^d)^J} \int_{ (I^d)_\xi^{l \chi_J}} 
|\Delta_\xi^{l \chi_J} f(x)|^p dx d\xi^J\biggr)^{1/p} dt^J.
\end{multline*}

Поскольку множества $ (2^{-\kappa} +2^{-\kappa} I^d)^J, \kappa \in
\Z_+^d: \s(\kappa) = J, $ попарно не пересекаются и $ \cup_{ \kappa
\in \Z_+^d: \s(\kappa) = J} (2^{-\kappa} +2^{-\kappa} I^d)^J
\subset (I^d)^J, $ то из (1.3.44) и (1.3.47) вытекает, что ряд
$$
\sum_{ \kappa \in \Z_+^d: \s(\kappa) = J} \| \D^\lambda \mathcal E_\kappa^{d,
l -\e,m}f \|_{L_q(I^d)}
$$
сходится для $ J \subset \Nu_{1,d}^1: J \ne \emptyset, $ \ 
а, следовательно, и ряд $ \sum_{ \kappa \in \Z_+^d} 
\| \D^\lambda \mathcal E_\kappa^{d, l -\e,m}f \|_{L_q(I^d)} $ сходится, и, 
значит, ряд
$$
\sum_{ \kappa \in \Z_+^d} \D^\lambda \mathcal E_\kappa^{d, l -\e,m}f
$$
сходится в $ L_q(I^d). $

Принимая во внимание это обстоятельство и равенство (1.3.46), для
любой функции $ \phi \in C_0^\infty (I^d) $ имеем
\begin{multline*}
\langle \D^\lambda f, \phi \rangle
= \int_{I^d} (\sum_{ \kappa \in \Z_+^d} \mathcal E_\kappa^{d, l -\e,m}f) 
(-1)^{(\lambda,\e)} \D^\lambda \phi dx
= \int_{I^d} (\sum_{ \kappa \in \Z_+^d} 
\D^\lambda (\mathcal E_\kappa^{d, l -\e,m}f)) \phi dx.
\end{multline*}

А это значит, что в $ L_q(I^d) $ верно равенство (1.3.45). $ \square $

Предложение 1.3.6

Пусть $ d \in \N, \alpha \in \R_+^d, l = l(\alpha), 1 \le p < \infty, 
1 \le q \le \infty, m \in \Z_+^d, \lambda \in \Z_+^d(m) $
и
\begin{equation*} 
\alpha -\lambda -(p^{-1} -q^{-1})_+ \e >0.
\end{equation*}
Тогда для любой функции $ f \in (S_p^\alpha H)^\prime(I^d) $ в $ L_q(I^d) $
имеет место равенство (1.3.45).

Доказательство предложения 1.3.6 почти дословно повторяет доказательство 
соответствующего утверждения из [16].
\bigskip

\centerline{\S 2. Продолжения функций из $ (S_{p,\theta}^\alpha B)^\prime(I^d) $ 
на всё пространство $ \R^d $}
\centerline{и их дифференциально-разностные свойства}
\bigskip

2.1. В этом пункте рассматриваются подпространства кусочно-полиномиальных 
функций, обладающих свойствами полезными для получения основных результатов 
работы.  

Введём в рассмотрение следующие пространства кусочно-полиномиальных 
функций.

При $ d \in \N, l \in \Z_+^d, m \in \N^d $ и $ \kappa \in \Z_+^d $ через 
$ \mathcal P_\kappa^{d,l,m} $ обозначим линейное пространство, состоящее из 
функций $ f: \R^d \mapsto \R, $ для каждой из которых существует набор  
полиномов $ \{f_\nu \in \mathcal P^{d,l}, \nu \in \Nu_{-m, 2^\kappa -\e}^d\} $ 
такой, что для $ x \in \R^d $ выполняется равенство
\begin{equation*} \tag{2.1.1}
f(x) = \sum_{\nu \in \Nu_{-m, 2^\kappa -\e}^d} f_\nu(x) g_{\kappa,\nu}^{d,m}(x).
\end{equation*}

При $ d \in \N, m \in \N^d $ обозначим через $ Q^{d,m} $ клетку
$ Q^{d, m} = -(m +\e) +(2m +3 \e) I^d. $ Отметим, что при 
$ d \in \N, l \in \Z_+^d, m \in \N^d $ и $ \kappa \in \Z_+^d $ для 
$ f \in \mathcal P_\kappa^{d,l,m} $ её носитель $ \supp f \subset Q^{d,m}, $
ибо $ \supp g_{\kappa,\nu}^{d,m} \subset Q^{d,m}, \nu \in \Nu_{-m, 2^\kappa -\e}^d. $ 

При $ d \in \N $ для $ j =1, \ldots, d $ зададим множество $ J_j, $ полагая
$$
J_j = \{1, \ldots, j -1, j +1, \ldots, d\}.
$$
При $ d \in \N, l \in \Z_+^d, m \in \N^d, \kappa \in \Z_+^d, $ пользуясь 
тем, что полиномы $ f_\nu \in \mathcal P^{d,l}, \nu \in \Nu_{-m, 2^\kappa -\e}^d, $ 
представляются в виде
$$
f_\nu(x) = \sum_{\lambda \in \Z_+^d(l)} a_\lambda^\nu x^\lambda, \ x \in \R^d, 
a_\lambda^\nu \in \R, \lambda \in \Z_+^d(l), \nu \in \Nu_{-m, 2^\kappa -\e}^d,
$$
для функции $ f, $ задаваемой равенством (2.1.1), при $ j =1, \ldots, d $ и
$ x \in \R^d $ имеем
\begin{multline*} \tag{2.1.2}
f(x) = \sum_{\nu \in \Nu_{-m, 2^\kappa -\e}^d} (\sum_{\lambda \in \Z_+^d(l)}
a_\lambda^\nu x^\lambda) g_{\kappa, \nu}^{d,m}(x) = \\
\sum_{\nu_j \in 
\Nu_{-m_j, 2^{\kappa_j} -1}^1, \nu^{J_j} \in (\Nu_{-m, 2^\kappa -\e}^d)^{J_j}} 
(\sum_{\lambda_j =0, \ldots, l_j, \lambda^{J_j} \in (\Z_+^d(l))^{J_j}}
a_\lambda^\nu x_j^{\lambda_j} (x^{J_j})^{\lambda^{J_j}}) g_{\kappa_j,
\nu_j}^{1, m_j}(x_j) (\prod_{i \in J_j} g_{\kappa_i, \nu_i}^{1, m_i}(x_i)) = \\
\sum_{\nu_j \in \Nu_{-m_j, 2^{\kappa_j} -1}^1, \nu^{J_j} \in 
(\Nu_{-m, 2^\kappa -\e}^d)^{J_j}} 
(\sum_{\lambda_j =0, \ldots, l_j} \sum_{ \lambda^{J_j} \in (\Z_+^d(l))^{J_j}}
a_\lambda^\nu (x^{J_j})^{\lambda^{J_j}} x_j^{\lambda_j})
(\prod_{i \in J_j} g_{\kappa_i, \nu_i}^{1, m_i}(x_i)) 
g_{\kappa_j, \nu_j}^{1, m_j}(x_j) = \\
\sum_{\nu_j \in \Nu_{-m_j, 2^{\kappa_j} -1}^1, \nu^{J_j} \in 
(\Nu_{-m, 2^\kappa -\e}^d)^{J_j}} 
(\sum_{\lambda_j =0, \ldots, l_j} x_j^{\lambda_j} (\sum_{ \lambda^{J_j} \in 
(\Z_+^d(l))^{J_j}} a_\lambda^\nu (x^{J_j})^{\lambda^{J_j}}))
(\prod_{i \in J_j} g_{\kappa_i, \nu_i}^{1, m_i}(x_i)) 
g_{\kappa_j, \nu_j}^{1, m_j}(x_j) = \\
\sum_{\nu_j \in \Nu_{-m_j, 2^{\kappa_j} -1}^1} \sum_{\nu^{J_j} \in 
(\Nu_{-m, 2^\kappa -\e}^d)^{J_j}} 
(\sum_{\lambda_j =0, \ldots, l_j} x_j^{\lambda_j} (\sum_{ \lambda^{J_j} \in 
(\Z_+^d(l))^{J_j}} a_\lambda^\nu (x^{J_j})^{\lambda^{J_j}}))
(\prod_{i \in J_j} g_{\kappa_i, \nu_i}^{1, m_i}(x_i)) 
g_{\kappa_j, \nu_j}^{1, m_j}(x_j) = \\
\sum_{\nu_j \in \Nu_{-m_j, 2^{\kappa_j} -1}^1} (\sum_{\nu^{J_j} \in 
(\Nu_{-m, 2^\kappa -\e}^d)^{J_j}} 
(\sum_{\lambda_j =0, \ldots, l_j} x_j^{\lambda_j} (\sum_{ \lambda^{J_j} \in 
(\Z_+^d(l))^{J_j}} a_\lambda^\nu (x^{J_j})^{\lambda^{J_j}}))
(\prod_{i \in J_j} g_{\kappa_i, \nu_i}^{1, m_i}(x_i))) 
g_{\kappa_j, \nu_j}^{1, m_j}(x_j) = \\
\sum_{\nu_j \in \Nu_{-m_j, 2^{\kappa_j} -1}^1} (\sum_{\nu^{J_j} \in 
(\Nu_{-m, 2^\kappa -\e}^d)^{J_j}} 
\sum_{\lambda_j =0, \ldots, l_j} x_j^{\lambda_j} (\sum_{ \lambda^{J_j} \in 
(\Z_+^d(l))^{J_j}} a_\lambda^\nu (x^{J_j})^{\lambda^{J_j}})
(\prod_{i \in J_j} g_{\kappa_i, \nu_i}^{1, m_i}(x_i))) 
g_{\kappa_j, \nu_j}^{1, m_j}(x_j) = \\
\sum_{\nu_j \in \Nu_{-m_j, 2^{\kappa_j} -1}^1} (\sum_{\lambda_j =0, \ldots, l_j} 
\sum_{\nu^{J_j} \in (\Nu_{-m, 2^\kappa -\e}^d)^{J_j}} x_j^{\lambda_j} 
(\sum_{ \lambda^{J_j} \in (\Z_+^d(l))^{J_j}} a_\lambda^\nu 
(x^{J_j})^{\lambda^{J_j}}) (\prod_{i \in J_j} g_{\kappa_i, \nu_i}^{1, m_i}(x_i))) 
g_{\kappa_j, \nu_j}^{1, m_j}(x_j) = \\
\sum_{\nu_j \in \Nu_{-m_j, 2^{\kappa_j} -1}^1} (\sum_{\lambda_j =0, \ldots, l_j} 
x_j^{\lambda_j} (\sum_{\nu^{J_j} \in (\Nu_{-m, 2^\kappa -\e}^d)^{J_j}} 
(\sum_{ \lambda^{J_j} \in (\Z_+^d(l))^{J_j}} a_\lambda^\nu 
(x^{J_j})^{\lambda^{J_j}}) (\prod_{i \in J_j} g_{\kappa_i, \nu_i}^{1, m_i}(x_i)))) 
g_{\kappa_j, \nu_j}^{1, m_j}(x_j) = \\
\sum_{\nu_j \in \Nu_{-m_j, 2^{\kappa_j} -1}^1} f_{\nu_j}^j(x)
g_{\kappa_j, \nu_j}^{1, m_j}(x_j), \ 
\end{multline*}
где 
\begin{multline*} \tag{2.1.3}
f_{\nu_j}^j(x) = \sum_{\lambda_j =0, \ldots, l_j} x_j^{\lambda_j} 
(\sum_{\nu^{J_j} \in (\Nu_{-m, 2^\kappa -\e}^d)^{J_j}} 
(\sum_{ \lambda^{J_j} \in (\Z_+^d(l))^{J_j}} a_\lambda^\nu  (x^{J_j})^{\lambda^{J_j}}) 
(\prod_{i \in J_j} g_{\kappa_i, \nu_i}^{1, m_i}(x_i))) = \\,
\sum_{\lambda_j =0, \ldots, l_j} x_j^{\lambda_j} 
f_{\nu_j}^{j, \lambda_j}(x^{J_j}), x \in \R^d,
\end{multline*}
причём
\begin{multline*} \tag{2.1.4}
f_{\nu_j}^{j, \lambda_j}(x^{J_j}) = \\
\sum_{\nu^{J_j} \in (\Nu_{-m, 2^\kappa -\e}^d)^{J_j}} 
(\sum_{ \lambda^{J_j} \in (\Z_+^d(l))^{J_j}} a_\lambda^\nu  (x^{J_j})^{\lambda^{J_j}}) 
(\prod_{i \in J_j} g_{\kappa_i, \nu_i}^{1, m_i}(x_i)) \in 
\mathcal P_{\kappa^{J_j}}^{d -1, l^{J_j}, m^{J_j}}, \\ \lambda_j =0, \ldots, l_j, 
\nu_j \in \Nu_{-m_j, 2^{\kappa_j} -1}^1.
\end{multline*} 

Несложно проверить, что при $ d \in \N, l \in \Z_+^d, m \in \N^d, 
\kappa \in \Z_+^d $ отображение, которое каждому набору полиномов
$ \{f_\nu \in \mathcal P^{d, l}, \nu \in \Nu_{-m, 2^\kappa -\e}^d\} $ ставит в 
соответствие функцию $ f, $ задаваемую равенством (2.1.1), является 
изоморфизмом прямого произведения $ (2^\kappa +m)^\e $ экземпляров пространства 
$ \mathcal P^{d,l} $ на пространство $ \mathcal P_\kappa^{d,l,m}. $ Линейность и 
сюръективность этого отображения очевидны, а инъективность легко 
установить при $ d =1, $ а при произвольном $ d \in \N: d \ge 2, $ без
труда выводится по индукции относительно $ d $ с помощью (2.1.2), (2.1.3), (2.1.4).

При $ d \in \N, l \in \Z_+^d, m \in \N^d, \kappa \in \Z_+^d $ обозначим через 
$ \mathcal P_\kappa^{\prime d,l,m} $ подпространство в пространстве 
непрерывных ограниченных в $ \R^d $ функций, 
состоящее из тех функций $ f \in \mathcal P_\kappa^{d,l,m}, $ для которых
при $ j =1, \ldots, d $ соблюдаются равенства 
\begin{equation*} \tag{2.1.5}
f_{\nu_j}^j(x) = f_0^j(x), \ x \in \R^d,  \nu_j \in \Nu_{-m_j, 0}^1; \\ 
\end{equation*}
и
\begin{equation*} \tag{2.1.6}
f_{\nu_j}^j(x) = f_{2^{\kappa_j} -m_j -1}^j(x), \ x \in \R^d, \nu_j \in 
\Nu_{2^{\kappa_j} -m_j -1, 2^{\kappa_j} -1}^1, \ (\text{см. } (2.1.2), (2.1.3)).
\end{equation*}

Замечание

Если для рассматриваемых функций при некотором $ j =1, \ldots, d $ условие 
(2.1.5) или (2.1.6) соблюдается при $ x \in Q^{d,m}, $ то оно выполняется 
при $ x \in \R^d. $

В самом деле, пусть при некотором $ j =1,\ldots,d $ условие (2.1.5) имеет 
место при $ x \in Q^{d,m}, $ т.е. при $ x = (x_1, \ldots, x_{j -1}, x_j, 
x_{j +1}, \ldots, x_d): x^{J_j} \in (Q^{d,m})^{J_j} = Q^{d -1,m^{J_j}} $ и 
$ x_j \in Q^{1,m_j}. $ 
Тогда, поскольку в силу (2.1.3) для каждого $ x^{J_j} \in \R^{d -1} $ верноо 
включение $ f_{\nu_j}^j(x_1, \ldots, x_{j -1}, \cdot, x_{j +1}, \ldots, x_d) 
\in \mathcal P^{1,l_j}, \nu_j \in \Nu_{-m_j, 2^{\kappa_j} -1}^1, $ то (2.1.5)
выполняется при $ x^{J_j} \in Q^{d -1,m^{J_j}} $ и $ x_j \in \R. $
Отсюда, пользуясь тем, что ввиду (2.1.3), (2.1.4) при любом фиксированном 
$ x_j \in \R $ функция $ f_{\nu_j}^j(x_1, \ldots, x_{j -1}, x_j, 
x_{j +1}, \ldots, x_d) \in \mathcal P_{\kappa^{J_j}}^{d -1, l^{J_j}, m^{J_j}}, $ 
а, значит, $ \supp f_{\nu_j}^j(x_1, \ldots, x_{j -1}, x_j, x_{j +1}, \ldots, x_d) \subset 
Q^{d -1,m^{J_j}}, \nu_j \in \Nu_{-m_j, 2^{\kappa_j} -1}^1. $
заключаем, что (2.1.5) имеет место при $ x_j \in \R, x^{J_j} \in \R^{d -1}, $  
или $ x \in \R^d. $ Точно так же устанавливается справедливость замечания 
в отношении условия (2.1.6).

Лемма 2.1.1 
 
При $ d \in \N, l \in \Z_+^d, m \in \N^d, \kappa \in \Z_+^d $ для того, чтобы 
для функции $ f: \R^d \mapsto \R $ имело место включение 
$ f \in \mathcal P_\kappa^{\prime d,l,m}, $ 
необходимо и достаточно, чтобы $ f \in \mathcal P_\kappa^{d,l,m} $ и при 
$ j =1, \ldots, d $ для каждой точки 
$ (x_1, \ldots, x_{j -1}, x_{j +1}, \ldots, x_d) \in \R^{d -1} $ соблюдалось включение
\begin{equation*} \tag{2.1.7} 
f(x_1, \ldots, x_{j -1}, \cdot, x_{j +1}, \ldots, x_d) \in 
\mathcal P_{\kappa_j}^{\prime 1, l_j, m_j}.
\end{equation*} 
    
Доказательство.

Если $ f \in \mathcal P_\kappa^{\prime d,l,m}, $ то $ f \in \mathcal P_\kappa^{d,l,m} $
и в силу (2.1.2), (2.1.3), (2.1.5), (2.1.6) при $ j =1, \ldots, d $ для 
каждой точки $ (x_1, \ldots, x_{j -1}, x_{j +1}, \ldots, x_d) \in \R^{d -1} $ 
на $ \R $ выполняется равенство
\begin{equation*} \tag{2.1.8}
f(x_1, \ldots, x_{j -1}, \cdot, x_{j +1}, \ldots, x_d) = \sum_{\nu_j \in 
\Nu_{-m_j, 2^{\kappa_j} -1}^1} f_{\nu_j}^j(x_1, \ldots, x_{j -1}, \cdot, 
x_{j +1}, \ldots, x_d) g_{\kappa_j, \nu_j}^{1, m_j}(\cdot),
\end{equation*}
а при $ \nu_j \in \Nu_{-m_j, 2^{\kappa_j} -1}^1 $ имеет место включение 
\begin{equation*} \tag{2.1.9} 
f_{\nu_j}^j(x_1, \ldots, x_{j -1}, \cdot, x_{j +1}, \ldots, x_d) \in 
\mathcal P^{1, l_j}, 
\end{equation*} 
и соблюдаются равенства
\begin{equation*} \tag{2.1.10}
f_{\nu_j}^j(x_1, \ldots, x_{j -1}, \cdot, x_{j +1}, \ldots, x_d) = 
f_0^j(x_1, \ldots, x_{j -1}, \cdot, x_{j +1}, \ldots, x_d), \   
\nu_j \in \Nu_{-m_j, 0}^1; \\ 
\end{equation*}
и
\begin{multline*} \tag{2.1.11}
f_{\nu_j}^j(x_1, \ldots, x_{j -1}, \cdot, x_{j +1}, \ldots, x_d) = \\
f_{2^{\kappa_j} -m_j -1}^j(x_1, \ldots, x_{j -1}, \cdot, x_{j +1}, \ldots, x_d), \  
\nu_j \in \Nu_{2^{\kappa_j} -m_j -1, 2^{\kappa_j} -1}^1, \ 
\end{multline*}
т.е. $ f(x_1, \ldots, x_{j -1}, \cdot, x_{j +1}, \ldots, x_d) \in 
\mathcal P_{\kappa_j}^{\prime 1, l_j, m_j}. $ 

Обратно, если $ f \in \mathcal P_\kappa^{d,l,m} $ и при 
$ j =1, \ldots, d $ для каждой точки 
$ (x_1, \ldots, x_{j -1}, x_{j +1}, \ldots, x_d) \in \R^{d -1} $ 
справедливо включение (2.1.7), то, рассматривая представления 
(2.1.2), (2.1.3), заключаем, что при $ j  =1, \ldots, d $ для каждой точки 
$ (x_1, \ldots, x_{j -1}, x_{j +1}, \ldots, x_d) \in  \R^{d -1} $ выполняются 
соотношения (2.1.8), (2.1.9), (2.1.10), (2.1.11), т.е. соблюдаются равенства 
(2.1.5), (2.1.6), тем самым, верно включение
$ f \in \mathcal P_\kappa^{\prime d, l, m}. \square $
 
Лемма 2.1.2   

Пусть $ l \in \Z_+, m \in \N, \lambda \in \Z_+^1(m), 1 \le q \le \infty.$ Тогда 
существуют константы $ c_1(l,m,\lambda,q) > 0, c_2(l,m,\lambda,q) > 0, 
c_3(l,m,\lambda,q) > 0 $ такие, что для $ f \in \mathcal P_0^{\prime 1, l, m} $ 
справедливы неравенства
\begin{equation*} \tag{2.1.12}
\| \D^\lambda f \|_{L_q(-m +(m +1) I)} \le c_1 \| f \|_{L_q(I)};
\end{equation*}
\begin{equation*} \tag{2.1.13}
\| \D^\lambda f \|_{L_q((m +1) I)} \le c_2 \| f \|_{L_q(I)}, \ 
\end{equation*}
а для $ f \in \mathcal P_0^{1, l, m} $ выполняется неравенство
\begin{equation*} \tag{2.1.14}
\| \D^\lambda f \|_{L_q(I)} \le c_3 \| f \|_{L_q(I)}.
\end{equation*}

Доказательство.

Заметим, что в условиях леммы функционал $ \| f \|_{L_q(I)}, f \in 
\mathcal P_0^{\prime 1, l, m}, $ является нормой на $ \mathcal P_0^{\prime 1,l,m}. $
То, что этот функционал является полунормой на $ \mathcal P_0^{\prime 1, l, m}, $ 
очевидно. Проверим, что если для $ f \in \mathcal P_0^{\prime 1, l, m} $ 
значение $ \| f \|_{L_q(I)} = 0, $ то $ f =0. $ Пусть для $ f \in 
\mathcal P_0^{\prime 1, l, m} $ выполнено равенство $ \| f \|_{L_q(I)} = 0. $ 
Тогда для $ x \in I $ справедливо равенство $ f(x) =0. $ Отсюда, учитывая, 
что в силу (1.3.4), (2.1.5), (2.1.6)  для $ x \in I $ имеет место соотношение 
$ 0 = f(x) = \sum_{\nu \in \Nu_{-m,0}^1} f_0(x) g_{0, \nu}^{1,m}(x) =
f_0(x) \sum_{\nu \in \Nu_{-m,0}^1} g_{0, \nu}^{1,m}(x) = f_0(x), \  
f_0 \in \mathcal P^{1,l}. $ 
Следовательно, $ f_0(x) =0 $ для $ x \in \R, $ и, значит, 
$ f(x) = \sum_{\nu \in \Nu_{-m,0}^1} f_0(x) g_{0, \nu}^{1,m}(x) =0 $ для $ x \in \R. $

Теперь, принимая во внимание, что линейное отображение конечномерного 
нормированного пространства в нормированное пространство является 
непрерывным, заключаем, что для отображения
$$
\mathcal P_0^{\prime 1,l,m} \cap L_q(I) \ni f \mapsto (\D^\lambda f) \mid_{-m +(m +1) I} \in
L_q(-m +(m +1) I), \ \lambda \in \Z_+^1(m), 1 \le q \le \infty, \ 
$$
существует константа $ c_1(l,m,\lambda,q) > 0 $ такая, что для 
$ f \in \mathcal P_0^{\prime 1, l, m} $ 
соблюдается неравенство (2.1.12).

Точно так же устанавливается справедливость (2.1.13).

В силу тех же соображений получаем, что для отображения
$$
\{ F \mid_I: F \in \mathcal P_0^{1,l,m}\} \cap L_q(I) \ni f \mid_I \mapsto 
\D^\lambda (f \mid_I) \in L_q(I), \ \lambda \in \Z_+^1(m), 1 \le q \le \infty, \ 
$$
существует константа $ c_3(l,m,\lambda,q) > 0 $ такая, что для 
$ f \in \mathcal P_0^{1, l, m} $ выполняется неравенство (2.1.14). $ \square $

Предложение 2.1.3

Пусть $ d \in \N, l \in \Z_+^d, m \in \N^d, \lambda \in \Z_+^d(m), 
1 \le q \le \infty. $ Тогда существует константа $ c_4(d,l,m,\lambda,q) >0 $ 
такая, что при $ \kappa \in \Z_+^d $ для $ f \in \mathcal P_\kappa^{\prime d,l,m} $ 
имеет место неравенство
\begin{equation*} \tag{2.1.15}
\| \D^\lambda f \|_{L_q(\R^d)} \le c_4 2^{(\kappa, \lambda)} \| f \|_{L_q(I^d)}.
\end{equation*}

Доказательство.

Справедливость предложения установим по индукции относительно $ d. $ Сначала 
проверим соблюдение (2.1.15) при $ d =1. $ 

При $ l \in \Z_+, m \in \N, \lambda \in \Z_+^1(m), \kappa \in \Z_+,  
1 \le q \le \infty $ Для $ f \in \mathcal P_\kappa^{\prime 1,l,m}, $ 
благодаря свойствам функции $ \psi^{1,m}, $ имеем
\begin{multline*} \tag{2.1.16}
\| \D^\lambda f \|_{L_q(\R)}^q = \int_\R | \D^\lambda f |^q dx =
\int_\R | \D^\lambda(\sum_{\nu \in \Nu_{-m, 2^\kappa -1}^1} f_\nu(x) 
g_{\kappa, \nu}^{1, m}(x))|^q dx = \\
\int_{-2^{-\kappa} m}^{1 +2^{-\kappa} m} 
| \D^\lambda (\sum_{\nu \in \Nu_{-m, 2^\kappa -1}^1} f_\nu(x) \psi^{1,m}(2^\kappa x -\nu))|^q dx =\\
\int_{-2^{-\kappa} m +2^{-\kappa} (m +1) I} 
| \D^\lambda(\sum_{\nu \in \Nu_{-m, 2^\kappa -1}^1} f_\nu(x) 
\psi^{1,m}(2^\kappa x -\nu))|^q dx +\\
\sum_{\nu^\prime \in \Nu_{1, 2^\kappa -2}^1} \int_{2^{-\kappa} \nu^\prime +2^{-\kappa} I}
| \D^\lambda(\sum_{\nu \in \Nu_{-m, 2^\kappa -1}^1} f_\nu(x) 
\psi^{1,m}(2^\kappa x -\nu))|^q dx +\\
\int_{1 -2^{-\kappa} +2^{-\kappa} (m +1) I} 
| \D^\lambda(\sum_{\nu \in \Nu_{-m, 2^\kappa -1}^1} f_\nu(x) 
\psi^{1,m}(2^\kappa x -\nu))|^q dx =\\
\int_{-2^{-\kappa} m +2^{-\kappa} (m +1) I} 
| \D^\lambda(\sum_{\nu \in \Nu_{-m, 0}^1} f_\nu(x) 
\psi^{1,m}(2^\kappa x -\nu))|^q dx +\\
\sum_{\nu^\prime \in \Nu_{1, 2^\kappa -2}^1} \int_{2^{-\kappa} \nu^\prime +2^{-\kappa} I}
| \D^\lambda(\sum_{\nu \in \nu^\prime +\Nu_{-m, 0}^1} f_\nu(x) 
\psi^{1,m}(2^\kappa x -\nu))|^q dx +\\
\int_{1 -2^{-\kappa} +2^{-\kappa} (m +1) I} 
| \D^\lambda(\sum_{\nu \in \Nu_{2^\kappa -m -1, 2^\kappa -1}^1} f_\nu(x) 
\psi^{1,m}(2^\kappa x -\nu))|^q dx.
\end{multline*}

Оценим слагаемые в правой части (2.1.16). Для этого, выбирая с учётом 
(2.1.5) полином $ F_0 \in \mathcal P^{1, l} $ такой, что при 
$ \nu \in \Nu_{-m, 0}^1 $ соблюдается равенство $ f_\nu(x) = F_0(2^\kappa x), \  
x \in \R, $ после дифференцирования, 
замены переменной в интеграле, применения (2.1.12) и обратной замены переменной в интеграле получаем
\begin{multline*} \tag{2.1.17}
\int_{-2^{-\kappa} m +2^{-\kappa} (m +1) I} 
| \D^\lambda(\sum_{\nu \in \Nu_{-m, 0}^1} f_\nu(x) \psi^{1,m}(2^\kappa x -\nu))|^q dx =\\
\int_{-2^{-\kappa} m +2^{-\kappa} (m +1) I} | \frac{d^\lambda}{dx^\lambda}
(\sum_{\nu \in \Nu_{-m, 0}^1} F_0(2^\kappa x) \psi^{1,m}(2^\kappa x -\nu))|^q dx =\\
\int_{-2^{-\kappa} m +2^{-\kappa} (m +1) I} | \frac{d^\lambda}{dy^\lambda}
(\sum_{\nu \in \Nu_{-m, 0}^1} F_0(y) \psi^{1,m}(y -\nu)) 
\mid_{y = 2^\kappa x} 2^{\kappa \lambda} |^q dx =\\
2^{\kappa \lambda q} \int_{-2^{-\kappa} m +2^{-\kappa} (m +1) I} 
| \frac{d^\lambda}{dy^\lambda}
(\sum_{\nu \in \Nu_{-m, 0}^1} F_0(y) g_{0, \nu}^{1,m}(y)) 
\mid_{y = 2^\kappa x} |^q dx =\\
2^{\kappa \lambda q} \int_{-m +(m +1) I} | \frac{d^\lambda}{dy^\lambda}
(\sum_{\nu \in \Nu_{-m, 0}^1} F_0(y) g_{0, \nu}^{1,m}(y))|^q 2^{-\kappa} dy \le \\
2^{\kappa \lambda q} 2^{-\kappa} c_1^q  \int_{I} | \sum_{\nu \in \Nu_{-m, 0}^1} 
F_0(y) g_{0, \nu}^{1,m}(y)|^q dy =\\
c_1^q 2^{\kappa \lambda q} 2^{-\kappa}  \int_{I} | \sum_{\nu \in \Nu_{-m, 0}^1} 
F_0(y) \psi^{1,m}(y -\nu)|^q dy =\\
c_1^q 2^{\kappa \lambda q} 2^{-\kappa}  \int_{2^{-\kappa} I} | \sum_{\nu \in \Nu_{-m, 0}^1} 
F_0(2^\kappa x) \psi^{1,m}(2^\kappa x -\nu)|^q 2^\kappa dx =\\
c_1^q 2^{\kappa \lambda q} \int_{2^{-\kappa} I} | \sum_{\nu \in \Nu_{-m, 0}^1} 
f_\nu(x) \psi^{1,m}(2^\kappa x -\nu)|^q dx.
\end{multline*}

Далее, так же подбирая полином $ \mathcal F_0 \in \mathcal P^{1,l}, $ 
для которого ввиду (2.1.6) при $ \nu \in \Nu_{2^\kappa -m -1, 2^\kappa -1}^1 $ 
выполняется равенство $ f_\nu(x) = \mathcal F_0(2^\kappa x -2^\kappa +1), $ 
проводя дифференцирование, делая замену переменной в интеграле, используя 
(2.1.13) и ещё раз делая замену переменной в интеграле, выводим
\begin{multline*} \tag{2.1.18}
\int_{1 -2^{-\kappa} +2^{-\kappa} (m +1) I} 
| \D^\lambda(\sum_{\nu \in \Nu_{2^\kappa -m -1, 2^\kappa -1}^1} f_\nu(x) 
\psi^{1,m}(2^\kappa x -\nu))|^q dx =\\
\int_{1 -2^{-\kappa} +2^{-\kappa} (m +1) I} | \frac{d^\lambda}{dx^\lambda}
(\sum_{\nu \in \Nu_{2^\kappa -m -1, 2^\kappa -1}^1} \mathcal F_0(2^\kappa x 
-2^\kappa +1) \psi^{1,m}(2^\kappa x -2^\kappa +1 -(\nu -2^\kappa +1)))|^q dx =\\
\int_{1 -2^{-\kappa} +2^{-\kappa} (m +1) I} | \frac{d^\lambda}{dx^\lambda}
(\sum_{\nu \in \Nu_{-m, 0}^1} \mathcal F_0(2^\kappa x -2^\kappa +1) 
\psi^{1,m}(2^\kappa x -2^\kappa +1 -\nu))|^q dx =\\
\int_{1 -2^{-\kappa} +2^{-\kappa} (m +1) I} | \frac{d^\lambda}{dy^\lambda}
(\sum_{\nu \in \Nu_{-m, 0}^1} \mathcal F_0(y) 
\psi^{1,m}(y -\nu)) \mid_{y = 2^\kappa x -2^\kappa +1} 2^{\kappa \lambda}|^q dx =\\
2^{\kappa \lambda q}
\int_{1 -2^{-\kappa} +2^{-\kappa} (m +1) I} | \frac{d^\lambda}{dy^\lambda}
(\sum_{\nu \in \Nu_{-m, 0}^1} \mathcal F_0(y) 
g_{0, \nu}^{1,m}(y)) \mid_{y = 2^\kappa x -2^\kappa +1} |^q dx =\\
2^{\kappa \lambda q}
\int_{(m +1) I} | \frac{d^\lambda}{dy^\lambda}
(\sum_{\nu \in \Nu_{-m, 0}^1} \mathcal F_0(y) g_{0, \nu}^{1,m}(y)) |^q 2^{-\kappa} dy \le\\
2^{\kappa \lambda q} 2^{-\kappa} c_2^q
\int_{I} | \sum_{\nu \in \Nu_{-m, 0}^1} \mathcal F_0(y) g_{0, \nu}^{1,m}(y)|^q dy =\\
c_2^q 2^{\kappa \lambda q} 2^{-\kappa} 
\int_{I} | \sum_{\nu \in \Nu_{-m, 0}^1} \mathcal F_0(y) \psi^{1,m}(y -\nu)|^q dy =\\
c_2^q 2^{\kappa \lambda q} 2^{-\kappa} 
\int_{1 -2^{-\kappa} +2^{-\kappa} I} | \sum_{\nu \in \Nu_{-m, 0}^1} 
\mathcal F_0(2^\kappa x -2^\kappa +1) 
\psi^{1,m}(2^\kappa x -(\nu +2^\kappa -1))|^q 2^\kappa dx =\\
c_2^q 2^{\kappa \lambda q} 
\int_{1 -2^{-\kappa} +2^{-\kappa} I} | \sum_{\nu \in \Nu_{2^\kappa -m -1, 2^\kappa -1}^1} 
\mathcal F_0(2^\kappa x -2^\kappa +1) 
\psi^{1,m}(2^\kappa x -\nu)|^q dx =\\
c_2^q 2^{\kappa \lambda q} 
\int_{1 -2^{-\kappa} +2^{-\kappa} I} | \sum_{\nu \in \Nu_{2^\kappa -m -1, 2^\kappa -1}^1} 
f_\nu(x) \psi^{1,m}(2^\kappa x -\nu)|^q dx.
\end{multline*}

Наконец, при $ \nu^\prime \in \Nu_{1, 2^\kappa -2}^1 $ и $ \nu \in \nu^\prime +
\Nu_{-m, 0}^1, $ беря полином $ F_\nu^{\nu^\prime} \in \mathcal P^{1, l} $ 
такой,  что имеет место равенство 
$ f_\nu(x) = F_\nu^{\nu^\prime}(2^\kappa x -\nu^\prime), x \in \R, $
как и выше, проводя дифференцирование, замены переменных в интегралах и 
пользуясь (2.1.14), приходим к неравенству
\begin{multline*} \tag{2.1.19} 
\int_{2^{-\kappa} \nu^\prime +2^{-\kappa} I}
| \D^\lambda(\sum_{\nu \in \nu^\prime +\Nu_{-m, 0}^1} f_\nu(x) \psi^{1,m}(2^\kappa x -\nu))|^q dx = \\
\int_{2^{-\kappa} \nu^\prime +2^{-\kappa} I} | \frac{d^\lambda}{dx^\lambda}
(\sum_{\nu \in \nu^\prime +\Nu_{-m, 0}^1} F_\nu^{\nu^\prime}(2^\kappa x -\nu^\prime) 
\psi^{1,m}(2^\kappa x -\nu^\prime -(\nu -\nu^\prime)))|^q dx = \\
\int_{2^{-\kappa} \nu^\prime +2^{-\kappa} I} | \frac{d^\lambda}{dx^\lambda}
(\sum_{\nu \in \Nu_{-m, 0}^1} F_{\nu^\prime +\nu}^{\nu^\prime}(2^\kappa x -\nu^\prime) 
\psi^{1,m}(2^\kappa x -\nu^\prime -\nu))|^q dx = \\
\int_{2^{-\kappa} \nu^\prime +2^{-\kappa} I} | \frac{d^\lambda}{dy^\lambda}
(\sum_{\nu \in \Nu_{-m, 0}^1} F_{\nu^\prime +\nu}^{\nu^\prime}(y) 
\psi^{1,m}(y -\nu)) \mid_{y = 2^\kappa x -\nu^\prime} 2^{\kappa \lambda}|^q dx =\\
2^{\kappa \lambda q} \int_{I} | \frac{d^\lambda}{dy^\lambda}
(\sum_{\nu \in \Nu_{-m, 0}^1} F_{\nu^\prime +\nu}^{\nu^\prime}(y) 
\psi^{1,m}(y -\nu))|^q 2^{-\kappa} dy = \\
2^{\kappa \lambda q} 2^{-\kappa} \int_{I} | \frac{d^\lambda}{dy^\lambda}
(\sum_{\nu \in \Nu_{-m, 0}^1} F_{\nu^\prime +\nu}^{\nu^\prime}(y) 
g_{0, \nu}^{1,m}(y))|^q dy \le \\
2^{\kappa \lambda q} 2^{-\kappa} c_3^q \int_{I}
| \sum_{\nu \in \Nu_{-m, 0}^1} F_{\nu^\prime +\nu}^{\nu^\prime}(y) 
g_{0, \nu}^{1,m}(y)|^q dy = \\
c_3^q 2^{\kappa \lambda q} 2^{-\kappa} \int_{I}
| \sum_{\nu \in \Nu_{-m, 0}^1} F_{\nu^\prime +\nu}^{\nu^\prime}(y) 
\psi^{1,m}(y -\nu)|^q dy = \\
c_3^q 2^{\kappa \lambda q} 2^{-\kappa} 
\int_{2^{-\kappa} \nu^\prime +2^{-\kappa} I}
| \sum_{\nu \in \Nu_{-m, 0}^1} F_{\nu^\prime +\nu}^{\nu^\prime}(2^\kappa x -\nu^\prime) 
\psi^{1,m}(2^\kappa x -\nu^\prime -\nu)|^q 2^\kappa dx = \\
c_3^q 2^{\kappa \lambda q} \int_{2^{-\kappa} \nu^\prime +2^{-\kappa} I}
| \sum_{\nu \in \Nu_{-m, 0}^1} F_{\nu^\prime +\nu}^{\nu^\prime}(2^\kappa x -\nu^\prime) 
\psi^{1,m}(2^\kappa x -(\nu^\prime +\nu))|^q dx = \\
c_3^q 2^{\kappa \lambda q} \int_{2^{-\kappa} \nu^\prime +2^{-\kappa} I}
| \sum_{\nu \in \nu^\prime +\Nu_{-m, 0}^1} F_\nu^{\nu^\prime}(2^\kappa x -\nu^\prime) 
\psi^{1,m}(2^\kappa x -\nu)|^q dx = \\
c_3^q 2^{\kappa \lambda q} \int_{2^{-\kappa} \nu^\prime +2^{-\kappa} I}
| \sum_{\nu \in \nu^\prime +\Nu_{-m, 0}^1} f_\nu(x) 
\psi^{1,m}(2^\kappa x -\nu)|^q dx.
\end{multline*}
 
Объединяя (2.1.16), (2.1.17), (2.1.18), (2.1.19) и учитывая (1.3.3), находим, 
что для $ f \in \mathcal P_\kappa^{\prime 1,l,m} $ справедливо неравенство 
\begin{multline*} \tag{2.1.20}
\| \D^\lambda f \|_{L_q(\R)}^q \le 
c_1^q 2^{\kappa \lambda q} \int_{2^{-\kappa} I} | \sum_{\nu \in \Nu_{-m, 0}^1} 
f_\nu(x) \psi^{1,m}(2^\kappa x -\nu)|^q dx + \\
\sum_{\nu^\prime \in \Nu_{1, 2^\kappa -2}^1} 
c_3^q 2^{\kappa \lambda q} \int_{2^{-\kappa} \nu^\prime +2^{-\kappa} I}
| \sum_{\nu \in \nu^\prime +\Nu_{-m, 0}^1} f_\nu(x) 
\psi^{1,m}(2^\kappa x -\nu)|^q dx + \\
c_2^q 2^{\kappa \lambda q} \int_{1 -2^{-\kappa} +2^{-\kappa} I} 
| \sum_{\nu \in \Nu_{2^\kappa -m -1, 2^\kappa -1}^1} 
f_\nu(x) \psi^{1,m}(2^\kappa x -\nu)|^q dx \le \\
c_4^q 2^{\kappa \lambda q} \sum_{\nu^\prime \in \Nu_{0, 2^\kappa -1}^1} 
\int_{2^{-\kappa} \nu^\prime +2^{-\kappa} I}
| \sum_{\nu \in \nu^\prime +\Nu_{-m, 0}^1} f_\nu(x) 
\psi^{1,m}(2^\kappa x -\nu)|^q dx = \\
c_4^q 2^{\kappa \lambda q} \sum_{\nu^\prime \in \Nu_{0, 2^\kappa -1}^1} 
\int_{2^{-\kappa} \nu^\prime +2^{-\kappa} I}
| \sum_{\nu \in \Nu_{-m, 2^\kappa -1}^1} f_\nu(x) 
g_{\kappa, \nu}^{1,m}(x)|^q dx = \\
c_4^q 2^{\kappa \lambda q} \sum_{\nu^\prime \in \Nu_{0, 2^\kappa -1}^1} 
\int_{2^{-\kappa} \nu^\prime +2^{-\kappa} I} | f(x)|^q dx = \\
c_4^q 2^{\kappa \lambda q} \int_I | f(x)|^q dx =
(c_4 2^{\kappa \lambda} \| f \|_{L_q(I)})^q. 
\end{multline*}
Из (2.1.20) следует (2.1.15) при $ d =1. $ 

Предположим, что при $ d \in \N: d \ge 2. $ предложение справедливо для 
размерности $ d -1 $ вместо $ d. $ Покажем, что тогда оно верно и для 
размерности $ d. $ Для этого в силу 
теоремы Фубини для $ f \in \mathcal P_\kappa^{\prime d,l,m}, \lambda 
\in \Z_+^d(m) $ имеем
\begin{multline*} \tag{2.1.21}
\| \D^\lambda f \|_{L_q(\R^d)}^q =  
\int_{\R^d} | \D^\lambda f(x)|^q dx = \\  
\int_{\R^{d -1}} \int_\R | \D^\lambda f(x_1, \ldots, x_{d -1}, x_d)|^q dx_d
dx_1 \ldots dx_{d -1} = \\ 
\int_{\R^{d -1}} \int_\R | (\D^{\lambda e_d} (\D^{\lambda \chi_{J_d}} 
f))(x_1, \ldots, x_{d -1}, x_d)|^q dx_d dx_1 \ldots dx_{d -1}.
\end{multline*}  

Далее, используя (2.1.2), при $ (x_1, \ldots, x_{d -1}, x_d) \in \R^d $ 
получаем
\begin{multline*} 
(\D^{\lambda \chi_{J_d}} f)(x_1, \ldots, x_{d -1}, x_d) =
\D^{\lambda \chi_{J_d}} (\sum_{\nu_d \in \Nu_{-m_d, 2^{\kappa_d} -1}^1} 
f_{\nu_d}^d(x_1, \ldots, x_{d -1}, x_d) g_{\kappa_d, \nu_d}^{1, m_d}(x_d)) = \\ 
\sum_{\nu_d \in \Nu_{-m_d, 2^{\kappa_d} -1}^1} \D^{\lambda \chi_{J_d}} 
(f_{\nu_d}^d(x_1, \ldots, x_{d -1}, x_d) g_{\kappa_d, \nu_d}^{1, m_d}(x_d)) = \\
\sum_{\nu_d \in \Nu_{-m_d, 2^{\kappa_d} -1}^1} 
(\D^{\lambda \chi_{J_d}} f_{\nu_d}^d)(x_1, \ldots, x_{d -1}, x_d)
g_{\kappa_d, \nu_d}^{1, m_d}(x_d).
\end{multline*}

Отсюда, принимая во внимание, что при $ \nu_d \in \Nu_{-m_d, 2^{\kappa_d} -1}^1 $ 
в каждой точке $ (x_1, \ldots, x_{d -1}) \in \R^{d -1} $ в силу (2.1.3), (2.1.4)
выполняется соотношение
\begin{multline*}
(\D^{\lambda \chi_{J_d}} f_{\nu_d}^d)(x_1, \ldots, x_{d -1}, x_d) =
\D^{\lambda \chi_{J_d}} (\sum_{\mu_d =0, \ldots, l_d} x_d^{\mu_d} 
f_{\nu_d}^{d, \mu_d}(x^{J_d})) = \\
\sum_{\mu_d =0, \ldots, l_d} \D^{\lambda \chi_{J_d}} (x_d^{\mu_d} 
f_{\nu_d}^{d, \mu_d}(x^{J_d})) = 
\sum_{\mu_d =0, \ldots, l_d} x_d^{\mu_d}
(\D^{\lambda \chi_{J_d}} f_{\nu_d}^{d, \mu_d})(x^{J_d}) \in \mathcal P^{1, l_d}, \ 
\end{multline*}
где
$$
f_{\nu_d}^{d, \mu_d}(x^{J_d}) \in 
\mathcal P_{\kappa^{J_d}}^{d -1, l^{J_d}, m^{J_d}}, \  
\mu_d =0, \ldots, l_d, \nu_d \in \Nu_{-m_d, 2^{\kappa_d} -1}^1, \ 
$$ 
и, благодаря (2.1.5), (2.1.6), для всех $ (x_1, \ldots, x_{d -1}) \in 
\R^{d -1} $ при любом $ x_d \in \R $ соблюдаются равенства
\begin{multline*}
(\D^{\lambda \chi_{J_d}} f_{\nu_d}^d)(x_1, \ldots, x_{d -1}, x_d) =
(\D^{\lambda \chi_{J_d}} f_{0}^d)(x_1, \ldots, x_{d -1}, x_d), \  
\nu_d \in \Nu_{-m_d, 0}^1; \\
(\D^{\lambda \chi_{J_d}} f_{\nu_d}^d)(x_1, \ldots, x_{d -1}, x_d) =
(\D^{\lambda \chi_{J_d}} f_{2^{\kappa_d} -m_d -1}^d)(x_1, \ldots, x_{d -1}, x_d), \  
\nu_d \in \Nu_{2^{\kappa_d} -m_d -1, 2^{\kappa_d} -1}^1, \ 
\end{multline*}
заключаем, что для $ (x_1, \ldots, x_{d -1}) \in \R^{d -1} $ функция 
$$ 
(\D^{\lambda \chi_{J_d}} f)(x_1, \ldots, x_{d -1}, \cdot) = \\
(\sum_{\nu_d \in \Nu_{-m_d, 2^{\kappa_d} -1}^1} 
(\D^{\lambda \chi_{J_d}} f_{\nu_d}^d)(x_1, \ldots, x_{d -1}, \cdot)
g_{\kappa_d, \nu_d}^{1, m_d}(\cdot)) \in \mathcal P_{\kappa_d}^{\prime 1, l_d, m_d}. 
$$
Поэтому, используя (2.1.15) при $ d =1, $ для $ (x_1, \ldots, x_{d -1}) \in 
\R^{d -1} $ имеем оценку
\begin{multline*} 
\int_\R | (\D^{\lambda e_d} (\D^{\lambda \chi_{J_d}} 
f))(x_1, \ldots, x_{d -1}, x_d)|^q dx_d = \\
\int_\R | \frac{\D^{\lambda_d}}{\D x_d^{\lambda_d}} (\D^{\lambda \chi_{J_d}} 
f)(x_1, \ldots, x_{d -1}, x_d)|^q dx_d \le \\ 
(c_4(1, l_d, m_d, \lambda_d, q))^q
2^{\kappa_d \lambda_d q} \int_I | \D^{\lambda \chi_{J_d}} 
f(x_1, \ldots, x_{d -1}, x_d)|^q dx_d.
\end{multline*}
Подставляя эту оценку в (2.1.21) и применяя теорему Фубини, находим, что
\begin{multline*} \tag{2.1.22} 
\| \D^\lambda f \|_{L_q(\R^d)}^q \le \\
\int_{\R^{d -1}} (c_4(1, l_d, m_d, \lambda_d, q))^q 2^{\kappa_d \lambda_d q} 
\int_I | \D^{\lambda \chi_{J_d}} 
f(x_1, \ldots, x_{d -1}, x_d)|^q dx_d dx_1 \ldots dx_{d -1} = \\ 
(c_4(1, l_d, m_d, \lambda_d, q))^q 2^{\kappa_d \lambda_d q} \int_I \int_{\R^{d -1}} 
| \D^{\lambda \chi_{J_d}} f(x_1, \ldots, x_{d -1}, x_d)|^q dx_1 \ldots dx_{d -1} dx_d.
\end{multline*}
  
Теперь заметим, что для $ f \in \mathcal P_\kappa^{\prime d, l, m} $ при 
$ x_d \in I $ ввиду (2.1.2), (2.1.3), (2.1.4) функция 
$$ 
f(x_1, \ldots, x_{d -1}, x_d) \in \mathcal P_{\kappa^{J_d}}^{d -1, l^{J_d}, m^{J_d}}, 
$$
и для каждого $ x_d \in I $ при $ j =1, \ldots, d -1 $ для любых 
$ (x_1, \ldots, x_{j -1}, x_{j +1}, \ldots, x_{d -1}) \in \R^{d -2} $ 
в соответствии с леммой 2.1.1 справедливо включение (2.1.7), а именно, 
\begin{equation*} 
f(x_1, \ldots, x_{j -1}, \cdot, x_{j +1}, \ldots, x_{d -1}, x_d) \in 
\mathcal P_{\kappa_j}^{\prime 1, l_j, m_j}, \ 
\end{equation*} 
и, значит, в силу леммы 2.1.1 при $ x_d \in I $ функция 
$$ 
f(x_1, \ldots, x_{d -1}, x_d) \in \mathcal P_{\kappa^{J_d}}^{\prime d -1, l^{J_d}, m^{J_d}}. 
$$
Учитывая это обстоятельство, на основании предположения индукции 
получаем, что для $ x_d \in I $ выполняется неравенство
\begin{multline*}
\int_{\R^{d -1}} | \D^{\lambda \chi_{J_d}} 
f(x_1, \ldots, x_{d -1}, x_d)|^q dx_1 \ldots dx_{d -1} = \\
\int_{\R^{d -1}} | \D^{\lambda^{J_d}} 
f(x_1, \ldots, x_{d -1}, x_d)|^q dx_1 \ldots dx_{d -1} \le \\ 
(c_4(d -1, l^{J_d}, m^{J_d}, \lambda^{J_d}, q))^q
2^{(\kappa^{J_d}, \lambda^{J_d}) q} \int_{I^{d -1}} | f(x_1, \ldots, x_{d -1}, x_d)|^q dx_1 \ldots dx_{d -1}. 
\end{multline*}
Подставляя эту оценку в (2.1.22), приходим к неравенству
\begin{multline*}
\| \D^\lambda f \|_{L_q(\R^d)}^q \le
(c_4(1, l_d, m_d, \lambda_d, q))^q 2^{\kappa_d \lambda_d q} \\
\times\int_I (c_4(d -1, l^{J_d}, m^{J_d}, \lambda^{J_d}, q))^q
2^{(\kappa^{J_d}, \lambda^{J_d}) q} \int_{I^{d -1}} 
| f(x_1, \ldots, x_{d -1}, x_d)|^q dx_1 \ldots dx_{d -1} dx_d = \\
(c_4(1, l_d, m_d, \lambda_d, q))^q 2^{\kappa_d \lambda_d q} 
(c_4(d -1, l^{J_d}, m^{J_d}, \lambda^{J_d}, q))^q
2^{(\kappa^{J_d}, \lambda^{J_d}) q}\\
\times \int_I \int_{I^{d -1}} 
| f(x_1, \ldots, x_{d -1}, x_d)|^q dx_1 \ldots dx_{d -1} dx_d = \\
(c_4(d, l, m, \lambda, q))^q 2^{(\kappa, \lambda) q} \int_{I^d} | f(x)|^q dx, \ 
\end{multline*}
из которого вытекает 
\begin{equation*}
\| \D^\lambda f \|_{L_q(\R^d)} \le c_4 2^{(\kappa, \lambda)} \| f \|_{L_q(I^d)}, \ 
\end{equation*}
что совпадает с (2.1.15). $ \square $
\bigskip

2.2. В этом пункте вводятся в рассмотрение линейные операторы, 
используемые в следующем пункте для доказательства основного результата 
работы.

При $ d \in \N, l \in \Z_+^d, m \in \N^d, \kappa \in \Z_+^d, 
\nu \in \Nu_{0, 2^\kappa -\e}^d $ определим линейный оператор 
$$ 
\mathcal S_{\kappa, \nu}^{d, l, m,\R}: L^{\loc}(\R^d) \mapsto L_\infty(\R^d), \ 
$$
полагая для $ f \in L^{\loc}(\R^d) $ значение
\begin{equation*} \tag{2.2.1}
(\mathcal S_{\kappa, \nu}^{d, l, m, \R} f)(x) = 
\chi_{Q^{d, m}}(x) (P_{\delta, x^0}^{d, l}(f \mid_{x^0 +\delta I^d}))(x), \ x \in \R^d, \ 
\end{equation*}
при $ \delta = 2^{-\kappa}, x^0 = 2^{-\kappa} \nu. $

Ясно, что при $ d \in \N, l \in \Z_+^d, m \in \N^d, \kappa \in \Z_+^d, 
\nu \in \Nu_{0, 2^\kappa -\e}^d, 1 \le p \le \infty $ сужение отображения 
$ \mathcal S_{\kappa, \nu}^{d,l,m,\R} $ на пространство $ L_p(\R^d) $ 
принадлежит $ \mathcal B(L_p(\R^d), L_p(\R^d)). $  
Из (1.2.8) вытекает, что при $ d \in \N, l \in \Z_+^d, m \in \N^d, 
\kappa \in \Z_+^d, \nu \in \Nu_{0, 2^\kappa -\e}^d, 1 \le p < \infty $ 
для $ f \in L_p(\R^d) $ имеет место равенство
\begin{equation*} \tag{2.2.2}
\mathcal S_{\kappa, \nu}^{d,l,m,\R} f = 
(\prod_{j =1}^d V_j(\mathcal S_{\kappa_j, \nu_j}^{1, l_j, m_j,\R})) f.
\end{equation*}

При $ d \in \N, l \in \Z_+^d, m \in \N^d, \kappa \in \Z_+^d $ определим 
линейный оператор $ \mathcal E_\kappa^{d,l,m,\R}: L^{\loc}(\R^d) \mapsto 
\mathcal P_\kappa^{d,l,m}, $ полагая для $ f \in L^{\loc}(\R^d) $ значение 

\begin{equation*} \tag{2.2.3}
(\mathcal E_\kappa^{d,l,m,\R} f)(x) 
= \sum_{ \nu \in \Nu_{-m, 2^\kappa -\e}^d} g_{\kappa, \nu}^{d,m}(x)
(U_{\kappa,\nu}^{d,l,m,\R} f)(x), \ x \in \R^d, \ 
\end{equation*}

где 
$ U_{\kappa,\nu}^{d,l,m,\R}: L^{\loc}(\R^d) \mapsto L_\infty(\R^d) $ -- 
линейный оператор, значение которого для $ f \in L^{\loc}(\R^d) $ определяется
равенством
\begin{multline*} \tag{2.2.4}
U_{\kappa,\nu}^{d,l,m,\R} f = \sum_{\epsilon \in \Upsilon^d:
\s(\epsilon) \subset \s(\kappa)} (-\e)^\epsilon \sum_{ \rho \in
\Rho_{\kappa,\nu,\epsilon}^{d,m}} (\prod_{j \in \s(\epsilon)}
a_{\nu_j -2\rho_j}^{m_j}) (\mathcal S_{\kappa -\epsilon,
\nu_{\kappa -\epsilon}^{d,m}(\rho)}^{d,l,m,\R} f), \\
d \in \N, l \in \Z_+^d, m \in \N^d, \kappa \in \Z_+^d, 
\nu \in \Nu_{-m, 2^\kappa -\e}^d \text{ ср. с } (1.3.21)).
\end{multline*}

Как и выше, при $ d \in \N, l \in \Z_+^d, m \in \N^d, \kappa \in \Z_+^d, 
\nu \in \Nu_{-m, 2^\kappa -\e}^d, 1 \le p \le \infty $ справедливы включения
\begin{equation*}
U_{\kappa,\nu}^{d,l,m,\R} \mid_{L_p(\R^d)} \in \mathcal B(L_p(\R^d), L_p(\R^d)), \\
\mathcal E_\kappa^{d,l,m,\R} \mid_{L_p(\R^d)} \in \mathcal B(L_p(\R^d), L_p(\R^d)).
\end{equation*}  

Обозначая через $ \mathcal I $ отображение, которое каждой функции 
$ f, $ заданной на $ I^d, $ сопоставляет функцию $ \mathcal I f, $ определяемую 
на $ \R^d $ равенством
\begin{equation*}
(\mathcal I f)(x) = \begin{cases} f(x), \text{ при } x \in I^d; \\
0, \text{ при } x \in \R^d \setminus I^d, \ 
\end{cases}
\end{equation*}
заметим, что для $ f \in L_1(I^d) $ при $ x \in I^d $ с учётом (2.2.1) 
выполняется равенство
\begin{multline*}
(\mathcal S_{\kappa, \nu}^{d,l,m,\R}(\mathcal I f))(x) = \chi_{Q^{d,m}}(x) 
(P_{2^{-\kappa}, 2^{-\kappa} \nu}^{d, l}((\mathcal I f) 
\mid_{2^{-\kappa} \nu +2^{-\kappa} I^d}))(x) = \\ 
(P_{2^{-\kappa}, 2^{-\kappa} \nu}^{d, l}(f \mid_{2^{-\kappa} \nu +2^{-\kappa} I^d}))(x) = 
(\mathcal S_{\kappa, \nu}^{d,l} f)(x),  \\
d \in \N, l \in \Z_+^d, m \in \N^d, \kappa \in \Z_+^d, \nu \in \Nu_{0, 2^\kappa -\e}^d.
\end{multline*}
  
Отсюда и из (2.2.4), (1.3.21) вытекает, что при $ d \in \N, l \in \Z_+^d, 
m \in \N^d, \kappa \in \Z_+^d, \nu \in \Nu_{-m, 2^\kappa -\e}^d $ 
для $ f \in L_1(I^d), x \in I^d $ соблюдается равенство
\begin{equation*}
(U_{\kappa,\nu}^{d,l,m,\R}(\mathcal I f))(x) = (U_{\kappa,\nu}^{d,l,m} f)(x), \ 
\end{equation*}
что в силу (2.2.3), (1.3.20) влечёт равенство
\begin{equation*} \tag{2.2.5}
(\mathcal E_\kappa^{d,l,m,\R}(\mathcal I f))(x) = (\mathcal E_\kappa^{d,l,m} f)(x).
\end{equation*}

Теперь установим справедливость следующего утверждения.

Лемма 2.2.1

Пусть $ d \in \N, l \in \Z_+^d, m \in \N^d, \kappa \in \Z_+^d, 
1 \le p < \infty. $ Тогда для $ f \in L_p(\R^d) $ справедливо включение
\begin{equation*} \tag{2.2.6}
(\mathcal E_\kappa^{d,l,m,\R} f) \in \mathcal P_\kappa^{\prime d,l,m}.
\end{equation*}

Доказательство.

В условиях леммы для $ f \in L_p(\R^d), $ принимая во внимание (2.2.1), (2.2.4), 
при $ \nu \in \Nu_{-m, 2^\kappa -\e}^d $ выберем полином 
$ \mathcal U_{\kappa,\nu}^{d,l,m} f \in \mathcal P^{d,l}, $ для которого при 
$ x \in Q^{d,m} $ соблюдается равенство
$ (\mathcal U_{\kappa,\nu}^{d,l,m} f)(x) =
(U_{\kappa,\nu}^{d,l,m,\R} f)(x). $  
Отсюда с учётом (2.2.3) и того факта, что $ \supp g_{\kappa,\nu}^{d,m} \subset 
Q^{d,m}, \ \nu \in \Nu_{-m, 2^\kappa -\e}^d, $ при $ j =1, \ldots, d $ 
для $ x \in \R^d $ имеем
\begin{multline*} \tag{2.2.7}
(\mathcal E_\kappa^{d,l,m,\R} f)(x) 
= \sum_{ \nu \in \Nu_{-m, 2^\kappa -\e}^d} g_{\kappa, \nu}^{d,m}(x)
(\mathcal U_{\kappa,\nu}^{d,l,m} f)(x) = \\
\sum_{ \nu_j \in \Nu_{-m_j, 2^{\kappa_j} -1}^1, \nu^{J_j} \in 
(\Nu_{-m, 2^\kappa -\e}^d)^{J_j}} 
g_{\kappa_j, \nu_j}^{1,m_j}(x_j)
g_{\kappa^{J_j}, \nu^{J_j}}^{d -1,m^{J_j}}(x^{J_j})
(\mathcal U_{\kappa,\nu}^{d,l,m} f)(x) = \\
\sum_{ \nu_j \in \Nu_{-m_j, 2^{\kappa_j} -1}^1}
\sum_{ \nu^{J_j} \in (\Nu_{-m, 2^\kappa -\e}^d)^{J_j}} 
g_{\kappa_j, \nu_j}^{1,m_j}(x_j)
g_{\kappa^{J_j}, \nu^{J_j}}^{d -1,m^{J_j}}(x^{J_j})
(\mathcal U_{\kappa,\nu}^{d,l,m} f)(x) = \\
\sum_{ \nu_j \in \Nu_{-m_j, 2^{\kappa_j} -1}^1}
g_{\kappa_j, \nu_j}^{1,m_j}(x_j)
(\sum_{ \nu^{J_j} \in (\Nu_{-m, 2^\kappa -\e}^d)^{J_j}} 
g_{\kappa^{J_j}, \nu^{J_j}}^{d -1,m^{J_j}}(x^{J_j})
(\mathcal U_{\kappa,\nu}^{d,l,m} f)(x)) = \\
\sum_{ \nu_j \in \Nu_{-m_j, 2^{\kappa_j} -1}^1}
g_{\kappa_j, \nu_j}^{1,m_j}(x_j) f_{\nu_j}^j(x), \  
\end{multline*}
где
\begin{equation*} \tag{2.2.8}
f_{\nu_j}^j(x) =
\sum_{ \nu^{J_j} \in (\Nu_{-m, 2^\kappa -\e}^d)^{J_j}} 
g_{\kappa^{J_j}, \nu^{J_j}}^{d -1,m^{J_j}}(x^{J_j})
(\mathcal U_{\kappa,\nu}^{d,l,m} f)(x), \ 
\nu_j \in \Nu_{-m_j, 2^{\kappa_j} -1}^1.
\end{equation*}

Сопоставляя (2.2.7), (2.2.8) с (2.1.1), (2.1.2), (2.1.3) и учитывая замечание 
перед леммой 2.1.1, для каждого $ j =1, \ldots, d $ проверим соблюдение условий 
(2.1.5), (2.1.6) при $ x \in Q^{d,m} $ для функций $ f_{\nu_j}^j, $ задаваемых 
равенствами (2.2.8).

Сначала рассмотрим случай, когда $ j \in \s(\kappa). $ В этом случае при 
$ \nu \in \Nu_{-m, 2^\kappa -\e}^d, x \in Q^{d,m}, $ используя (2.2.4), 
(2.2.2), (1.2.4), выводим
\begin{multline*} \tag{2.2.9}   
(\mathcal U_{\kappa,\nu}^{d,l,m} f)(x) = 
(U_{\kappa,\nu}^{d,l,m,\R} f)(x) = \\
\sum_{\epsilon \in \Upsilon^d:
\s(\epsilon) \subset \s(\kappa)} (-\e)^\epsilon \sum_{ \rho \in
\Rho_{\kappa,\nu,\epsilon}^{d,m}} (\prod_{i \in \s(\epsilon)}
a_{\nu_i -2\rho_i}^{m_i}) (\mathcal S_{\kappa -\epsilon,
\nu_{\kappa -\epsilon}^{d,m}(\rho)}^{d,l,m,\R} f)(x) = \\
\sum_{\substack{(\epsilon_j, \epsilon^{J_j}) \mid \epsilon_j \in \{0, 1\}, \\
\epsilon^{J_j} \in (\Upsilon^d)^{J_j}:
\s(\epsilon^{J_j}) \subset \s(\kappa^{J_j})}}
(-1)^{\epsilon_j} (-\e^{J_j})^{\epsilon^{J_j}} \sum_{\substack{ (\rho_j, 
\rho^{J_j}): \rho_j \in \Rho_{\kappa_j,\nu_j,\epsilon_j}^{1,m_j},\\
 \rho^{J_j} \in\Rho_{\kappa^{J_j},\nu^{J_j},\epsilon^{J_j}}^{d -1,m^{J_j}}}} 
(\prod_{i \in \s(\epsilon)}
a_{\nu_i -2\rho_i}^{m_i}) \\
\times\biggl(V_j(\mathcal S_{\kappa_j -\epsilon_j,
\nu_{\kappa_j -\epsilon_j}^{1,m_j}(\rho_j)}^{1,l_j,m_j,\R}) \biggl(\biggl(\prod_{i \in J_j} 
V_i(\mathcal S_{\kappa_i -\epsilon_i,
\nu_{\kappa_i -\epsilon_i}^{1,m_i}(\rho_i)}^{1,l_i,m_i,\R})\biggr) f\biggr)\biggr)(x) = \\
\sum_{\substack{\epsilon^{J_j} \in (\Upsilon^d)^{J_j}:\\
\s(\epsilon^{J_j}) \subset \s(\kappa^{J_j})}}
\sum_{\epsilon_j =0}^1 (-1)^{\epsilon_j} (-\e^{J_j})^{\epsilon^{J_j}} 
\sum_{ \rho_j \in \Rho_{\kappa_j,\nu_j,\epsilon_j}^{1,m_j}} 
\sum_{\rho^{J_j} \in \Rho_{\kappa^{J_j},\nu^{J_j},\epsilon^{J_j}}^{d -1,m^{J_j}}} 
(\prod_{i \in \s(\epsilon)}
a_{\nu_i -2\rho_i}^{m_i})\\
\times \biggl(V_j(\mathcal S_{\kappa_j -\epsilon_j,
\nu_{\kappa_j -\epsilon_j}^{1,m_j}(\rho_j)}^{1,l_j,m_j,\R}) \biggl(\biggl(\prod_{i \in J_j} 
V_i(\mathcal S_{\kappa_i -\epsilon_i,
\nu_{\kappa_i -\epsilon_i}^{1,m_i}(\rho_i)}^{1,l_i,m_i,\R})\biggr) f\biggr)\biggr)(x) = \\
\sum_{\substack{\epsilon^{J_j} \in (\Upsilon^d)^{J_j}:\\
\s(\epsilon^{J_j}) \subset \s(\kappa^{J_j})}}
(-\e^{J_j})^{\epsilon^{J_j}} 
\Biggl(\sum_{\rho^{J_j} \in \Rho_{\kappa^{J_j},\nu^{J_j},\epsilon^{J_j}}^{d -1,m^{J_j}}} 
(\prod_{i \in \s(\epsilon^{J_j})}
a_{\nu_i -2\rho_i}^{m_i})\\
\times \biggl(V_j(\mathcal S_{\kappa_j,
\nu_{\kappa_j}^{1,m_j}(\nu_j)}^{1,l_j,m_j,\R}) \biggl(\biggl(\prod_{i \in J_j} 
V_i(\mathcal S_{\kappa_i -\epsilon_i,
\nu_{\kappa_i -\epsilon_i}^{1,m_i}(\rho_i)}^{1,l_i,m_i,\R})\biggr) f\biggr)\biggr)(x) - \\
\sum_{ \rho_j \in \Rho_{\kappa_j,\nu_j,1}^{1,m_j}} a_{\nu_j -2\rho_j}^{m_j} 
\sum_{\rho^{J_j} \in \Rho_{\kappa^{J_j},\nu^{J_j},\epsilon^{J_j}}^{d -1,m^{J_j}}} 
(\prod_{i \in \s(\epsilon^{J_j})}
a_{\nu_i -2\rho_i}^{m_i}) \\
\times\biggl(V_j(\mathcal S_{\kappa_j -1, 
\nu_{\kappa_j -1}^{1,m_j}(\rho_j)}^{1,l_j,m_j,\R}) \biggl(\biggl(\prod_{i \in J_j} 
V_i(\mathcal S_{\kappa_i -\epsilon_i,
\nu_{\kappa_i -\epsilon_i}^{1,m_i}(\rho_i)}^{1,l_i,m_i,\R})\biggr) f\biggr)\biggr)(x)\Biggr) = \\
\sum_{\epsilon^{J_j} \in (\Upsilon^d)^{J_j}:
\s(\epsilon^{J_j}) \subset \s(\kappa^{J_j})}
(-\e^{J_j})^{\epsilon^{J_j}} 
\Biggl(\Biggl(V_j(\mathcal S_{\kappa_j, \nu_{\kappa_j}^{1,m_j}(\nu_j)}^{1,l_j,m_j,\R}) \\
\times\biggl(\sum_{\rho^{J_j} \in \Rho_{\kappa^{J_j},\nu^{J_j},\epsilon^{J_j}}^{d -1,m^{J_j}}} 
(\prod_{i \in \s(\epsilon^{J_j})}
a_{\nu_i -2\rho_i}^{m_i}) (\prod_{i \in J_j} 
V_i(\mathcal S_{\kappa_i -\epsilon_i,
\nu_{\kappa_i -\epsilon_i}^{1,m_i}(\rho_i)}^{1,l_i,m_i,\R})) f\biggr)\Biggr)(x) - \\
\sum_{ \rho_j \in \Rho_{\kappa_j,\nu_j,1}^{1,m_j}} a_{\nu_j -2\rho_j}^{m_j} 
\Biggl(V_j(\mathcal S_{\kappa_j -1, \nu_{\kappa_j -1}^{1,m_j}(\rho_j)}^{1,l_j,m_j,\R})\\ 
\times\biggl(\sum_{\rho^{J_j} \in \Rho_{\kappa^{J_j},\nu^{J_j},\epsilon^{J_j}}^{d -1,m^{J_j}}} 
(\prod_{i \in \s(\epsilon^{J_j})}
a_{\nu_i -2\rho_i}^{m_i}) (\prod_{i \in J_j} 
V_i(\mathcal S_{\kappa_i -\epsilon_i,
\nu_{\kappa_i -\epsilon_i}^{1,m_i}(\rho_i)}^{1,l_i,m_i,\R})) f\biggr)\Biggr)(x)\Biggr).
\end{multline*} 

Замечая, что при $ \nu_j \in \Nu_{-m_j,0}^1 $ значение 
\begin{equation*}
\nu_{\kappa_j}^{1,m_j}(\nu_j) = (2^{\kappa_j} -m_j -1)_+ -(2^{\kappa_j} -m_j -1 -\nu_{+ j})_+ =
(2^{\kappa_j} -m_j -1)_+ -(2^{\kappa_j} -m_j -1 -0)_+ =0, \ 
\end{equation*}
а при $ \nu_j \in \Nu_{-m_j,0}^1, \rho_j \in \Rho_{\kappa_j,\nu_j,1}^{1,m_j} $ 
имеют место соотношения
\begin{equation*}
-m_j \le \rho_j = (\nu_j -\mu_j) /2 \le (0 -0) /2 =0, 
\text{где} \mu_j \in \Nu_{0, m_j +1}^1,
\end{equation*}
и, значит,
$$
\nu_{\kappa_j -1}^{1,m_j}(\rho_j) = (2^{\kappa_j -1} -m_j -1)_+ -(2^{\kappa_j -1} -m_j -1 -\rho_{+ j})_+ =0,
$$
заключаем, что для $ \nu \in \Nu_{-m, 2^\kappa -\e}^d: \nu_j \in \Nu_{-m_j,0}^1, $ 
в соответствии с (2.2.9), (1.3.5) имеет место равенство
\begin{multline*} \tag{2.2.10}
(\mathcal U_{\kappa,\nu}^{d,l,m} f)(x) = 
\sum_{\epsilon^{J_j} \in (\Upsilon^d)^{J_j}:
\s(\epsilon^{J_j}) \subset \s(\kappa^{J_j})}
(-\e^{J_j})^{\epsilon^{J_j}} 
\Biggl(\Biggl(V_j(\mathcal S_{\kappa_j, 0}^{1,l_j,m_j,\R}) \\
\times\biggl(\sum_{\rho^{J_j} \in \Rho_{\kappa^{J_j},\nu^{J_j},\epsilon^{J_j}}^{d -1,m^{J_j}}} 
(\prod_{i \in \s(\epsilon^{J_j})}
a_{\nu_i -2\rho_i}^{m_i}) (\prod_{i \in J_j} 
V_i(\mathcal S_{\kappa_i -\epsilon_i,
\nu_{\kappa_i -\epsilon_i}^{1,m_i}(\rho_i)}^{1,l_i,m_i,\R})) f\biggr)\Biggr)(x) - \\
\biggl(\sum_{ \rho_j \in \Rho_{\kappa_j,\nu_j,1}^{1,m_j}} a_{\nu_j -2\rho_j}^{m_j}\biggr) 
\biggl(V_j(\mathcal S_{\kappa_j -1, 0}^{1,l_j,m_j,\R}) \\
\times\biggl(\sum_{\rho^{J_j} \in \Rho_{\kappa^{J_j},\nu^{J_j},\epsilon^{J_j}}^{d -1,m^{J_j}}} 
(\prod_{i \in \s(\epsilon^{J_j})}
a_{\nu_i -2\rho_i}^{m_i}) (\prod_{i \in J_j} 
V_i(\mathcal S_{\kappa_i -\epsilon_i,
\nu_{\kappa_i -\epsilon_i}^{1,m_i}(\rho_i)}^{1,l_i,m_i,\R})) f\biggr)\biggr)(x)\Biggr) = \\
\sum_{\epsilon^{J_j} \in (\Upsilon^d)^{J_j}:
\s(\epsilon^{J_j}) \subset \s(\kappa^{J_j})}
(-\e^{J_j})^{\epsilon^{J_j}} 
\Biggl((V_j(\mathcal S_{\kappa_j, 0}^{1,l_j,m_j,\R}) -
V_j(\mathcal S_{\kappa_j -1, 0}^{1,l_j,m_j,\R})) \\
\times\biggl(\sum_{\rho^{J_j} \in \Rho_{\kappa^{J_j},\nu^{J_j},\epsilon^{J_j}}^{d -1,m^{J_j}}} 
(\prod_{i \in \s(\epsilon^{J_j})}
a_{\nu_i -2\rho_i}^{m_i}) (\prod_{i \in J_j} 
V_i(\mathcal S_{\kappa_i -\epsilon_i,
\nu_{\kappa_i -\epsilon_i}^{1,m_i}(\rho_i)}^{1,l_i,m_i,\R})) f\biggr)\Biggr)(x), \ x \in Q^{d,m}.  
\end{multline*}
Соединяя (2.2.8) и (2.2.10), видим, что соблюдается (2.1.5) при 
$ x \in Q^{d,m}, j \in \s(\kappa). $

Кроме того, принимая во внимание, что при $ \nu_j \in \Nu_{2^{\kappa_j} -m_j -1, 
2^{\kappa_j} -1}^1 $ значение 
\begin{equation*}
\nu_{\kappa_j}^{1,m_j}(\nu_j) = (2^{\kappa_j} -m_j -1)_+ -(2^{\kappa_j} -m_j -1 -\nu_{+ j})_+ =
(2^{\kappa_j} -m_j -1)_+, \  
\end{equation*}
а при $ \nu_j \in \Nu_{2^{\kappa_j} -m_j -1, 2^{\kappa_j} -1}^1,  
\rho_j \in \Rho_{\kappa_j,\nu_j,1}^{1,m_j}, $ справедливо 
соотношение
\begin{equation*}
2^{\kappa_j -1} -1 \ge \rho_j = (\nu_j -\mu_j) /2 \ge (2^{\kappa_j} -m_j -1 
-(m_j +1)) /2 = 2^{\kappa_j -1} -m_j -1, \\
\text{где } \mu_j \in \Nu_{0, m_j +1}^1, \ 
\end{equation*}
и, следовательно,
$$
\nu_{\kappa_j -1}^{1,m_j}(\rho_j) = (2^{\kappa_j -1} -m_j -1)_+ -(2^{\kappa_j -1} -m_j -1 -\rho_{+ j})_+ = 
(2^{\kappa_j -1} -m_j -1)_+, \ 
$$
получаем, что для $ \nu \in \Nu_{-m, 2^\kappa -\e}^d: \nu_j \in 
\Nu_{2^{\kappa_j} -m_j -1, 2^{\kappa_j} -1}^1, $ 
согласно (2.2.9), (1.3.5), соблюдается равенство
\begin{multline*} \tag{2.2.11}
(\mathcal U_{\kappa,\nu}^{d,l,m} f)(x) = 
\sum_{\epsilon^{J_j} \in (\Upsilon^d)^{J_j}:
\s(\epsilon^{J_j}) \subset \s(\kappa^{J_j})}
(-\e^{J_j})^{\epsilon^{J_j}} 
\Biggl(\biggl(V_j(\mathcal S_{\kappa_j, (2^{\kappa_j} -m_j -1)_+}^{1,l_j,m_j,\R}) \\
\times\biggl(\sum_{\rho^{J_j} \in \Rho_{\kappa^{J_j},\nu^{J_j},\epsilon^{J_j}}^{d -1,m^{J_j}}} 
(\prod_{i \in \s(\epsilon^{J_j})}
a_{\nu_i -2\rho_i}^{m_i}) (\prod_{i \in J_j} 
V_i(\mathcal S_{\kappa_i -\epsilon_i,
\nu_{\kappa_i -\epsilon_i}^{1,m_i}(\rho_i)}^{1,l_i,m_i,\R})) f\biggr)\biggr)(x) - \\
\biggl(\sum_{ \rho_j \in \Rho_{\kappa_j,\nu_j,1}^{1,m_j}} a_{\nu_j -2\rho_j}^{m_j}\biggr) 
\biggl(V_j(\mathcal S_{\kappa_j -1, (2^{\kappa_j -1} -m_j -1)_+}^{1,l_j,m_j,\R}) \\
\times\biggl(\sum_{\rho^{J_j} \in \Rho_{\kappa^{J_j},\nu^{J_j},\epsilon^{J_j}}^{d -1,m^{J_j}}} 
(\prod_{i \in \s(\epsilon^{J_j})}
a_{\nu_i -2\rho_i}^{m_i}) (\prod_{i \in J_j} 
V_i(\mathcal S_{\kappa_i -\epsilon_i,
\nu_{\kappa_i -\epsilon_i}^{1,m_i}(\rho_i)}^{1,l_i,m_i,\R})) f\biggr)\biggr)(x)\Biggr) = \\
\sum_{\epsilon^{J_j} \in (\Upsilon^d)^{J_j}:
\s(\epsilon^{J_j}) \subset \s(\kappa^{J_j})}
(-\e^{J_j})^{\epsilon^{J_j}} 
\Biggl(\biggl(V_j(\mathcal S_{\kappa_j, (2^{\kappa_j} -m_j -1)_+}^{1,l_j,m_j,\R}) -
V_j(\mathcal S_{\kappa_j -1, (2^{\kappa_j -1} -m_j -1)_+}^{1,l_j,m_j,\R})\biggr) \\
\times\biggl(\sum_{\rho^{J_j} \in \Rho_{\kappa^{J_j},\nu^{J_j},\epsilon^{J_j}}^{d -1,m^{J_j}}} 
(\prod_{i \in \s(\epsilon^{J_j})}
a_{\nu_i -2\rho_i}^{m_i}) (\prod_{i \in J_j} 
V_i(\mathcal S_{\kappa_i -\epsilon_i,
\nu_{\kappa_i -\epsilon_i}^{1,m_i}(\rho_i)}^{1,l_i,m_i,\R})) f\biggr)\Biggr)(x), \ x \in Q^{d,m}.  
\end{multline*}

Соединяя (2.2.8) с (2.2.11), приходим к выводу, что имеет место (2.1.6) при 
$ x \in Q^{d,m}, j \in \s(\kappa). $

Теперь установим соблюдение условий (2.1.5), (2.1.6) при $ x \in Q^{d,m} $
для функций $ f_{\nu_j}^j, $ задаваемых равенствами (2.2.8), в случае, когда 
$ j \in \Nu_{1, d}^1 \setminus \s(\kappa), $ т.е. $ \kappa_j =0. $ 
В этом случае $ \Nu_{-m_j, 2^{\kappa_j} -1}^1 = \Nu_{-m_j, 0}^1 = 
\Nu_{2^{\kappa_j} -m_j -1, 2^{\kappa_j} -1}^1, $ и при 
$ \nu \in \Nu_{-m, 2^\kappa -\e}^d, x \in Q^{d,m}, $ используя (2.2.4), 
(2.2.2), (1.2.4) и учитывая, что $ \Rho_{\kappa_j,\nu_j,\epsilon_j}^{1,m_j} =
\Rho_{0,\nu_j,0}^{1,m_j} = \{\nu_j\}, \nu_0^{1,m_j}(\nu_j) = 
(1 -m_j -1)_+ -(1 -m_j -1 -\nu_{+ j})_+ =0, \ \nu_j \in \Nu_{-m_j, 2^{\kappa_j} -1}^1, \ $ 
имеем
\begin{multline*} \tag{2.2.12}
(\mathcal U_{\kappa,\nu}^{d,l,m} f)(x) = 
(U_{\kappa,\nu}^{d,l,m,\R} f)(x) = \\ 
\sum_{\epsilon \in \Upsilon^d:
\s(\epsilon) \subset \s(\kappa)} (-\e)^\epsilon \sum_{ \rho \in
\Rho_{\kappa,\nu,\epsilon}^{d,m}} (\prod_{i \in \s(\epsilon)}
a_{\nu_i -2\rho_i}^{m_i}) (\mathcal S_{\kappa -\epsilon,
\nu_{\kappa -\epsilon}^{d,m}(\rho)}^{d,l,m,\R} f)(x) = \\
\sum_{\substack{(\epsilon_j, \epsilon^{J_j}) \mid \epsilon_j =0, \\
\epsilon^{J_j} \in (\Upsilon^d)^{J_j}:
\s(\epsilon^{J_j}) \subset \s(\kappa^{J_j})}}
(-1)^{\epsilon_j} 
(-\e^{J_j})^{\epsilon^{J_j}} \sum_{\substack{ (\rho_j, \rho^{J_j}): 
\rho_j \in \Rho_{\kappa_j,\nu_j,\epsilon_j}^{1,m_j},\\ 
\rho^{J_j} \in\Rho_{\kappa^{J_j},\nu^{J_j},\epsilon^{J_j}}^{d -1,m^{J_j}}}} 
(\prod_{i \in \s(\epsilon)}
a_{\nu_i -2\rho_i}^{m_i}) \\
\times\biggl(V_j(\mathcal S_{\kappa_j -\epsilon_j,
\nu_{\kappa_j -\epsilon_j}^{1,m_j}(\rho_j)}^{1,l_j,m_j,\R}) ((\prod_{i \in J_j} 
V_i(\mathcal S_{\kappa_i -\epsilon_i,
\nu_{\kappa_i -\epsilon_i}^{1,m_i}(\rho_i)}^{1,l_i,m_i,\R})) f)\biggr)(x) = \\
\sum_{\epsilon^{J_j} \in (\Upsilon^d)^{J_j}:
\s(\epsilon^{J_j}) \subset \s(\kappa^{J_j})}
(-\e^{J_j})^{\epsilon^{J_j}} 
\sum_{\rho^{J_j} \in \Rho_{\kappa^{J_j},\nu^{J_j},\epsilon^{J_j}}^{d -1,m^{J_j}}} 
(\prod_{i \in \s(\epsilon^{J_j})}
a_{\nu_i -2\rho_i}^{m_i}) \\
\times\biggl(V_j(\mathcal S_{\kappa_j,
\nu_{\kappa_j}^{1,m_j}(\nu_j)}^{1,l_j,m_j,\R}) ((\prod_{i \in J_j} 
V_i(\mathcal S_{\kappa_i -\epsilon_i,
\nu_{\kappa_i -\epsilon_i}^{1,m_i}(\rho_i)}^{1,l_i,m_i,\R})) f)\biggr)(x) = \\
\sum_{\epsilon^{J_j} \in (\Upsilon^d)^{J_j}:
\s(\epsilon^{J_j}) \subset \s(\kappa^{J_j})}
(-\e^{J_j})^{\epsilon^{J_j}} 
\Biggl(V_j(\mathcal S_{0, \nu_0^{1,m_j}(\nu_j)}^{1,l_j,m_j,\R}) 
\biggl(\sum_{\rho^{J_j} \in \Rho_{\kappa^{J_j},\nu^{J_j},\epsilon^{J_j}}^{d -1,m^{J_j}}} 
(\prod_{i \in \s(\epsilon^{J_j})}
a_{\nu_i -2\rho_i}^{m_i})\\
\times (\prod_{i \in J_j} 
V_i(\mathcal S_{\kappa_i -\epsilon_i,
\nu_{\kappa_i -\epsilon_i}^{1,m_i}(\rho_i)}^{1,l_i,m_i,\R})) f\biggr)\Biggr)(x) = \\
\sum_{\epsilon^{J_j} \in (\Upsilon^d)^{J_j}:
\s(\epsilon^{J_j}) \subset \s(\kappa^{J_j})}
(-\e^{J_j})^{\epsilon^{J_j}} 
\Biggl(V_j(\mathcal S_{0, 0}^{1,l_j,m_j,\R}) 
\biggl(\sum_{\rho^{J_j} \in \Rho_{\kappa^{J_j},\nu^{J_j},\epsilon^{J_j}}^{d -1,m^{J_j}}} 
(\prod_{i \in \s(\epsilon^{J_j})}
a_{\nu_i -2\rho_i}^{m_i}) \\
\times(\prod_{i \in J_j} 
V_i(\mathcal S_{\kappa_i -\epsilon_i,
\nu_{\kappa_i -\epsilon_i}^{1,m_i}(\rho_i)}^{1,l_i,m_i,\R})) f\biggr)\Biggr)(x).
\end{multline*} 

Объединяя (2.2.8) и (2.2.12), заключаем, что для рассматриваемых функций условия 
(2.1.5) и (2.1.6) выполняются при $ x \in Q^{d,m} $ в случае 
$ j \in \Nu_{1, d}^1 \setminus \s(\kappa). \square $

Предложение 2.2.2

Пусть $ d \in \N, l \in \N^d, m \in \N^d, \lambda \in \Z_+^d(m), 
1 \le p < \infty, 1 \le q \le \infty. $ 
Тогда существуют константы $ c_1(d,l,m,\lambda,p,q) > 0, c_2(d,m) > 0 $ 
такие, что при $ \kappa \in \Z_+^d $ для $ f \in L_p(I^d) $ имеет место 
неравенство
\begin{multline*} \tag{2.2.13}
\| \D^\lambda \mathcal E_\kappa^{d,l -\e,m,\R} (\mathcal I f) \|_{L_q(\R^d)} \le
c_1 2^{(\kappa, \lambda +(p^{-1} -q^{-1})_+ \e)} \\
\Omega^{\prime l \chi_{\s(\kappa)}}(f, (c_2 2^{-\kappa})^{\s(\kappa)})_{L_p(I^d)}.
\end{multline*}

Доказательство.

На основании (2.2.6) применяя (2.1.15), а затем используя (2.2.5), (1.3.43), 
получаем, что в условиях предложения для $ f \in L_p(I^d) $ при 
$ \kappa \in \Z_+^d $ справедливо неравенство
\begin{multline*}
\| \D^\lambda \mathcal E_\kappa^{d,l -\e,m,\R} (\mathcal I f) \|_{L_q(\R^d)} \le 
c_3 2^{(\kappa, \lambda)} \| \mathcal E_\kappa^{d,l -\e,m,\R} (\mathcal I f) \|_{L_q(I^d)} = \\
c_3 2^{(\kappa, \lambda)} \| \mathcal E_\kappa^{d,l -\e,m} (f) \|_{L_q(I^d)} \le
c_3 2^{(\kappa, \lambda)} c_4 2^{(\kappa, (p^{-1} -q^{-1})_+ \e)} 
\Omega^{\prime l \chi_{\s(\kappa)}}(f, (c_2 2^{-\kappa})^{\s(\kappa)})_{L_p(I^d)} = \\
c_1 2^{(\kappa, \lambda +(p^{-1} -q^{-1})_+ \e)} 
\Omega^{\prime l \chi_{\s(\kappa)}}(f, (c_2 2^{-\kappa})^{\s(\kappa)})_{L_p(I^d)}. \square
\end{multline*}
\bigskip

2.3. В этом пункте доказывается основной результат работы

Теорема 2.3.1

Пусть $ d \in \N, \alpha \in \R_+^d, l = l(\alpha), 1 \le p < \infty, 
1 \le \theta \le \infty. $ Тогда существует непрерывное линейное отображение 
$ \mathcal E^{d,\alpha,p,\theta}: (S_{p,\theta}^\alpha B)^\prime(I^d) \mapsto L_p(\R^d), $  
обладающее следующими свойствами: 

1) для $ f \in (S_{p,\theta}^\alpha B)^\prime(I^d) $ выполняется равенство
\begin{equation*} \tag{2.3.1}
f = (\mathcal E^{d,\alpha,p,\theta} f) \mid_{I^d},
\end{equation*} 

2) существует константа $ C(d,\alpha,p,\theta) >0 $ такая, что 
при $ \lambda \in \Z_+^d: \lambda < \alpha, $ для $ f \in (S_{p,\theta}^\alpha B)^\prime(I^d) $ 
и любого множества $ J \subset \{1,\ldots,d\} $ верно неравенство
\begin{multline*} \tag{2.3.2}
\biggl(\int_{(\R_+^d)^J} (t^J)^{-\e^J -\theta (\alpha^J -\lambda^J)} 
(\Omega^{(l -\lambda) \chi_J}(\D^\lambda (\mathcal E^{d,\alpha,p,\theta} f), 
t^J)_{L_p(\R^d)})^\theta dt^J\biggr)^{1/\theta} \le \\
C \| f\|_{(S_{p,\theta}^\alpha B)^\prime(I^d)}, \ \theta \ne \infty; \\
\sup_{t^J \in (\R_+^d)^J} (t^J)^{-(\alpha^J -\lambda^J)} 
\Omega^{(l -\lambda) \chi_J}(\D^\lambda (\mathcal E^{d,\alpha,p,\theta} f), 
t^J)_{L_p(\R^d)} \le\\
 C \| f\|_{(S_{p,\theta}^\alpha B)^\prime(I^d)}, \ \theta = \infty.
\end{multline*}

Доказательство.

Проведём доказательство при $ \theta \ne \infty. $
В условиях теоремы фиксируем $ m \in \N^d $ так, чтобы $ l \in \Z_+^d(m). $ 
Ориентируясь на выкладку перед (1.1.8), так же, как (2.2.7) из [12], выводится 
используемое ниже соотношение 
\begin{multline*} \tag{2.3.3}
\sum_{\kappa \in \Z_+^d \setminus \{0\}} (2^{(\kappa, \alpha)} 
\Omega^{\prime l \chi_{\s(\kappa)}}(f,
c_0 (2^{-\kappa})^{\s(\kappa)})_{L_p(I^d)})^\theta \le\\
c_1(d,\alpha,p,\theta,c_0) \sum_{ J \subset \Nu_{1,d}^1: J \ne \emptyset} 
\int_{(\R_+^d)^J} (t^J)^{-\e^J -\theta \alpha^J} (\Omega^{\prime l \chi_J}(f,
t^J)_{L_p(I^d)})^\theta dt^J, \\ 
f \in (S_{p, \theta}^\alpha B)^\prime(I^d), c_0 \in \R_+. 
\end{multline*}

Для построения оператора $ \mathcal  E^{d, \alpha, p, \theta} $ с требуемыми
свойствами заметим, что для $ f \in (S_{p, \theta}^\alpha B)^\prime(I^d) $ при 
$ \lambda \in \Z_+^d: \lambda < \alpha, $ ряд 
\begin{equation*}
\sum_{\kappa \in \Z_+^d} \| \D^\lambda (\mathcal E_\kappa^{d,l -\e,m,\R} 
(\mathcal I f))\|_{L_p(\R^d)}
\end{equation*}
сходится, ибо в силу (2.2.13), неравенства Гёльдера, а также 
соотношений (2.3.3) и (1.1.1) при $ a = 1/\theta $ имеет место неравенство
\begin{multline*} \tag{2.3.4} 
\sum_{ \kappa \in \Z_+^d: n \le (\kappa, \e) \le n +k} 
\| \D^\lambda (\mathcal E_\kappa^{d,l -\e,m,\R}(\mathcal I f)) \|_{L_p(\R^d)} \le \\
\sum_{ \kappa \in \Z_+^d: n \le (\kappa, \e) \le n +k} c_2 2^{(\kappa, \lambda)} 
\Omega^{\prime l \chi_{\s(\kappa)}}(f, (c_3 2^{-\kappa})^{\s(\kappa)})_{L_p(I^d)} = \\
\sum_{ \kappa \in \Z_+^d: n \le (\kappa, \e) \le n +k} 2^{-(\kappa, \alpha -\lambda)} 
2^{(\kappa, \alpha)} c_2 
\Omega^{\prime l \chi_{\s(\kappa)}}(f, (c_3 2^{-\kappa})^{\s(\kappa)})_{L_p(I^d)} \le\\
 c_2 \biggl( \sum_{ \kappa \in \Z_+^d: n \le (\kappa, \e) \le n +k}
2^{-\theta^\prime (\kappa, \alpha -\lambda)}\biggr)^{1 / \theta^\prime} 
\biggl(\sum_{\kappa \in \Z_+^d: n \le (\kappa, \e) \le n +k} (2^{(\kappa, \alpha)} 
\Omega^{\prime l \chi_{\s(\kappa)}}(f, (c_3 2^{-\kappa})^{\s(\kappa)})_{L_p(I^d)})^\theta \biggr)^{1 / \theta} \le \\
c_2 \biggl( \sum_{ \kappa \in \Z_+^d: n \le (\kappa, \e) \le n +k}
2^{-\theta^\prime (\kappa, \alpha -\lambda)}\biggr)^{1 / \theta^\prime} 
\biggl(\sum_{\kappa \in \Z_+^d \setminus \{0\}} (2^{(\kappa, \alpha)} 
\Omega^{\prime l \chi_{\s(\kappa)}}(f, (c_3 2^{-\kappa})^{\s(\kappa)})_{L_p(I^d)})^\theta \biggr)^{1 / \theta} \le \\
c_4 \biggl( \sum_{ \kappa \in \Z_+^d: n \le (\kappa, \e) \le n +k}
2^{-\theta^\prime (\kappa, \alpha -\lambda)}\biggr)^{1 / \theta^\prime} \\
\sum_{ J \subset \Nu_{1,d}^1: J \ne \emptyset} 
\biggl(\int_{(\R_+^d)^J} (t^J)^{-\e^J -\theta \alpha^J} (\Omega^{\prime l \chi_J}(f,
t^J)_{L_p(I^d)})^\theta dt^J\biggr)^{1 /\theta}, \ n,k \in \N, \ 
\end{multline*}
где $ \theta^\prime = \theta/(\theta -1), $ а для семейства чисел 
$ \{a_\kappa \in \R: a_\kappa \ge 0, \kappa \in \Z_+^d\} $ 
сходимость ряда $ \sum_{\kappa \in \Z_+^d} a_\kappa $ эквивалентна 
сходимости последовательности $ \{\sum_{\kappa \in \Z_+^d: (\kappa, \e) \le n} 
a_\kappa, \ n \in \Z_+\}. $
Учитывая замечание после леммы 1.3.1, видим, что для $ f \in 
(S_{p,\theta}^\alpha B)^\prime(I^d) $ при $ \lambda \in \Z_+^d: 
\lambda < \alpha, $ ряд 
$ \sum_{ \kappa \in \Z_+^d} \D^\lambda (\mathcal E_\kappa^{d,l -\e,m,\R}(\mathcal I f)) $ 
сходится в $ L_p(\R^d). $ 
Определим отображение $ \mathcal E^{d,\alpha,p,\theta}: (S_{p,\theta}^\alpha B)^\prime(I^d) 
\mapsto L_p(\R^d), $ полагая для $ f \in (S_{p,\theta}^\alpha B)^\prime(I^d) $ значение
\begin{equation*} \tag{2.3.5}
\mathcal E^{d,\alpha,p,\theta} f = \sum_{ \kappa \in \Z_+^d} 
\mathcal E_\kappa^{d,l -\e,m,\R}(\mathcal I f).
\end{equation*}

Исходя из (2.3.5), с учётом сказанного выше заметим, что для 
$ f \in (S_{p,\theta}^\alpha B)^\prime(I^d) $ 
при $ \lambda \in \Z_+^d: \lambda < \alpha, $ в $ L_p(\R^d) $ имеет 
место равенство
\begin{equation*} \tag{2.3.6}
\D^\lambda (\mathcal E^{d,\alpha,p,\theta} f) = \sum_{ \kappa \in \Z_+^d} 
\D^\lambda (\mathcal E_\kappa^{d,l -\e,m,\R}(\mathcal I f)).
\end{equation*}
Ясно, что отображение $ \mathcal E^{d,\alpha,p,\theta} $ линейно, а вследствие 
(2.3.5), (2.2.5), (1.3.45) (см. предложение 1.3.6 и (1.1.8)) в 
$ L_p(I^d) $ справедливо равенство
\begin{equation*} 
(\mathcal E^{d,\alpha,p,\theta} f) \mid_{I^d} = \sum_{ \kappa \in \Z_+^d} 
(\mathcal E_\kappa^{d,l -\e,m,\R}(\mathcal I f)) \mid_{I^d} = 
\sum_{ \kappa \in \Z_+^d} \mathcal E_\kappa^{d, l-\e,m} f = f, \ 
\end{equation*}
т.е. соблюдается (2.3.1).

При проверке справедливости п. 2 теоремы, прежде всего, отметим, что 
для $ f \in (S_{p,\theta}^\alpha B)^\prime(I^d) $ 
при $ \lambda \in \Z_+^d: \lambda < \alpha, $ ввиду (2.3.6), (2.3.4) 
выполняется неравенство
\begin{multline*} \tag{2.3.7}
\| \D^\lambda (\mathcal E^{d,\alpha,p,\theta} f) \|_{L_p(\R^d)} =
\| \sum_{ \kappa \in \Z_+^d } 
\D^\lambda (\mathcal E_\kappa^{d,l -\e,m,\R}(\mathcal I f)) \|_{L_p(\R^d)} \le \\
\sum_{ \kappa \in \Z_+^d} 
\| \D^\lambda (\mathcal E_\kappa^{d,l -\e,m,\R}(\mathcal I f)) \|_{L_p(\R^d)} \le \\
c_5 (\| f \|_{L_p(I^d)} +\sum_{ J \subset \Nu_{1,d}^1: J \ne \emptyset} 
\biggl(\int_{(\R_+^d)^J} (t^J)^{-\e^J -\theta \alpha^J} (\Omega^{\prime l \chi_J}(f,
t^J)_{L_p(I^d)})^\theta dt^J\biggr)^{1 /\theta}).
\end{multline*}
 
Далее, при $ \lambda \in \Z_+^d: \lambda < \alpha, $ для функции $ F =
\mathcal E^{d,\alpha,p,\theta} f, \ f \in (S_{p,\theta}^\alpha B)^\prime(I^d) $ 
и любого непустого множества $ J \subset \{1,\ldots,d\} $ оценим \\
$ (\int_{(\R_+^d)^J} (t^J)^{-\e^J -\theta (\alpha^J -\lambda^J)} 
(\Omega^{(l -\lambda) \chi_J}(\D^\lambda F, 
t^J)_{L_p(\R^d)})^\theta dt^J)^{1/\theta}. $ Поскольку
для любого непустого множества $ J \subset \{1,\ldots,d\} $
справедливо представление $ (\R_+^d)^J = \cup_{ \J \subset J}
(I^d)^\J \times ([1, \infty)^d)^{J \setminus \J}, $ причём
множества в правой части последнего равенства попарно не пересекаются, то
\begin{multline*} \tag{2.3.8}
\int_{(\R_+^d)^J} (t^J)^{-\e^J -\theta (\alpha^J -\lambda^J)}
(\Omega^{(l -\lambda) \chi_J}(\D^\lambda F, t^J)_{L_p(\R^d)})^\theta dt^J =\\
\sum_{ \J \subset J} \int_{ (I^d)^{\J} \times ([1, \infty)^d)^{J
\setminus \J}} (t^J)^{-\e^J -\theta (\alpha^J -\lambda^J)} 
(\Omega^{(l -\lambda) \chi_J}(\D^\lambda F, t^J)_{L_p(\R^d)})^\theta dt^J.
\end{multline*}

Учитывая, что для $  J \subset \{1,\ldots,d\}: J \ne \emptyset, \J
\subset J $ при $ t^J \in (\R_+^d)^J $ выполняется неравенство
\begin{multline*}
\Omega^{(l -\lambda) \chi_J} (\D^\lambda F, t^J)_{L_p(\R^d)} = 
\supvrai_{ \{ h \in \R^d: h^J \in t^J (B^d)^J \}} \| \Delta_h^{(l -\lambda) \chi_J}
\D^\lambda F \|_{L_p(\R^d)} =\\
\supvrai_{ \{ h \in \R^d: h^J \in t^J (B^d)^J \}} 
\| \Delta_h^{(l -\lambda) \chi_{J \setminus \J}} (\Delta_h^{(l -\lambda) \chi_{\J}}
\D^\lambda F) \|_{L_p(\R^d)} \le \\
\supvrai_{ \{ h \in \R^d: h^J \in t^J (B^d)^J \}} \!c_6
\| \Delta_h^{(l -\lambda) \chi_{\J}} \D^\lambda F \|_{L_p(\R^d)}\! = c_6
\!\supvrai_{ \{ h \in \R^d: h^J \in t^J (B^d)^J \}} \!
\| \Delta_h^{(l -\lambda) \chi_{\J}}
\D^\lambda F \|_{L_p(\R^d)} \le\\
 c_6 \supvrai_{ \{ h \in \R^d: h^{\J} \in t^{\J} (B^d)^{\J} \}} 
\| \Delta_h^{(l -\lambda) \chi_{\J}} \D^\lambda F \|_{L_p(\R^d)} = c_6
\Omega^{(l -\lambda) \chi_{\J}} (\D^\lambda F, t^{\J})_{L_p(\R^d)}, 
\end{multline*}
находим, что
\begin{multline*}
\int_{ (I^d)^{\J} \times ([1, \infty)^d)^{J \setminus \J}}
(t^J)^{-\e^J -\theta (\alpha^J -\lambda^J)}
(\Omega^{(l -\lambda) \chi_J}(\D^\lambda F, t^J)_{L_p(\R^d)})^\theta dt^J \le\\
c_6^\theta \int_{ (I^d)^{\J} \times ([1, \infty)^d)^{J \setminus \J}} 
(t^J)^{-\e^J -\theta (\alpha^J -\lambda^J)}
(\Omega^{(l -\lambda) \chi_{\J}}(\D^\lambda F, 
t^{\J})_{L_p(\R^d)})^\theta dt^J =\\
c_6^\theta \int_{ (I^d)^{\J} \times ([1, \infty)^d)^{J \setminus \J}} 
(t^{\J}\!)^{\!-\! \e^{\J} \!-\! \theta (\alpha^{\J} -\lambda^{\J})} \!
(t^{J \setminus \J}\!)^{\!-\! \e^{J \setminus \J} \!-\! \theta 
(\alpha^{J \setminus \J} -\lambda^{J \setminus \J})}\\
\times(\Omega^{(l -\lambda) \chi_{\J}}(\D^\lambda F, 
t^{\J}\!)_{L_p(\R^d)})^\theta dt^{\J} dt^{J \setminus \J}\! =\\
c_6^\theta ( \int_{ ([1, \infty)^d)^{J \setminus \J}} \!
(t^{J \setminus \J}\!)^{\!-\! \e^{J \setminus \J} \!-\! \theta 
(\alpha^{J \setminus \J} -\lambda^{J \setminus \J})} \!dt^{J \setminus \J}\!) \!\int_{ (I^d)^{\J} }
\!(t^{\J}\!)^{\!-\! \e^{\J} \!-\! \theta (\alpha^{\J} -\lambda^{\J})}\\
\times(\Omega^{(l -\lambda) \chi_{\J}}(\D^\lambda F, t^{\J}\!)_{L_p(\R^d)})^\theta 
\!dt^{\J} \!\le\\
c_7 \int_{ (I^d)^{\J} } (t^{\J})^{-\e^{\J} -\theta (\alpha^{\J} -\lambda^{\J})} 
(\Omega^{(l -\lambda) \chi_{\J}}(\D^\lambda F, 
t^{\J})_{L_p(\R^d)})^\theta dt^{\J}.
\end{multline*}

Подставляя последнюю оценку в (2.3.8), получаем, что
\begin{multline*} \tag{2.3.9}
\int_{(\R_+^d)^J} (t^J)^{-\e^J -\theta (\alpha^J -\lambda^J)}
(\Omega^{(l -\lambda) \chi_J}(\D^\lambda F, t^J)_{L_p(\R^d)})^\theta dt^J \le\\
 c_7 \sum_{ \J \subset J} \int_{ (I^d)^{\J} } (t^{\J})^{-\e^{\J}
-\theta (\alpha^{\J} -\lambda^{\J})} 
(\Omega^{(l -\lambda) \chi_{\J}}(\D^\lambda F, 
t^{\J})_{L_p(\R^d)})^\theta dt^{\J} = \\
 c_7 ( \| \D^\lambda F \|_{L_p(\R^d)}^\theta +
 \sum_{ \J \subset J: \J \ne \emptyset} \int_{ (I^d)^{\J} }
(t^{\J})^{-\e^{\J} -\theta (\alpha^{\J} -\lambda^{\J})} \\
\times(\Omega^{(l -\lambda) \chi_{\J}}(\D^\lambda F,
t^{\J})_{L_p(\R^d)})^\theta dt^{\J}).
\end{multline*}

Таким образом, приходим к необходимости оценки $ \int_{ (I^d)^{\J} }
(t^{\J})^{-\e^{\J} -\theta (\alpha^{\J} -\lambda^{\J})}$ 
$(\Omega^{(l -\lambda) \chi_{\J}}(\D^\lambda F,
t^{\J})_{L_p(\R^d)})^\theta dt^{\J} $ для $ \J \subset \{1, \ldots,d\}: 
\J \ne \emptyset. $  Проведём эту оценку.

Пользуясь тем, что для любого непустого множества $ J \subset \{1,\ldots, d\} $ 
имеет место представление $ (I^d)^J = (\cup_{k^J \in (\N^d)^J} 
(2^{-k^J} +2^{-k^J} (I^d)^J)) \cup M^J, \ $
причём множества в правой части последнего равенства попарно не
пересекаются и $ \mes M^J =0, $  получаем, что для любого непустого
множества $ J \subset \{1,\ldots,d\}$ имеет место неравенство
\begin{multline*} \tag{2.3.10}
\int_{ (I^d)^J } (t^J)^{-\e^J -\theta (\alpha^J -\lambda^J)}
(\Omega^{(l -\lambda) \chi_J}(\D^\lambda F, t^J)_{L_p(\R^d)})^\theta dt^J  =\\
\sum_{ k^J \in (\N^d)^J} \!\int_{(2^{-k^J} \!+2^{-k^J} (I^d)^J)} \!
(t^J)^{-\e^J -\theta (\alpha^J -\lambda^J)}
(\Omega^{(l -\lambda) \chi_J}(\D^\lambda F, t^J)_{L_p(\R^d)})^\theta dt^J \le\\
\sum_{ k^J \in (\N^d)^J} \!\int_{(2^{-k^J} \!+2^{-k^J} (I^d)^J)} \!
(2^{-k^J}\!)^{-\e^J -\theta (\alpha^J -\lambda^J)}
(\Omega^{(l -\lambda) \chi_J}(\D^\lambda F, 
2 \!\cdot\! 2^{-k^J}\!)_{L_p(\R^d)})^\theta dt^J\! =\\
\sum_{ k^J \in (\N^d)^J} 2^{\theta (k^J, \alpha^J -\lambda^J)}
(\Omega^{(l -\lambda) \chi_J}(\D^\lambda F, 2 \cdot 2^{-k^J})_{L_p(\R^d)})^\theta.
\end{multline*}

Оценивая $ \Omega^{(l -\lambda) \chi_J} (\D^\lambda F, 2 \cdot 2^{-k^J})_{L_p(\R^d)} $
при $ k^J \in (\N^d)^J, J \subset \{1,\ldots,d\}: J \ne \emptyset, $ 
для $ h \in \R^d: h^J \in 2 \cdot 2^{-k^J} (B^d)^J $ с учётом (2.3.6),
(1.1.7) и (2.2.13) имеем
\begin{multline*}
\| \Delta_h^{(l -\lambda) \chi_J} \D^\lambda F \|_{L_p(\R^d)} = 
\| \Delta_h^{(l -\lambda) \chi_J} ( \sum_{ \kappa \in \Z_+^d} 
\D^\lambda (\mathcal E_\kappa^{d,l -\e,m,\R}(\mathcal I f))) \|_{L_p(\R^d)} = \\
\| \sum_{ \kappa \in \Z_+^d} \Delta_h^{(l -\lambda) \chi_J} 
(\D^\lambda (\mathcal E_\kappa^{d,l -\e,m,\R}(\mathcal I f))) \|_{L_p(\R^d)} \le \\
\sum_{ \kappa \in \Z_+^d} \| \Delta_h^{(l -\lambda) \chi_J} 
(\D^\lambda (\mathcal E_\kappa^{d,l -\e,m,\R}(\mathcal I f))) \|_{L_p(\R^d)} = \\
\sum_{ \J \subset J} \sum_{\substack{ \kappa \in \Z_+^d: 
\kappa^{\J} \le k^{\J},\\ \kappa^{J \setminus \J} > k^{J \setminus \J}}}
\| \Delta_h^{(l -\lambda) \chi_J} 
(\D^\lambda (\mathcal E_\kappa^{d,l -\e,m,\R}(\mathcal I f))) \|_{L_p(\R^d)} \le \\
\sum_{ \J \subset J} \sum_{\substack{ \kappa \in \Z_+^d: 
\kappa^{\J} \le k^{\J},\\ \kappa^{J \setminus \J} > k^{J \setminus \J}}}
c_6 \| \Delta_h^{(l -\lambda) \chi_{\J}} 
(\D^\lambda (\mathcal E_\kappa^{d,l -\e,m,\R}(\mathcal I f))) \|_{L_p(\R^d)} \le \\
\sum_{ \J \subset J} \sum_{\substack{ \kappa \in \Z_+^d: 
\kappa^{\J} \le k^{\J},\\ \kappa^{J \setminus \J} > k^{J \setminus \J}} }
c_6 (\prod_{ j \in \J} | h_j |^{ l_j -\lambda_j}) \| \D^{(l -\lambda) \chi_{\J}} 
(\D^\lambda (\mathcal E_\kappa^{d,l -\e,m,\R}(\mathcal I f))) \|_{L_p(\R^d)} \le \\
\sum_{ \J \subset J} \sum_{\substack{ \kappa \in \Z_+^d: 
\kappa^{\J} \le k^{\J},\\ \kappa^{J \setminus \J} > k^{J \setminus \J}}} 
c_6 (\prod_{ j \in \J} (2^{-k_j +1})^{ l_j -\lambda_j})
\| \D^{(l -\lambda) \chi_{\J} +\lambda} 
(\mathcal E_\kappa^{d,l -\e,m,\R}(\mathcal I f)) \|_{L_p(\R^d)} \le \\
 \sum_{ \J \subset J} \sum_{\substack{ \kappa \in \Z_+^d: 
\kappa^{\J} \le k^{\J},\\ \kappa^{J \setminus \J} > k^{J \setminus \J}} }
c_8 (\prod_{ j \in \J} 2^{-k_j (l_j -\lambda_j)}) 
2^{ ( \kappa, (l -\lambda) \chi_{\J} +\lambda)}
\Omega^{\prime l \chi_{\s(\kappa)}}(f, c_3 (2^{-\kappa})^{\s(\kappa)})_{L_p(I^d)} =\\
c_8 \sum_{ \J \subset J} 2^{ -( k^{\J}, l^{\J} -\lambda^{\J})} 
\sum_{\substack{ \kappa \in \Z_+^d: \kappa^{\J} \le k^{\J},\\ 
\kappa^{J \setminus \J} > k^{J \setminus \J}} }
2^{ ( \kappa^{\J}, l^{\J} -\lambda^{\J})} 2^{(\kappa, \lambda)} 
\Omega^{\prime l \chi_{\s(\kappa)}}(f, c_3 (2^{-\kappa})^{\s(\kappa)})_{L_p(I^d)}.
\end{multline*}

Отсюда получаем, что для $ k^J \in (\N^d)^J, J \subset
\{1,\ldots,d\}: J \ne \emptyset, $ выполняется неравенство
\begin{multline*}
\Omega^{(l -\lambda) \chi_J} (\D^\lambda F, 2 \cdot 2^{-k^J})_{L_p(\R^d)} \le\\
 c_8 \sum_{ \J \subset J} 2^{ -( k^{\J}, l^{\J} -\lambda^{\J})} 
\sum_{\substack{ \kappa \in \Z_+^d: \kappa^{\J} \le k^{\J},\\ 
\kappa^{J \setminus \J} > k^{J \setminus \J}}}
2^{ ( \kappa^{\J}, l^{\J} -\lambda^{\J})} 2^{(\kappa, \lambda)} 
\Omega^{\prime l \chi_{\s(\kappa)}}(f, c_3 (2^{-\kappa})^{\s(\kappa)})_{L_p(I^d)}.
\end{multline*}

Подставляя эту оценку  в (2.3.10) и применяя неравенство Гёльдера,
для $  J \subset \{1,\ldots,d\}: J \ne \emptyset, $ выводим
\begin{multline*} \tag{2.3.11}
\int_{ (I^d)^J } (t^J)^{-\e^J -\theta (\alpha^J -\lambda^J)} 
(\Omega^{(l -\lambda) \chi_J}(\D^\lambda F, t^J)_{L_p(\R^d)})^\theta dt^J \le 
\sum_{ k^J \in (\N^d)^J} 2^{\theta (k^J, \alpha^J -\lambda^J)} \\
\times\biggl( c_8 \sum_{ \J \subset J} 2^{ -( k^{\J}, l^{\J} -\lambda^{\J})} 
\sum_{\substack{ \kappa \in \Z_+^d: \kappa^{\J} \le k^{\J}, \\
\kappa^{J \setminus \J} > k^{J \setminus \J}} }
2^{ ( \kappa^{\J}, l^{\J} -\lambda^{\J})} 2^{(\kappa, \lambda)}
\Omega^{\prime l \chi_{\s(\kappa)}}(f, c_3 (2^{-\kappa})^{\s(\kappa)})_{L_p(I^d)} \biggr)^\theta \le \\
 c_{9} \sum_{ k^J \in (\N^d)^J} 2^{\theta (k^J, \alpha^J -\lambda^J)}
\sum_{ \J \subset J} \biggl(2^{ -( k^{\J}, l^{\J} -\lambda^{\J})} 
\sum_{ \kappa \in \Z_+^d: \kappa^{\J} \le k^{\J}, 
\kappa^{J \setminus \J} > k^{J \setminus \J}} 
2^{ ( \kappa^{\J}, l^{\J} -\lambda^{\J})} \\
\times2^{(\kappa, \lambda)}
\Omega^{\prime l \chi_{\s(\kappa)}}(f, c_3 (2^{-\kappa})^{\s(\kappa)})_{L_p(I^d)} \biggr)^\theta = \\
c_{9} \sum_{ k^J \in (\N^d)^J} \sum_{ \J \subset J} 
2^{\theta (k^J, \alpha^J -\lambda^J) -\theta ( k^{\J}, l^{\J} -\lambda^{\J})} \\
\times\biggl(\sum_{ \kappa \in \Z_+^d: \kappa^{\J} \le k^{\J}, 
\kappa^{J \setminus \J} > k^{J \setminus \J}}
2^{ ( \kappa^{\J}, l^{\J} -\lambda^{\J})} 2^{(\kappa, \lambda)}
\Omega^{\prime l \chi_{\s(\kappa)}}(f, c_3 (2^{-\kappa})^{\s(\kappa)})_{L_p(I^d)} \biggr)^\theta =\\
c_{9} \sum_{ \J \subset J} \sum_{ k^J \in (\N^d)^J} 
2^{\theta (k^J, \alpha^J -\lambda^J) -\theta ( k^{\J}, l^{\J} -\lambda^{\J})} \\
\times\biggl(\sum_{ \kappa \in \Z_+^d: \kappa^{\J} \le k^{\J}, 
\kappa^{J \setminus \J} > k^{J \setminus \J}}
2^{ ( \kappa^{\J}, l^{\J} -\lambda^{\J})} 2^{(\kappa, \lambda)} 
\Omega^{\prime l \chi_{\s(\kappa)}}(f, c_3 (2^{-\kappa})^{\s(\kappa)})_{L_p(I^d)} \biggr)^\theta.
\end{multline*}

Для оценки правой части (2.3.11) при $ \lambda \in \Z_+^d: \lambda < \alpha, $ 
фиксируем $ \epsilon \in \R_+^d, $ для которого соблюдаются условия 
$ \alpha -\lambda -\epsilon >0, l -\alpha -\epsilon >0, $ и положим 
$ \overline J = \{1, \ldots, d\} \setminus J. $ Используя неравенство 
Гёльдера, для $  J \subset \{1,\ldots,d\}: J \ne \emptyset, \J \subset J, 
k^J \in (\N^d)^J $ приходим к неравенству
\begin{multline*} \tag{2.3.12}
\biggl(\sum_{ \kappa \in \Z_+^d: \kappa^{\J} \le k^{\J}, 
\kappa^{J \setminus \J} > k^{J \setminus \J}} 
2^{ ( \kappa^{\J}, l^{\J} -\lambda^{\J})} 2^{(\kappa, \lambda)} 
\Omega^{\prime l \chi_{\s(\kappa)}}(f,
c_3 (2^{-\kappa})^{\s(\kappa)})_{L_p(I^d)}\biggr)^\theta =\\
\biggl(\sum_{ \kappa \in \Z_+^d: \kappa^{\J} \le k^{\J},
\kappa^{J \setminus \J} > k^{J \setminus \J}}
2^{-\!(\kappa, \alpha)} 2^{ (\kappa^{\J} \!, l^{\J} -\lambda^{\J})} 
2^{(\kappa, \lambda)} 2^{(\kappa^{J \setminus \J} \!, \epsilon^{J \setminus \J})} \\
\times 2^{-\!(\kappa^{\J} \!, \epsilon^{\J})} 2^{-\!(\kappa^{J \setminus \J} \!, 
\epsilon^{J \setminus \J})} 2^{ (\kappa^{\J} \!, \epsilon^{\J})} 
2^{(\kappa, \alpha)} \Omega^{\prime l \chi_{\s(\kappa)}}(f, 
c_3 (2^{-\kappa})^{\s(\kappa)})_{L_p(I^d)} \biggr)^\theta \le\\
\biggl(\sum_{ \kappa \in \Z_+^d: \kappa^{\J} \le k^{\J}, 
\kappa^{J \setminus \J} > k^{J \setminus \J}}\! 2^{-\theta^\prime (\kappa,
\alpha -\lambda)} 2^{\theta^\prime (\kappa^{J \setminus \J}, 
\epsilon^{J \setminus \J})}
2^{ \theta^\prime (\kappa^{\J}, l^{\J} -\lambda^{\J} -\epsilon^{\J})}\biggr)^{\theta / \theta^\prime} \\
\times\biggl(\sum_{ \kappa \in \Z_+^d: \kappa^{\J} \le k^{\J}, 
\kappa^{J \setminus \J} > k^{J \setminus \J}} (2^{-(\kappa^{J \setminus \J},
\epsilon^{J \setminus \J})} 2^{ (\kappa^{\J}, \epsilon^{\J})}
2^{(\kappa, \alpha)} \Omega^{\prime l \chi_{\s(\kappa)}}(f, 
c_3 (2^{-\kappa})^{\s(\kappa)})_{L_p(I^d)} )^\theta\biggr).
\end{multline*}

Оценивая правую часть (2.3.12), для $ J \subset \{1,\ldots,d\}: J \ne \emptyset, 
\J \subset J, k^J \in (\N^d)^J $ имеем
\begin{multline*} \tag{2.3.13}
\biggl(\sum_{ \kappa \in \Z_+^d: \kappa^{\J} \le k^{\J}, 
\kappa^{J \setminus \J} > k^{J \setminus \J}} 2^{-\theta^\prime (\kappa, \alpha -\lambda)} 
2^{\theta^\prime (\kappa^{J \setminus \J}, \epsilon^{J \setminus \J})} 
2^{ \theta^\prime (\kappa^{\J}, l^{\J} -\lambda^{\J} -\epsilon^{\J})}\biggr)^{\theta / \theta^\prime} = \\
\biggl(\sum_{ \kappa^{\overline J} \in (\Z_+^d)^{\overline J},
\kappa^{\J} \in (\Z_+^d)^{\J}: \kappa^{\J} \le k^{\J}, 
\kappa^{J \setminus \J} \in (\Z_+^d)^{J \setminus \J}: 
\kappa^{J \setminus \J} > k^{J \setminus \J}}
2^{-\theta^\prime (\kappa^{\overline J}, \alpha^{\overline J} -\lambda^{\overline J})}  \\
\times 2^{-\theta^\prime (\kappa^{J \setminus \J}, \alpha^{J \setminus \J} -
\lambda^{J \setminus \J} -\epsilon^{J \setminus \J})} 
2^{ \theta^\prime (\kappa^{\J}, l^{\J} -\alpha^{\J} -\epsilon^{\J})} \biggr)^{\theta /
\theta^\prime} = 
\biggl(\biggl(\sum_{ \kappa^{\overline J} \in
(\Z_+^d)^{\overline J}} 2^{-\theta^\prime (\kappa^{\overline J},
 \alpha^{\overline J} -\lambda^{\overline J})} \biggr)  \\
\times\biggl(\sum_{\kappa^{J \setminus \J} \in (\Z_+^d)^{J \setminus \J}: 
\kappa^{J \setminus \J} > k^{J \setminus \J}}
2^{-\theta^\prime (\kappa^{J \setminus \J}, \alpha^{J \setminus \J} 
-\lambda^{J \setminus \J} -\epsilon^{J \setminus \J})} \biggr)  \\
\times\biggl(\sum_{\kappa^{\J} \in (\Z_+^d)^{\J}: \kappa^{\J} \le k^{\J}} 
2^{ \theta^\prime (\kappa^{\J}, l^{\J} -\alpha^{\J} -\epsilon^{\J})} \biggr)\biggr)^{\theta / \theta^\prime} = 
\biggl((\prod_{j \in \overline J}
(\sum_{\kappa_j =0}^\infty
2^{-\theta^\prime \kappa_j (\alpha_j -\lambda_j)}))  \\
\times (\prod_{j \in J \setminus \J} (\sum_{\kappa_j = k_j +1}^\infty 
2^{-\theta^\prime \kappa_j (\alpha_j -\lambda_j -\epsilon_j)})) \cdot
(\prod_{j \in \J} (\sum_{\kappa_j =0}^{ k_j} 
2^{ \theta^\prime \kappa_j (l_j -\alpha_j -\epsilon_j)}))\biggr)^{\theta / \theta^\prime} \le \\
 \biggl( c_{10} (\prod_{j \in J \setminus \J} 
2^{-\theta^\prime k_j (\alpha_j -\lambda_j -\epsilon_j)}) 
(\prod_{j \in \J} 2^{ \theta^\prime k_j
(l_j -\alpha_j -\epsilon_j)} )\biggr)^{\theta / \theta^\prime} = \\
 c_{11} \biggl(2^{-\theta^\prime (k^{J \setminus \J}, \alpha^{J \setminus \J} -
\lambda^{J \setminus \J} -\epsilon^{J \setminus \J})} 
2^{ \theta^\prime (k^{\J}, l^{\J} -\alpha^{\J} -\epsilon^{\J})}\biggr)^{\theta / \theta^\prime} = \\
c_{11} 2^{-\theta (k^{J \setminus \J}, \alpha^{J \setminus \J}
-\lambda^{J \setminus \J} -\epsilon^{J \setminus \J})} 
2^{ \theta (k^{\J}, l^{\J} -\alpha^{\J} -\epsilon^{\J})}.
\end{multline*}

Соединяя (2.3.11), (2.3.12), (2.3.13), получаем неравенство
\begin{multline*} \tag{2.3.14}
\int_{ (I^d)^J } (t^J)^{-\e^J -\theta (\alpha^J -\lambda^J)} 
(\Omega^{(l -\lambda) \chi_J}(\D^\lambda F, t^J)_{L_p(\R^d)})^\theta dt^J \le \\
c_{9} \sum_{ \J \subset J} \sum_{ k^J \in (\N^d)^J} 
2^{\theta (k^J, \alpha^J -\lambda^J) -\theta ( k^{\J}, l^{\J} -\lambda^{\J})} \\
\times c_{11} 2^{-\theta (k^{J \setminus \J}, \alpha^{J \setminus \J} 
-\lambda^{J \setminus \J} -\epsilon^{J \setminus \J})}
2^{ \theta (k^{\J}, l^{\J} -\alpha^{\J} -\epsilon^{\J})} \\
\biggl(\sum_{ \kappa \in \Z_+^d: \kappa^{\J} \le k^{\J}, 
\kappa^{J \setminus \J} > k^{J \setminus \J}} (2^{-(\kappa^{J \setminus \J},
\epsilon^{J \setminus \J})} 2^{ (\kappa^{\J}, \epsilon^{\J})}
2^{(\kappa, \alpha)}
\Omega^{\prime l \chi_{\s(\kappa)}}(f, c_3 (2^{-\kappa})^{\s(\kappa)})_{L_p(I^d)} )^\theta\biggr) = \\
c_{12} \sum_{ \J \subset J} \sum_{ k^J \in (\N^d)^J} 
(2^{\theta (k^{J \setminus \J}, \epsilon^{J \setminus \J})}
2^{ -\theta (k^{\J}, \epsilon^{\J})}) \\
\sum_{ \kappa \in \Z_+^d: \kappa^{\J} \le k^{\J}, 
\kappa^{J \setminus \J} > k^{J \setminus \J}} \biggl(2^{-(\kappa^{J \setminus \J},
\epsilon^{J \setminus \J})} 2^{ (\kappa^{\J}, \epsilon^{\J})}
2^{(\kappa, \alpha)}
\Omega^{\prime l \chi_{\s(\kappa)}}(f, c_3 (2^{-\kappa})^{\s(\kappa)})_{L_p(I^d)} \biggr)^\theta = \\
 c_{12} \sum_{ \J \subset J} \sum_{ k^J \in (\N^d)^J} 
\sum_{\kappa \in \Z_+^d: \kappa^{\J} \le k^{\J}, 
\kappa^{J \setminus \J} > k^{J \setminus \J}} 
(2^{\theta (k^{J \setminus \J}, \epsilon^{J \setminus \J})}
2^{ -\theta (k^{\J}, \epsilon^{\J})})  \\
\times (2^{-(\kappa^{J \setminus \J}, \epsilon^{J \setminus \J})}
2^{ (\kappa^{\J}, \epsilon^{\J})} 2^{(\kappa, \alpha)}
\Omega^{\prime l \chi_{\s(\kappa)}}(f, c_3 (2^{-\kappa})^{\s(\kappa)})_{L_p(I^d)} )^\theta = \\
 c_{12} \sum_{ \J \subset J} \sum_{ k^J \in (\N^d)^J, \kappa \in \Z_+^d: 
\kappa^{\J} \le k^{\J}, \kappa^{J \setminus \J} > k^{J \setminus \J}} 
(2^{\theta (k^{J \setminus \J}, \epsilon^{J \setminus \J})}
2^{ -\theta (k^{\J}, \epsilon^{\J})})  \\
\times (2^{-(\kappa^{J \setminus \J}, \epsilon^{J \setminus \J})}
2^{ (\kappa^{\J}, \epsilon^{\J})} 2^{(\kappa, \alpha)}
\Omega^{\prime l \chi_{\s(\kappa)}}(f, c_3 (2^{-\kappa})^{\s(\kappa)})_{L_p(I^d)} )^\theta = \\
c_{12} \sum_{ \J \subset J} \sum_{\kappa \in \Z_+^d} \sum_{ k^J \in (\N^d)^J: 
\kappa^{\J} \le k^{\J}, \kappa^{J \setminus \J} > k^{J \setminus \J}} 
(2^{\theta (k^{J \setminus \J}, \epsilon^{J \setminus \J})}
2^{ -\theta (k^{\J}, \epsilon^{\J})})  \\
\times (2^{-(\kappa^{J \setminus \J}, \epsilon^{J \setminus \J})}
2^{ (\kappa^{\J}, \epsilon^{\J})} 2^{(\kappa, \alpha)}
\Omega^{\prime l \chi_{\s(\kappa)}}(f, c_3 (2^{-\kappa})^{\s(\kappa)})_{L_p(I^d)} )^\theta = \\
c_{12} \sum_{ \J \subset J} \sum_{\kappa \in \Z_+^d}
(2^{-(\kappa^{J \setminus \J}, \epsilon^{J \setminus \J})} 
2^{(\kappa^{\J}, \epsilon^{\J})} 2^{(\kappa, \alpha)}
\Omega^{\prime l \chi_{\s(\kappa)}}(f, c_3 (2^{-\kappa})^{\s(\kappa)})_{L_p(I^d)} )^\theta \\
\sum_{ k^J \in (\N^d)^J: \kappa^{\J} \le k^{\J}, 
\kappa^{J \setminus \J} > k^{J \setminus \J}} 
(2^{\theta (k^{J \setminus \J}, \epsilon^{J \setminus \J})}
2^{ -\theta (k^{\J}, \epsilon^{\J})}), \\
 J \subset \{1,\ldots,d\}: J \ne \emptyset.
\end{multline*}

Оценивая правую часть (2.3.14), для $  J \subset \{1,\ldots,d\}: J
\ne \emptyset, \J \subset J, \kappa \in \Z_+^d $ имеем
\begin{multline*}
\sum_{ k^J \in (\N^d)^J: \kappa^{\J} \le k^{\J}, 
\kappa^{J \setminus \J} > k^{J \setminus \J}} 
(2^{\theta (k^{J \setminus \J}, \epsilon^{J \setminus \J})}
2^{ -\theta (k^{\J}, \epsilon^{\J})}) = \\
 \sum_{ k^{\J} \in (\N^d)^{\J}: \kappa^{\J} \le k^{\J}, 
k^{J \setminus \J} \in (\N^d)^{J \setminus \J}: 
\kappa^{J \setminus \J} > k^{J \setminus \J}} 
2^{\theta (k^{J \setminus \J}, \epsilon^{J \setminus \J})}
2^{ -\theta (k^{\J}, \epsilon^{\J})} = \\
\biggl(\sum_{ k^{\J} \in (\N^d)^{\J}: \kappa^{\J} \le k^{\J}} 
2^{ -\theta (k^{\J}, \epsilon^{\J})}\biggr) 
\biggl(\sum_{k^{J \setminus \J} \in (\N^d)^{J \setminus \J}: 
\kappa^{J \setminus \J} > k^{J \setminus \J}}
2^{\theta (k^{J \setminus \J}, \epsilon^{J \setminus \J})}\biggr) = \\
\biggl( \prod_{j \in \J} (\sum_{ k_j = \kappa_j}^\infty 
2^{ -\theta k_j \epsilon_j})\biggr) \biggl(\prod_{j \in J \setminus \J} 
(\sum_{k_j =1}^{\kappa_j -1} 2^{\theta k_j \epsilon_j})\biggr) \le \\
c_{13} (\prod_{j \in \J} 2^{ -\theta \kappa_j \epsilon_j}) 
(\prod_{j \in J \setminus \J}
2^{\theta \kappa_j \epsilon_j}) = 
c_{13} 2^{ -\theta (\kappa^{\J}, \epsilon^{\J})} 
2^{\theta(\kappa^{J \setminus \J}, \epsilon^{J \setminus \J})}.
\end{multline*}

Подставляя эту оценку в (2.3.14) и учитывая (2.3.3), получаем, что для 
$ J \subset \{1,\ldots,d\}: J \ne \emptyset, $ справедливо неравенство
\begin{multline*} \tag{2.3.15}
\int_{ (I^d)^J } (t^J)^{-\e^J -\theta (\alpha^J -\lambda^J)}
(\Omega^{(l -\lambda) \chi_J}(\D^\lambda F, t^J)_{L_p(\R^d)})^\theta dt^J \le \\
c_{12} \!\sum_{ \J \subset J} \sum_{\kappa \in \Z_+^d}
\!(2^{-(\kappa^{J \setminus \J}, \epsilon^{J \setminus \J})} 
2^{(\kappa^{\J}, \epsilon^{\J})} 2^{(\kappa, \alpha)} 
\Omega^{\prime l \chi_{\s(\kappa)}}(f, 
c_3 (2^{-\kappa})^{\s(\kappa)})_{L_p(I^d)})^\theta \\
\times c_{13} 2^{ -\theta (\kappa^{\J}, \epsilon^{\J})}
2^{\theta (\kappa^{J \setminus \J}, \epsilon^{J \setminus \J})} = \\
c_{14} \!\sum_{ \J \subset J} \sum_{\kappa \in \Z_+^d}
\!(2^{(\kappa, \alpha)} \Omega^{\prime l \chi_{\s(\kappa)}}(f, 
c_3 (2^{-\kappa})^{\s(\kappa)})_{L_p(I^d)} )^\theta \le\\
c_{15} \sum_{\kappa \in \Z_+^d} (2^{(\kappa, \alpha)} 
\Omega^{\prime l \chi_{\s(\kappa)}}(f, c_3 (2^{-\kappa})^{\s(\kappa)})_{L_p(I^d)})^\theta = \\
c_{15} (\| f \|_{L_p(I^d)}^\theta +\sum_{\kappa \in \Z_+^d \setminus \{0\}} 
(2^{(\kappa, \alpha)} \Omega^{\prime l \chi_{\s(\kappa)}}(f, 
c_3 (2^{-\kappa})^{\s(\kappa)})_{L_p(I^d)})^\theta) \le \\
c_{16} \biggl(\| f\|_{L_p(I^d)}^\theta +\sum_{ \J \subset \Nu_{1,d}^1: \J \ne \emptyset} 
\int_{(\R_+^d)^{\J}} (t^{\J})^{-\e^{\J} -\theta \alpha^{\J}} (\Omega^{\prime l \chi_{\J}}(f,
t^{\J})_{L_p(I^d)})^\theta dt^{\J}\biggr) \le\\
c_{17} \| f \|_{(S_{p, \theta}^\alpha B)^\prime(I^d)}^\theta. 
\end{multline*}

Соединяя (2.3.7), (2.3.9), (2.3.15), приходим к (2.3.2), чем завершаем
доказательство теоремы при $ \theta \ne \infty. $ При $ \theta = \infty $ 
доказательство теоремы проводится по той же схеме с заменой в соответствующих 
местах операции суммирования на операцию взятия супремума или максимума.
$ \square $

Следствие

В условиях теоремы 2.3.1 при $ \l \in \Z_+^d: \l < \alpha, $ имеет место 
включение
\begin{equation*} \tag{2.3.16}
(S_{p, \theta}^\alpha B)^\prime(I^d) \subset (S_{p, \theta}^\alpha B)^{\l}(I^d), \ 
\end{equation*}
и для любой функции $ f \in (S_{p, \theta}^\alpha B)^\prime(I^d) $ 
выполняется неравенство
\begin{equation*} \tag{2.3.17} 
\| f\|_{(S_{p, \theta}^\alpha B)^{\l}(I^d)} \le C
\| f\|_{(S_{p, \theta}^\alpha B)^\prime(I^d)}.
\end{equation*}

Для получения (2.3.16), (2.3.17) достаточно применить (2.3.1), оценку
\begin{equation*}
\| (\mathcal E^{d,\alpha,p,\theta} f) \mid_{I^d} \|_{(S_{p, \theta}^\alpha B)^{\l}(I^d)} \le 
\| \mathcal E^{d,\alpha,p,\theta} f\|_{(S_{p, \theta}^\alpha B)^{\l}(\R^d)}
\end{equation*}
и неравенство (2.3.2) при 
$ \lambda = \l \chi_J, J \subset \{1,\ldots,d\}. $
Из (1.1.9), (1.1.10) и (2.3.16), (2.3.17)
вытекает, что $ (S_{p, \theta}^\alpha B)^\prime(I^d) = 
(S_{p, \theta}^\alpha B)^{\l}(I^d), $ и нормы на этих пространствах 
$ \| \cdot \|_{(S_{p, \theta}^\alpha B)^\prime(I^d)},
\| \cdot \|_{(S_{p, \theta}^\alpha B)^{\l}(I^d)} $ эквивалентны 
при $ \l \in \Z_+^d: \l < \alpha. $

Отметим также, что из теоремы 2.3.1 и (1.1.6) следует, что при 
$ \lambda \in \Z_+^d: \lambda < \alpha, $ имеет место включение
$ \D^\lambda \mid_{(S_{p, \theta}^\alpha B)^\prime(I^d)} \in
\mathcal B((S_{p, \theta}^\alpha B)^\prime(I^d),
(S_{p, \theta}^{\alpha -\lambda} B)^\prime(I^d)). $ 

\end{document}